\documentclass[11pt,a4paper,oneside]{article}

\usepackage[english]{babel}
\usepackage[T1]{fontenc}
\usepackage[utf8]{inputenc}

\usepackage{lmodern}
\usepackage{indentfirst}
\usepackage{microtype}

\usepackage{graphicx}
\usepackage{caption}
\usepackage{subcaption}
\usepackage{import}

\usepackage{amsmath}
\usepackage{amsthm}
\usepackage{amsfonts}
\usepackage{amssymb}
\usepackage{amsbsy}
\usepackage{amstext}
\usepackage{amscd}
\usepackage{xfrac}
\usepackage{mathtools}
\usepackage{thmtools}
\usepackage{thm-restate}

\usepackage{cite}
\usepackage{hyperref}
\usepackage[capitalise, noabbrev]{cleveref}

\usepackage{authblk}
\usepackage{lipsum}
\usepackage{color}
\usepackage{enumitem}
\usepackage[normalem]{ulem}
\usepackage{hyphenat}

\hypersetup{
  pdftitle={Dynamics and Stability of Non-Smooth Dynamical Systems with Two Switches},
  pdfauthor={Guilherme Tavares da Silva and Ricardo Miranda Martins},
  pdfsubject={Academic paper},
  pdfcreator={LaTeX},
  pdfkeywords={Dynamical systems. Filippov systems. Structural stability.
      Singular perturbations (Mathematics).},
  colorlinks=true,
  linkcolor=blue,
  citecolor=blue,
  urlcolor=blue,
  filecolor=magenta}

\providecommand{\keywords}[1]
{
  \small
  \textbf{Keywords:} #1
}

\crefformat{equation}{(#2#1#3)}
\crefmultiformat{equation}%
{(#2#1#3)}{ and~(#2#1#3)}{, (#2#1#3)}{ and~(#2#1#3)}

\crefformat{condition}{#2#1#3}
\crefmultiformat{condition}%
{#2#1#3}{ and~#2#1#3}{, #2#1#3}{ and~#2#1#3}

\theoremstyle{plain}
\newtheorem{theorem}{Theorem}[section]

\newtheorem{proposition}{Proposition}[section]
\newtheorem{corollary}{Corollary}[section]

\theoremstyle{definition}
\newtheorem{definition}{Definition}[section]
\newtheorem{example}{Example}[section]

\declaretheorem[name=Theorem,numberwithin=section]{rtheorem}

\DeclareMathOperator{\sgn}{\mathrm{sgn}}

\newcommand{\norm}[1]{\left\lVert #1 \right\rVert}

\newcommand\restr[2]{{
\left.\kern-\nulldelimiterspace
#1
\vphantom{\big|}
\right|_{#2}}
}

\title{\bf Dynamics and Stability of Non-Smooth Dynamical Systems with Two Switches}

\author{Guilherme Tavares da Silva\thanks{g119509@dac.unicamp.br} \quad}
\author{\quad Ricardo Miranda Martins\thanks{Corresponding author: rmiranda@unicamp.br}}

\affil{\footnotesize Departamento de Matem\'atica,
  Universidade Estadual de Campinas,
  Rua S\'ergio Buarque de Holanda, 651,
  Cidade Universitária Zeferino Vaz,
  13083--859, Campinas, SP, Brazil}

\date{}

\begin{document}

\frenchspacing

\maketitle

\begin{abstract}
  One of the most common hypotheses on the theory of non-smooth dynamical
  systems is a regular surface as switching manifold, at which case there is at
  least well-defined and established Filippov dynamics. However, systems with
  singular switching manifolds still lack such well-established dynamics,
  although present in many relevant models of phenomena where multiple switches
  or multiple abrupt changes occur. At this work, we leverage a methodology
  that, through blow-ups and singular perturbation, allows the extension of
  Filippov dynamics to the singular case. Specifically, tridimensional systems
  whose switching manifold consists of an algebraic manifold with transversal
  self-intersection are considered. This configuration, known as double
  discontinuity, represents systems with two switches and whose singular part
  consists of a straight line, where ordinary Filippov dynamics is not directly
  appliable. For the general, non-linear case, beyond defining the so-called
  fundamental dynamics over the singular part, general theorems on its
  qualitative behavior are provided. For the affine case, however, theorems
  fully describing the fundamental dynamics are obtained. Finally, this
  fine-grained control over the dynamics is leveraged to derive Peixoto-like
  theorems characterizing semi-local structural stability.
\end{abstract}


\clearpage

\section{Introduction}%
\label{sec:introduction}

The theory of dynamical systems given by ordinary differential equations

\begin{equation}
  \label{sec:problem:eq:smooth_system}
  \dot{\mathbf{x}} = \mathbf{F}(\mathbf{x}),
\end{equation}

\noindent where $\mathbf{F}: \mathbb{R}^n \to \mathbb{R}^n$ is at least a
continuous vector field evolved naturally with the birth of Calculus itself,
with \cite{Arrowsmith1990} and \cite{Katok1997} being exceptional modern
references on the subject. In fact, the machinery provided by this theory has
been used in the study of models all around science: from classical Newtonian
Mechanics to modern Machine Learning \cite{Weinan2017}.

However, either naturally or due to simplifications and practicality, many of
these phenomena are better approached with non-smooth models, i.e., where the
vector field $\mathbf{F}$ above has discontinuities. More specifically, given $U
  \subset \mathbb{R}^n$ open, $h: U \to \mathbb{R}$ continually differentiable
having $0$ as a \textbf{regular value} and two vector fields $\mathbf{F}_{\pm}:
  U \to \mathbb{R}^n$ of class $C^k(U)$ with $k \geq 1$, we understand as a
\textbf{non-smooth dynamical system} that given by a differential equation as
\cref{sec:problem:eq:smooth_system} where

\begin{equation}
  \label{sec:problem:eq:nonsmooth_system}
  \mathbf{F}(\mathbf{x}) =
  \begin{cases}
    \mathbf{F}_+(\mathbf{x}), & \text{ if } \mathbf{x} \in \Sigma_{+}, \\
    \mathbf{F}_-(\mathbf{x}), & \text{ if } \mathbf{x} \in \Sigma_{-},
  \end{cases}
\end{equation}

\noindent with $\Sigma_{+} = \left\{\mathbf{x} \in U ;~ h(\mathbf{x}) \ge
  0\right\}$ and $\Sigma_{-} = \left\{\mathbf{x} \in U ;~ h(\mathbf{x}) \le
  0\right\}$ intersecting at a regular surface $\Sigma$ called \textbf{switching
  manifold}. We denote the set of vector fields $\mathbf{F}$ defined as above by

\begin{equation*}
  \mathcal{R}^k(U) \equiv C^k(U,\mathbb{R}^n) \times C^k(U,\mathbb{R}^n)
\end{equation*}

\noindent which we consider equipped with the Whitney topology. Generally, we
write $\mathbf{F} = (\mathbf{F_+}, \mathbf{F_-})$ to denote the elements of this
set. These systems arise frequently, for instance, in the study of mechanical
systems with impact or friction \cite{Brogliato1999, Hinrichs1998,
  Kowalczyk2008, Wojewoda2008}, electronic circuits with switches
\cite{Bernardo1998, Bernardo2001, Cristiano2018, Wang2017}, biological and
climate models with abrupt changes \cite{Barry2017, Carvalho2020, Leifeld2018,
  Piltz2014, Prokopiou2014}, economics and politics \cite{Amador2013,
  ValenciaCalvo2020, Wang2018}, etc. Hence, not only due to its applications, but
also, its mathematical beauty, the theory of non-smooth dynamical system is a
very active field, attracting and mobilizing scientists from all around the
world.

In this endeavor, the establishment of definitions is one of the main
challenges. The definition of solution for a non-smooth system, for instance, is
not always clear. Nevertheless, in this context, one of the greatest
contributions came from Filippov in \cite{Filippov1988}, which introduced a
convention to define such solutions in such a way that, apparently, is both
geometrically beautiful and consistent with the physical world\footnote{It is
  important to remark, however, that other conventions exist with equal beauty.
  Filippov's convention just happens to be the most accepted one nowadays. For
  instance, we cite here Carathéodory \cite{Spraker1996} and Utkin
  \cite{Utkin1992} conventions. See also \cite{Glendinning2019} for some
  historical aspects.}. More specifically, for points $\mathbf{x} \in U \setminus
  \Sigma$, the usual local dynamics of the fields $\mathbf{F}_{\pm}$ is
considered. On the other hand, roughly speaking, for points $\mathbf{x} \in
  \Sigma$ and considering the Lie derivative $\mathbf{F}_{\pm}h(\mathbf{x})
  \coloneqq \nabla h(\mathbf{x}) \cdot \mathbf{F}_{\pm}(\mathbf{x})$, the
switching manifold $\Sigma$ splits into three regions:

\begin{itemize}
  \item \textbf{Crossing Region}: $\Sigma^{cr} = \left\{\mathbf{x} \in \Sigma;~
          \mathbf{F}_+h(\mathbf{x}) \mathbf{F}_-h(\mathbf{x}) > 0\right\}$,
        where touching trajectories cross $\Sigma$ through concatenation.

  \item \textbf{Sliding Region}: $\Sigma^{sl} = \left\{\mathbf{x} \in \Sigma;~
          \mathbf{F}_+h(\mathbf{x}) > 0,~ \mathbf{F}_-h(\mathbf{x}) <
          0\right\}$, where touching trajectories remains tangent to $\Sigma$
        for positive time.

  \item \textbf{Escaping Region}: $\Sigma^{es} = \left\{\mathbf{x} \in \Sigma;~
          \mathbf{F}_+h(\mathbf{x}) < 0,~ \mathbf{F}_-h(\mathbf{x}) >
          0\right\}$, where touching trajectories remains tangent to $\Sigma$
        for negative time.
\end{itemize}

Due to the continuity, all regions above are open sets separated by the
so-called \textbf{tangency points} $\mathbf{x} \in \Sigma$ where
$\mathbf{F}_+h(\mathbf{x}) \mathbf{F}_-h(\mathbf{x}) = 0$ which, dynamically,
acts as singularities. Moreover, for points $\mathbf{x} \in \Sigma^s \coloneqq
  \Sigma^{sl} \cup \Sigma^{es}$, the trajectory slides tangent to $\Sigma$
according to a well-defined \textbf{sliding vector field} $\mathbf{F}^s:\Sigma^s
  \to T\Sigma^s$ given by

\begin{equation}
  \label{sec:problem:eq:sliding}
  \mathbf{F}^s(\mathbf{x}) =
  \frac{\mathbf{F}_-h(\mathbf{x}) \mathbf{F}_+(\mathbf{x}) -
    \mathbf{F}_+h(\mathbf{x}) \mathbf{F}_-(\mathbf{x})}
  {\mathbf{F}_-h(\mathbf{x}) - \mathbf{F}_+h(\mathbf{x})},
\end{equation}

\noindent which consists of the single vector in the intersection
\[\text{Conv}(\{\mathbf{F}_+(\mathbf{x}),\mathbf{F}_-(\mathbf{x})\}) \cap
  \Sigma,\] where $\text{Conv}(\cdot)$ represents \textbf{convex hull}.

Using the above construction, many advances have been achieved on this class of
systems concerning, for instance, its bifurcations \cite{Guardia2011},
regularization \cite{Panazzolo2017, Sotomayor1996, Teixeira2012}, structural
stability \cite{Broucke2001, Gomide2020, Teixeira1990} and uncountable works
regarding minimal sets. However, as previously observed, the theory established
by Filippov's convention has a fundamental hypothesis: a regular surface as
switching manifold between the smooth parts of the system, i.e., a surface
$\Sigma = h^{-1}(\left\{0\right\})$ where $0$ is a regular value of $h$. Many
relevant phenomena, however, require a model where $\Sigma$ is, actually, the
preimage of a singular value. Generally speaking, models where two or more
abrupt changes might occur. See, for instance, the ``On or Off Genes'' section
in \cite[p.~28]{Jeffrey2020}, where a model is presented for two genes
interacting in an organic cell of a living system in order to produce proteins.
At that same reference, \cite[p.~30]{Jeffrey2020}, the section ``Jittery
Investments'' presents another interesting model for a game with two players
buying or selling stocks of a company.

\begin{figure}[ht]
  \centering
  \begin{subfigure}[b]{0.4\linewidth}
    \def\svgwidth{\linewidth}
\begingroup%
  \makeatletter%
  \providecommand\color[2][]{%
    \errmessage{(Inkscape) Color is used for the text in Inkscape, but the package 'color.sty' is not loaded}%
    \renewcommand\color[2][]{}%
  }%
  \providecommand\transparent[1]{%
    \errmessage{(Inkscape) Transparency is used (non-zero) for the text in Inkscape, but the package 'transparent.sty' is not loaded}%
    \renewcommand\transparent[1]{}%
  }%
  \providecommand\rotatebox[2]{#2}%
  \newcommand*\fsize{\dimexpr\f@size pt\relax}%
  \newcommand*\lineheight[1]{\fontsize{\fsize}{#1\fsize}\selectfont}%
  \ifx\svgwidth\undefined%
    \setlength{\unitlength}{219.97349536bp}%
    \ifx\svgscale\undefined%
      \relax%
    \else%
      \setlength{\unitlength}{\unitlength * \real{\svgscale}}%
    \fi%
  \else%
    \setlength{\unitlength}{\svgwidth}%
  \fi%
  \global\let\svgwidth\undefined%
  \global\let\svgscale\undefined%
  \makeatother%
  \begin{picture}(1,0.88836483)%
    \lineheight{1}%
    \setlength\tabcolsep{0pt}%
    \put(0,0){\includegraphics[width=\unitlength,page=1]{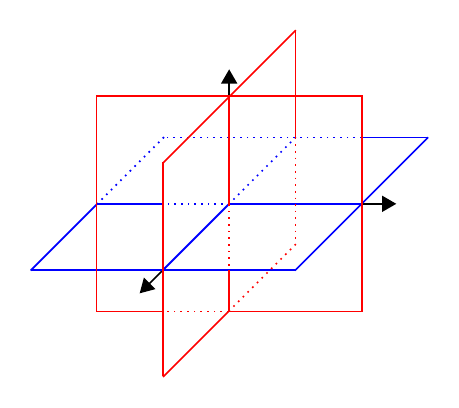}}%
    \put(0.06216741,0.34268715){\color[rgb]{0,0,0}\makebox(0,0)[lt]{\lineheight{0}\smash{\begin{tabular}[t]{l}$\mathcal{R}$\end{tabular}}}}%
    \put(0.65489788,0.77154051){\color[rgb]{0,0,0}\makebox(0,0)[lt]{\lineheight{0}\smash{\begin{tabular}[t]{l}$\mathcal{D}$\end{tabular}}}}%
    \put(0.50429777,0.74172633){\color[rgb]{0,0,0}\makebox(0,0)[lt]{\lineheight{0}\smash{\begin{tabular}[t]{l}$z$\end{tabular}}}}%
    \put(0.87098392,0.43679323){\color[rgb]{0,0,0}\makebox(0,0)[lt]{\lineheight{0}\smash{\begin{tabular}[t]{l}$y$\end{tabular}}}}%
    \put(0.26794436,0.21420459){\color[rgb]{0,0,0}\makebox(0,0)[lt]{\lineheight{0}\smash{\begin{tabular}[t]{l}$x$\end{tabular}}}}%
  \end{picture}%
\endgroup%

    \caption{$\mathcal{R}$ in blue, $\mathcal{D}$ in red and $\mathcal{T} =
        \mathcal{R} \cup \mathcal{D}$.}
    \label{sec:problem:fig:manifolds:subfig:rdt}
  \end{subfigure}
  \qquad
  \begin{subfigure}[b]{0.35\linewidth}
    \def\svgwidth{\linewidth}
    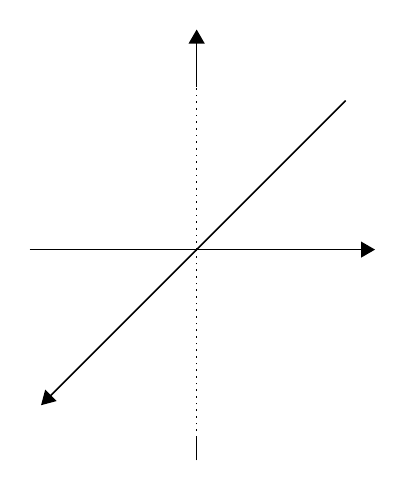
    \caption{$\mathcal{C}$ in red.}
    \label{sec:problem:fig:manifolds:subfig:c}
  \end{subfigure}
  \qquad
  \begin{subfigure}[b]{0.4\linewidth}
    \def\svgwidth{\linewidth}
    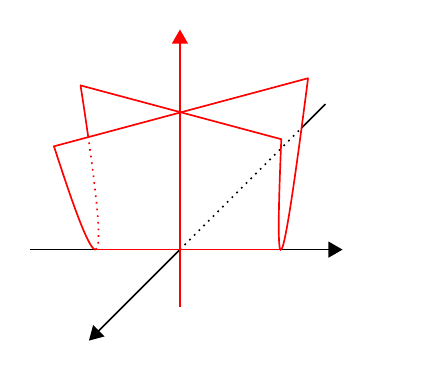
    \caption{$\mathcal{W}$ in red.}
    \label{sec:problem:fig:manifolds:subfig:w}
  \end{subfigure}
  \caption{Gutierrez-Sotomayor algebraic manifolds.}
  \label{sec:problem:fig:manifolds}
\end{figure}

In this context, an important class of non-smooth dynamical systems in
$\mathbb{R}^3$ with singular switching manifold, known as Gutierrez-Sotomayor
and described in \cite{Gutierrez1982}, is obtained when the regularity condition
is broken in a dynamically stable manner. More precisely, in order to avoid
non-trivial recurrence on non-orientable manifolds, a restriction to $\Sigma$ is
imposed so that its smooth parts are either orientable or diffeomorphic to an
open set of $\mathbb{P}^2$ (projective plane), $\mathbb{K}^2$ (Klein's bottle)
or $G^2 = \mathbb{T}^2 \# \mathbb{P}^2$ (torus with cross-cap). After a proper
coordinates normalization, this restriction leads to five algebraic manifolds,
with a regular configuration

\begin{equation}
  \label{sec:problem:eq:regular_manifold}
  \mathcal{R} = \left\{(x,y,z) \in \mathbb{R}^3;~ z = 0 \right\}
\end{equation}

\noindent known as \textbf{regular discontinuity} and four singular
configurations given by

\begin{equation}
  \label{sec:problem:eq:singular_manifolds}
  \begin{split}
    \mathcal{D} & = \left\{(x,y,z) \in \mathbb{R}^3;~ xy = 0 \right\},          \\
    \mathcal{T} & = \left\{(x,y,z) \in \mathbb{R}^3;~ xyz = 0 \right\},         \\
    \mathcal{C} & = \left\{(x,y,z) \in \mathbb{R}^3;~ z^2-x^2-y^2 = 0 \right\}, \\
    \mathcal{W} & = \left\{(x,y,z) \in \mathbb{R}^3;~ zx^2-y^2 = 0 \right\},
  \end{split}
\end{equation}

\noindent and known as \textbf{double}, \textbf{triple}, \textbf{cone} and
\textbf{Whitney discontinuities}, respectively. See
\cref{sec:problem:fig:manifolds}.

\begin{figure}[ht]
  \centering
  \begin{subfigure}[b]{0.45\linewidth}
    \def\svgwidth{\linewidth}
\begingroup%
  \makeatletter%
  \providecommand\color[2][]{%
    \errmessage{(Inkscape) Color is used for the text in Inkscape, but the package 'color.sty' is not loaded}%
    \renewcommand\color[2][]{}%
  }%
  \providecommand\transparent[1]{%
    \errmessage{(Inkscape) Transparency is used (non-zero) for the text in Inkscape, but the package 'transparent.sty' is not loaded}%
    \renewcommand\transparent[1]{}%
  }%
  \providecommand\rotatebox[2]{#2}%
  \newcommand*\fsize{\dimexpr\f@size pt\relax}%
  \newcommand*\lineheight[1]{\fontsize{\fsize}{#1\fsize}\selectfont}%
  \ifx\svgwidth\undefined%
    \setlength{\unitlength}{321.95512426bp}%
    \ifx\svgscale\undefined%
      \relax%
    \else%
      \setlength{\unitlength}{\unitlength * \real{\svgscale}}%
    \fi%
  \else%
    \setlength{\unitlength}{\svgwidth}%
  \fi%
  \global\let\svgwidth\undefined%
  \global\let\svgscale\undefined%
  \makeatother%
  \begin{picture}(1,0.70540142)%
    \lineheight{1}%
    \setlength\tabcolsep{0pt}%
    \put(0,0){\includegraphics[width=\unitlength,page=1]{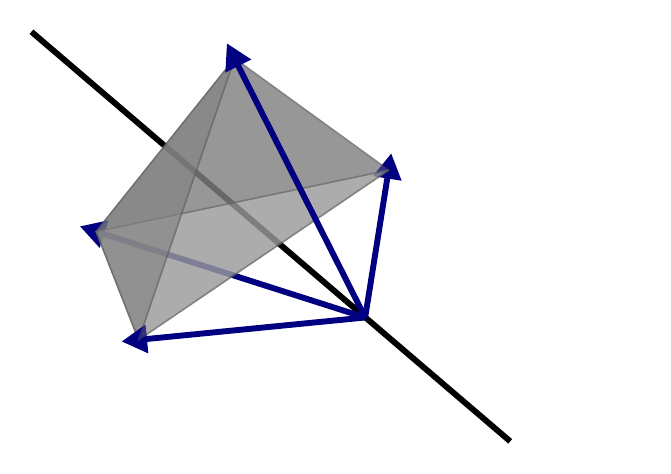}}%
    \put(0.74589882,0.09262335){\color[rgb]{0,0,0}\makebox(0,0)[lt]{\lineheight{1.25}\smash{\begin{tabular}[t]{l}$\Sigma_{\mathcal{S}}$\end{tabular}}}}%
    \put(0,0){\includegraphics[width=\unitlength,page=2]{hull.pdf}}%
  \end{picture}%
\endgroup%

    \caption{Hull.}
    \label{sec:introduction:subfig:hull}
  \end{subfigure}
  \quad
  \begin{subfigure}[b]{0.45\linewidth}
    \def\svgwidth{\linewidth}
\begingroup%
  \makeatletter%
  \providecommand\color[2][]{%
    \errmessage{(Inkscape) Color is used for the text in Inkscape, but the package 'color.sty' is not loaded}%
    \renewcommand\color[2][]{}%
  }%
  \providecommand\transparent[1]{%
    \errmessage{(Inkscape) Transparency is used (non-zero) for the text in Inkscape, but the package 'transparent.sty' is not loaded}%
    \renewcommand\transparent[1]{}%
  }%
  \providecommand\rotatebox[2]{#2}%
  \newcommand*\fsize{\dimexpr\f@size pt\relax}%
  \newcommand*\lineheight[1]{\fontsize{\fsize}{#1\fsize}\selectfont}%
  \ifx\svgwidth\undefined%
    \setlength{\unitlength}{321.95512426bp}%
    \ifx\svgscale\undefined%
      \relax%
    \else%
      \setlength{\unitlength}{\unitlength * \real{\svgscale}}%
    \fi%
  \else%
    \setlength{\unitlength}{\svgwidth}%
  \fi%
  \global\let\svgwidth\undefined%
  \global\let\svgscale\undefined%
  \makeatother%
  \begin{picture}(1,0.70540142)%
    \lineheight{1}%
    \setlength\tabcolsep{0pt}%
    \put(0,0){\includegraphics[width=\unitlength,page=1]{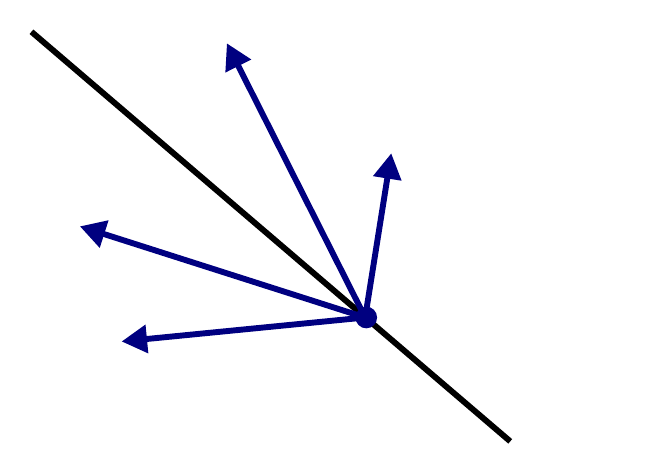}}%
    \put(0.74589882,0.09262335){\color[rgb]{0,0,0}\makebox(0,0)[lt]{\lineheight{1.25}\smash{\begin{tabular}[t]{l}$\Sigma_{\mathcal{S}}$\end{tabular}}}}%
    \put(0,0){\includegraphics[width=\unitlength,page=2]{canopy.pdf}}%
  \end{picture}%
\endgroup%

    \caption{Canopy.}
    \label{sec:introduction:subfig:canopy}
  \end{subfigure}
  \caption{Convex hull and canopy.}
  \label{sec:introduction:fig:hull_and_canopy}
\end{figure}

For systems whose switching manifold is homeomorphic to
\cref{sec:problem:eq:regular_manifold}, the Filippov dynamics described above is
fully applicable. In fact, these systems exactly corresponds to those described
at the beginning of this text. However, for systems whose switching manifold is
homeomorphic to one of the singular configurations
\cref{sec:problem:eq:singular_manifolds}, Filippov dynamics is not directly
applicable to the whole manifold $\Sigma$. More precisely, the switching
manifold can be decomposed in the following disjoint union:

\begin{equation}
  \label{sec:problem:eq:split_manifold}
  \Sigma = \Sigma_{\mathcal{R}} \cup \Sigma_{\mathcal{S}}
\end{equation}

\noindent where $\Sigma_{\mathcal{R}}$ consists, locally, of regular
discontinuities; and $\Sigma_{\mathcal{S}}$ consists of points where $\Sigma$
self-intersects, difficulting direct application of the usual Filippov dynamics.
In fact, an attempt to directly generalize the Filippov convention to points in
$\Sigma_{\mathcal{S}}$, leads to the existence of up to infinite possible
sliding fields, as proved at Lemma 2.4 of \cite[p.~1087]{Jeffrey2014}. See
\cref{sec:introduction:subfig:hull}.

In other words, the class of non-smooth dynamical systems $\mathbf{F} =
  \left(\mathbf{F}_i \right)$ whose switching manifold is homeomorphic to one of
those at \cref{sec:problem:eq:singular_manifolds}, in the sense of
Gutierrez-Sotomayor, represents the simplest singular systems. However,
besides its many applications, knowledge of its dynamics is scarce. In
particular, over the last decade, three main methodologies arose to study
these systems. Not necessarily in chronological order, these methodologies are
briefly presented below.

The first one, presented in \cite{Jeffrey2014} by \textit{Jeffrey}, propose an
extension of the Filippov dynamics to $\Sigma_{\mathcal{S}}$ through the
so-called ``canopy'', a convex-like ruled surface, built with the convex hull
$\text{Conv}(\{\mathbf{F}_i\})$, which can be proved to intersect $\Sigma$ at a
finite number of points, see \cref{sec:introduction:subfig:canopy}. Each one of
these intersections represents a sliding vector and, therefore, this methodology
leads again to non-uniqueness of the sliding field. To deal with this lack of
uniqueness, the author there conjectures the so-called ``dummy dynamics'' acting
over the canopy. This idea led to many results such as, for instance,
\cite{Jeffrey2018, Kaklamanos2019, Webber2020}. However, as stated at
\cite[p.~1102]{Jeffrey2014}, a justification for the dummy dynamics remains an
open problem.

The next one, presented in \cite{Dieci2009} by \textit{Dieci et al.}, although
older than the previous methodology, proposes a similar construction where,
again, non-uniqueness of sliding vectors happens. Here, however, the authors
show that, imposing certain attractivity hypothesis on the switching manifold
$\Sigma$, many conclusions can be proved on the behavior of the dynamics. In
fact, this idea led to the sequence of works \cite{Dieci2011,
  Dieci2013,Dieci2014, Dieci2015, Dieci2015b, Dieci2015c, Dieci2016} where several
aspects of the dynamics are explored under different types of attractivity: from
minimal sets to structural stability. However, imposing conditions on $\Sigma$
is a fundamental and restrictive hypothesis here.

Finally, \cite{Buzzi2012} by \textit{Buzzi et al.}, propose an extension of the
Filippov dynamics to $\Sigma_{\mathcal{S}}$ through the application of a proper
blow-up and use of Geometrical Singular Perturbation Theory (see
\cite{Fenichel1979,Teixeira2012}), or GSP-Theory for short, to study the
resulting slow-fast systems. Although distant from a direct generalization of
Filippov's convention, this methodology is also a natural approach with
advantages over the previous ones. In fact, while the non-uniqueness of the
sliding field is also predicted, here it is managed naturally, as will be seen
over the text. Moreover, yet in comparison with the previous ones, due to the
blow-up, this methodology provides a broader view of the dynamics -- appropriate
for the codimension of the problem. Even more, no assumptions neither on
$\Sigma$ or the underlying vector fields $\mathbf{F}_i$ are required here.
However, both \cite{Buzzi2012} and the posterior works \cite{Llibre2015,
  Panazzolo2017, Teixeira2012} lack a clear presentation and justification for the
dynamics induced over $\Sigma_{\mathcal{S}}$. Up to our knowledge, there are no
works focused on the study of this dynamics for any class of fields
$\mathbf{F}_i$ such as, for instance, linear ones. In fact, the main focus of
the above-cited works lies on the verification that, after the blow-up, the
resulting system contains only regular discontinuities or, in other words,
\cref{sec:framework:sub:blowup:thm:regular_discontinuities} with some sparse
examples lacking proper justification for the underlying dynamics: see, for
instance, the examples in \cite[p.~449]{Buzzi2012}, \cite[p.
  1952]{Teixeira2012}, and the remark and example in \cite[p.~501]{Llibre2015}.

Given the arguments above, for this text, we embrace and leverage \emph{Buzzi's}
blow-up based methodology to study the dynamics associated with singular
switching manifolds, since it

\begin{enumerate}
  \item does not depend on imposing conditions on $\Sigma$;
  \item deals naturally with the non-uniqueness of sliding vectors; and
  \item provide a broader view of the dynamics over $\Sigma_{\mathcal{S}}$.
\end{enumerate}

\begin{figure}[ht]
  \centering
  \begin{subfigure}[b]{0.42\linewidth}
    \def\svgwidth{\linewidth}
\begingroup%
  \makeatletter%
  \providecommand\color[2][]{%
    \errmessage{(Inkscape) Color is used for the text in Inkscape, but the package 'color.sty' is not loaded}%
    \renewcommand\color[2][]{}%
  }%
  \providecommand\transparent[1]{%
    \errmessage{(Inkscape) Transparency is used (non-zero) for the text in Inkscape, but the package 'transparent.sty' is not loaded}%
    \renewcommand\transparent[1]{}%
  }%
  \providecommand\rotatebox[2]{#2}%
  \newcommand*\fsize{\dimexpr\f@size pt\relax}%
  \newcommand*\lineheight[1]{\fontsize{\fsize}{#1\fsize}\selectfont}%
  \ifx\svgwidth\undefined%
    \setlength{\unitlength}{234.06075785bp}%
    \ifx\svgscale\undefined%
      \relax%
    \else%
      \setlength{\unitlength}{\unitlength * \real{\svgscale}}%
    \fi%
  \else%
    \setlength{\unitlength}{\svgwidth}%
  \fi%
  \global\let\svgwidth\undefined%
  \global\let\svgscale\undefined%
  \makeatother%
  \begin{picture}(1,1)%
    \lineheight{1}%
    \setlength\tabcolsep{0pt}%
    \put(0,0){\includegraphics[width=\unitlength,page=1]{slice.pdf}}%
    \put(0.73501026,0.7354913){\color[rgb]{0,0,0}\makebox(0,0)[lt]{\lineheight{0}\smash{\begin{tabular}[t]{l}$\mathbf{F}_1$\end{tabular}}}}%
    \put(0.24033606,0.24061491){\color[rgb]{0,0,0}\makebox(0,0)[lt]{\lineheight{0}\smash{\begin{tabular}[t]{l}$\mathbf{F}_3$\end{tabular}}}}%
    \put(0.73646534,0.2407504){\color[rgb]{0,0,0}\makebox(0,0)[lt]{\lineheight{0}\smash{\begin{tabular}[t]{l}$\mathbf{F}_4$\end{tabular}}}}%
    \put(0.24033606,0.73681099){\color[rgb]{0,0,0}\makebox(0,0)[lt]{\lineheight{0}\smash{\begin{tabular}[t]{l}$\mathbf{F}_2$\end{tabular}}}}%
    \put(0.52429262,0.41253198){\color[rgb]{0,0,0}\makebox(0,0)[lt]{\lineheight{0}\smash{\begin{tabular}[t]{l}$\Sigma_x$\end{tabular}}}}%
  \end{picture}%
\endgroup%

    \caption{$\mathcal{D}_2^k$}
    \label{sec:introduction:subfig:planar_double}
  \end{subfigure}
  \quad
  \begin{subfigure}[b]{0.53\linewidth}
    \def\svgwidth{\linewidth}
\begingroup%
  \makeatletter%
  \providecommand\color[2][]{%
    \errmessage{(Inkscape) Color is used for the text in Inkscape, but the package 'color.sty' is not loaded}%
    \renewcommand\color[2][]{}%
  }%
  \providecommand\transparent[1]{%
    \errmessage{(Inkscape) Transparency is used (non-zero) for the text in Inkscape, but the package 'transparent.sty' is not loaded}%
    \renewcommand\transparent[1]{}%
  }%
  \providecommand\rotatebox[2]{#2}%
  \newcommand*\fsize{\dimexpr\f@size pt\relax}%
  \newcommand*\lineheight[1]{\fontsize{\fsize}{#1\fsize}\selectfont}%
  \ifx\svgwidth\undefined%
    \setlength{\unitlength}{232.74925412bp}%
    \ifx\svgscale\undefined%
      \relax%
    \else%
      \setlength{\unitlength}{\unitlength * \real{\svgscale}}%
    \fi%
  \else%
    \setlength{\unitlength}{\svgwidth}%
  \fi%
  \global\let\svgwidth\undefined%
  \global\let\svgscale\undefined%
  \makeatother%
  \begin{picture}(1,0.8396019)%
    \lineheight{1}%
    \setlength\tabcolsep{0pt}%
    \put(0,0){\includegraphics[width=\unitlength,page=1]{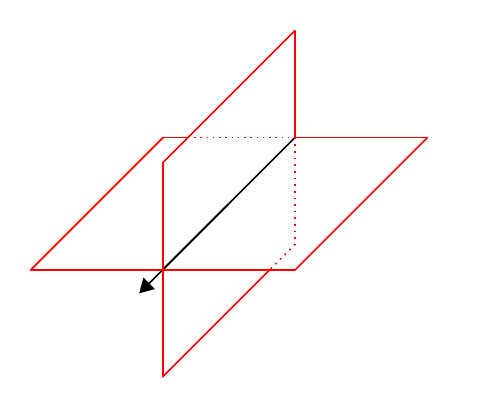}}%
    \put(0.23767726,0.18552944){\color[rgb]{0,0,0}\makebox(0,0)[lt]{\lineheight{0}\smash{\begin{tabular}[t]{l}$\Sigma_x$\end{tabular}}}}%
    \put(0.69686239,0.63021399){\color[rgb]{0,0,0}\makebox(0,0)[lt]{\lineheight{0}\smash{\begin{tabular}[t]{l}$\mathbf{F}_1$\end{tabular}}}}%
    \put(0.34482043,0.63021399){\color[rgb]{0,0,0}\makebox(0,0)[lt]{\lineheight{0}\smash{\begin{tabular}[t]{l}$\mathbf{F}_2$\end{tabular}}}}%
    \put(0.11764466,0.14525004){\color[rgb]{0,0,0}\makebox(0,0)[lt]{\lineheight{0}\smash{\begin{tabular}[t]{l}$\mathbf{F}_3$\end{tabular}}}}%
    \put(0.5687739,0.14847239){\color[rgb]{0,0,0}\makebox(0,0)[lt]{\lineheight{0}\smash{\begin{tabular}[t]{l}$\mathbf{F}_4$\end{tabular}}}}%
  \end{picture}%
\endgroup%

    \caption{$\mathcal{D}_3^k$}
    \label{sec:introduction:subfig:tridimensional_double}
  \end{subfigure}
  \caption{Double discontinuity.}
  \label{sec:introduction:fig:double_discontinuity}
\end{figure}

Specifically, we deal essentially with the Gutierrez\hyp{}Sotomayor algebraic
manifold $\mathcal{D}$, the double discontinuity, both an equivalent in
$\mathbb{R}^2$ and the traditional in $\mathbb{R}^3$, whose classes of vector
fields are henceforth denoted $\mathcal{D}_2^k$ and $\mathcal{D}_3^k$,
respectively, see \cref{sec:introduction:fig:double_discontinuity}.
Geometrically, these configurations represent transversal self-intersections of
the switching manifold. After a general qualitative study on the non-linear
case, we focus on systems given by affine vector fields

\begin{equation*}
  \mathbf{F}_i(\mathbf{x}) =
  \mathbf{A}_i \mathbf{x} + \mathbf{b}_i,
\end{equation*}

\noindent where $\mathbf{A}_i$ and $\mathbf{b}_i$ are real matrices for every $i
  \in \{1,2,3,4\}$ respectively of sizes $j \times j$ and $j \times 1$, with $j
  \in \left\{ 2, 3 \right\}$ representing the dimension and progressively
increasing its complexity: starting at the constant case ($\mathbf{A}_i =
  \mathbf{0}$), linear ($\mathbf{b}_i = \mathbf{0}$) and then, finally, the
complete affine case. We denote these classes of constant, linear and affine
vector fields, respectively, as $\mathcal{C}_j$, $\mathcal{L}_j$ and
$\mathcal{A}_j$. This program assures a progressive and effective increase on
the intuition and understanding of the dynamics.

For the class $\mathcal{D}_3^k$ of piecewise vector fields having the
traditional Double Discontinuity as switching manifold, see
\cref{sec:introduction:subfig:tridimensional_double}, we tackle the main problem
of defining a dynamics over $\Sigma_{\mathcal{S}}$, here given by a straight
line $\Sigma_{x}$. In particular, we use the cylindrical blow-up suggested in
\cite{Llibre2015} to induce a dynamics over $\Sigma_{x}$, with GSP-Theory
playing a major role. As a result, we obtain the \textbf{Fundamental
  \cref{sec:framework:blowup:thm:dynamics}} stated below, which improves
\emph{Buzzi's} methodology by clarifying the issues raised above:

\begin{restatable*}[Fundamental Dynamics]{rlemma}{LemmaFundamental}
  \label{sec:framework:blowup:thm:dynamics}
  Given $\mathbf{F} \in \mathcal{D}_3^k$ with components $\mathbf{F}_i = (w_i,
    p_i, q_i)$, let $\mathbf{\tilde{F}} \in \tilde{\mathcal{D}}_3^k$ be the
  vector field induced by the blow-up

  \begin{equation*}
    \phi_1(x,\theta,r) = (x, r\cos{\theta}, r\sin{\theta}).
  \end{equation*}

  \noindent Then, this blow-up associates the dynamics over $\Sigma_x$ with the
  following dynamics over the cylinder $C = \mathbb{R} \times S^1 = S_1 \cup
    \ldots \cup S_4$: over each stripe $S_i$ acts a slow-fast dynamics whose
  reduced dynamics is given by

  \begin{equation}
    \label{sec:framework:blowup:thm:dynamics:eq:reduced}
    \left\{\begin{matrix*}[l]
      \dot{x} = w_i\\
      0 = q_i \cos{\theta} - p_i \sin{\theta}\\
    \end{matrix*}\right.,
  \end{equation}

  \noindent with slow radial dynamics $\dot{r} = p_i \cos{\theta} + q_i
    \sin{\theta}$; and layer dynamics given by

  \begin{equation}
    \label{sec:framework:blowup:thm:dynamics:eq:layer}
    \left\{\begin{matrix*}[l]
      x' = 0\\
      \theta' = q_i \cos{\theta} - p_i \sin{\theta}
    \end{matrix*}\right.,
  \end{equation}

  \noindent with fast radial dynamics $r' = 0$. Finally, at every equation above
  the functions $w_i$, $p_i$ and $q_i$ must be calculated at the point
  $\phi_1(x,\theta,0) = (x,0,0)$.
\end{restatable*}

We note here that no blowing-down is ever carried out over this text, i.e., once
we have the blow-up induced (fundamental) dynamics above over the cylinder $C$,
the inverse operation is never performed to recover a dynamics over $\Sigma_{x}$
with the original coordinates. Actually, this operation would make little to no
sense most of the times given the higher codimension of $\Sigma_{x}$, i.e., most
of the information on the dynamics would be lost. For instance, under the
so-called \textbf{fundamental hypothesis}
\cref{sec:framework:sub:blowup:eq:weak_fund_hypo} or
\cref{sec:framework:sub:blowup:eq:strong_fund_hypo}, the fundamental lemma
assures not only the sequence of qualitative theorems bellow for the general,
non-linear case, but also most of the original results in this text and, hence,
its name.

\begin{restatable*}{rtheorem}{TheoremTransversality}
  \label{sec:framework:sub:blowup:cor:radial_dynamics}
  The radial dynamics can only be transversal ($\dot{r} \neq 0$)  to the
  cylinder $C$ over the slow manifold $\mathcal{M}_i$. More over, under
  \cref{sec:framework:sub:blowup:eq:weak_fund_hypo}, it is in fact transversal.
\end{restatable*}

\begin{restatable*}{rtheorem}{TheoremGraph}
  \label{sec:framework:sub:blowup:cor:explicit_slow_manifold}
  The slow manifold $\mathcal{M}_i$ is locally a graph \\ $\left( x, \theta(x)
    \right)$ under \cref{sec:framework:sub:blowup:eq:weak_fund_hypo}. However,
  if $\norm{(f_i)_{\theta}}$ admits a global positive minimum, then
  $\mathcal{M}_i$ is globally a graph $\left( x, \theta(x) \right)$. Either
  way, $\theta(x)$ is of class $C^k$.
\end{restatable*}

\begin{restatable*}{rtheorem}{TheoremNormalHiperbolicity}
  \label{sec:framework:sub:blowup:cor:normal_hyperbolic}
  The slow manifold $\mathcal{M}_i$ is normally hyperbolic at every point that
  satisfies \cref{sec:framework:sub:blowup:eq:weak_fund_hypo}.
\end{restatable*}

\begin{restatable*}{rtheorem}{TheoremHiperbolicSingularities}
  \label{sec:framework:sub:blowup:cor:singularities}
  The hyperbolic singularities of the reduced system
  \cref{sec:framework:blowup:thm:dynamics:eq:reduced} acts as hyperbolic saddle
  or node singularities of $S_i$ under
  \cref{sec:framework:sub:blowup:eq:weak_fund_hypo}.
\end{restatable*}

Then, using these fundamental dynamics, we focus on the constant
($\mathcal{C}_3$) and affine ($\mathcal{A}_3$) cases, to fully describe the
respective induced dynamics over the cylinder as stated in
\cref{sec:constant:thm:dynamics} for the constant case and
\cref{sec:affine:thm:dynamics} below for the affine case.

\begin{restatable*}[Affine Dynamics]{rtheorem}{TheoremAffineDynamics}
  \label{sec:affine:thm:dynamics}
  Given $\mathbf{F} \in \mathcal{A}_3$ with affine components $\mathbf{F}_i$
  given by \eqref{sec:affine:eq:system} and such that $\gamma_i \neq 0$, let
  $\mathbf{\tilde{F}} \in \tilde{\mathcal{A}_3}$ be the vector field induced by
  the blow-up $\phi_1(x,\theta,r) = (x, r\cos{\theta}, r\sin{\theta})$. Then,
  this blow-up associates the dynamics over $\Sigma_x$ with the following
  fundamental dynamics over the cylinder $C = \mathbb{R} \times S^1 = S_1 \cup
    \ldots \cup S_4$: over each stripe $S_i$ acts a slow-fast dynamics whose slow
  manifold is given by $\mathcal{M}_i = A_{i} \cup A_{i}^{\pi}$, where
  $A_{i}^{\pi}$ is a $\pi$-translation of $A_{i}$ in $\theta$ and

  \begin{enumerate}
    \item \label{sec:affine:thm:dynamics:a2_nonzero} case $a_{i2} \neq 0$, then

          \begin{align*}
            A_{i} & = \left\{(x,\theta) \in [-\infty,\alpha_i] \times
            \left[0,2\pi\right];~
            \theta = \theta_i(x) + \pi \right\} \cup                     \\
                  & \cup \left\{(x,\theta) \in [\alpha_i,+\infty] \times
            \left[0, 2\pi\right];~ \theta = \theta_i(x) \right\}
          \end{align*}

          \noindent with $\theta_i(x) = \arctan{\left(\frac{a_{i3}x +
                d_{i3}}{a_{i2}x + d_{i2}}\right)}$, which consists of an
          arctangent\hyp{}like curve inside the cylinder C with $\theta =
            \beta_i + \pi$ and $\theta = \beta_i$ as negative and positive
          horizontal asymptotes, respectively;

    \item \label{sec:affine:thm:dynamics:a2_zero} case $a_{i2} = 0$, then

          \begin{equation*}
            A_{i} = \left\{(x,\theta) \in \mathbb{R} \times \left[0,2\pi\right] ;~
            \theta = \theta_i(x) \right\}
          \end{equation*}

          \noindent with $\theta_i(x) = \arctan{\left(\frac{a_{i3}x +
                d_{i3}}{d_{i2}}\right)}$, which consists of an
          arctangent\hyp{}like curve inside the cylinder C with $\theta =
            \sigma_{i-}$ and $\theta = \sigma_{i+}$ as negative and positive
          horizontal asymptotes, respectively.
  \end{enumerate}

  Both arctangents are increasing if $\gamma_i > 0$ and decreasing if $\gamma_i
    < 0$. Over them act the reduced dynamics $\dot{x} = a_{i1}x + d_{i1}$ and,
  around them, acts the layer dynamics described in
  \cref{sec:affine:tab:fast}, but exchanging $a_{i2}$ with $d_{i2}$ if $a_{i2}
    = 0$. Finally, the new parameters above are given by $\alpha_i =
    -\frac{d_{i2}}{a_{i2}}$, $\beta_i = \arctan{\left(\frac{a_{i3}}{a_{i2}}
      \right)}$, $\gamma_i = a_{i3}d_{i2} - d_{i3}a_{i2}$, $\delta_i =
    -\frac{d_{i1}}{a_{i1}}$ and $\sigma_{i\pm} = \pm\sgn{\left(\gamma_i
      \right)}\frac{\pi}{2}$.
\end{restatable*}

Finally, combining this fine-grained control of the fundamental dynamics with
the structural stability characterization provided by \cite{Broucke2001}, we
also derive Peixoto-like theorems characterizing semi-local structural stability
of the dynamics over the cylinder for both the constant
(\cref{sec:stability:sub:constant:thm:conditions}) and affine
(\cref{sec:stability:sub:affine:thm:conditions} as stated below) cases.

\begin{restatable*}[Affine Dynamics Stability]{rtheorem}{TheoremAffineStability}
  \label{sec:stability:sub:affine:thm:conditions}
  Let $\mathbf{F} \in \mathcal{A}_3$ be given by
  \cref{sec:stability:sub:affine:eq:system} with $\gamma_{i} \neq 0$. Given
  $\Sigma_{\theta_0} \in \tilde{\mathcal{I}}_C$, let $\mathbf{X} =
    (\mathbf{X}_{-}, \mathbf{X}_{+})$ be the Filippov system induced around
  $\Sigma_{\theta_0}$ and inside a convex compact set $K \subset C_{+} \cup
    C_{-}$, where $C_{+}$ and $C_{-}$ are two consecutive stripes meeting at
  $\Sigma_{\theta_0}$. Then, $\mathbf{F}$ is $(\Sigma_{\theta_0}, K)$-semi-local
  structurally stable in $\mathcal{A}_3$ if, and only if, $\mathbf{X}_{+}$ and
  $\mathbf{X}_{-}$ satisfies

  \begin{enumerate}
    \item $a_{i1} \neq 0$ and $\mathbf{P} \not\in \Sigma_{\theta_0}$, where
          $\mathbf{P}$ is the only singularity of $\mathbf{X}_{\pm}$;
    \item conditions \cref{sec:stability:conditions:sg_colinear} —
          \cref{sec:stability:conditions:cr_recurrent} of
          \cref{sec:stability:prop:conditions}.
  \end{enumerate}
\end{restatable*}

For clarification on the technicalities involved, especially in the statement
above, we encourage the interested reader to please consult the respective
section in the text, which is structured as follows: preliminaries are presented
in \cref{sec:preliminaries}, followed by a formal statement of the problem in
\cref{sec:problem} and then the methodology in \cref{sec:framework} -- general
results for the non-linear case are also contained here as consequences of the
fundamental dynamics obtained. Next, when generated by constant and affine
vector fields, a complete qualitative description of the fundamental dynamics
are obtained in \cref{sec:constant,sec:affine}, respectively. Finally, in
\cref{sec:stability}, this fine-grained control of the fundamental dynamics is
leveraged to derive Peixoto-like theorems characterizing semi-local structural
stability. Conclusion, in \cref{sec:conclusion}, contains some final thoughts
and further directions of investigation.

\section{Preliminaries}%
\label{sec:preliminaries}

In this section, we introduce the concept of regularization developed by Jorge
Sotomayor and Marco Antonio Teixeira in \cite{Sotomayor1996}. First, given a
piecewise vector field, we construct its regularization, which consists of a
1-parameter family of smooth vector fields which converges to the given
piecewise field. Next, we introduce part of the Geometrical Singular
Perturbation Theory developed by Neil Fenichel in \cite{Fenichel1979} and its
connection with Filippov dynamics through regularization.

\subsection{Sotomayor-Teixeira regularization}
\label{sec:sotomayor_teixeira}

\begin{figure}[ht]
  \centering
  \begin{subfigure}[b]{0.45\linewidth}
    \def\svgwidth{\linewidth}
\begingroup%
  \makeatletter%
  \providecommand\color[2][]{%
    \errmessage{(Inkscape) Color is used for the text in Inkscape, but the package 'color.sty' is not loaded}%
    \renewcommand\color[2][]{}%
  }%
  \providecommand\transparent[1]{%
    \errmessage{(Inkscape) Transparency is used (non-zero) for the text in Inkscape, but the package 'transparent.sty' is not loaded}%
    \renewcommand\transparent[1]{}%
  }%
  \providecommand\rotatebox[2]{#2}%
  \newcommand*\fsize{\dimexpr\f@size pt\relax}%
  \newcommand*\lineheight[1]{\fontsize{\fsize}{#1\fsize}\selectfont}%
  \ifx\svgwidth\undefined%
    \setlength{\unitlength}{385.64645842bp}%
    \ifx\svgscale\undefined%
      \relax%
    \else%
      \setlength{\unitlength}{\unitlength * \real{\svgscale}}%
    \fi%
  \else%
    \setlength{\unitlength}{\svgwidth}%
  \fi%
  \global\let\svgwidth\undefined%
  \global\let\svgscale\undefined%
  \makeatother%
  \begin{picture}(1,0.67061215)%
    \lineheight{1}%
    \setlength\tabcolsep{0pt}%
    \put(0,0){\includegraphics[width=\unitlength,page=1]{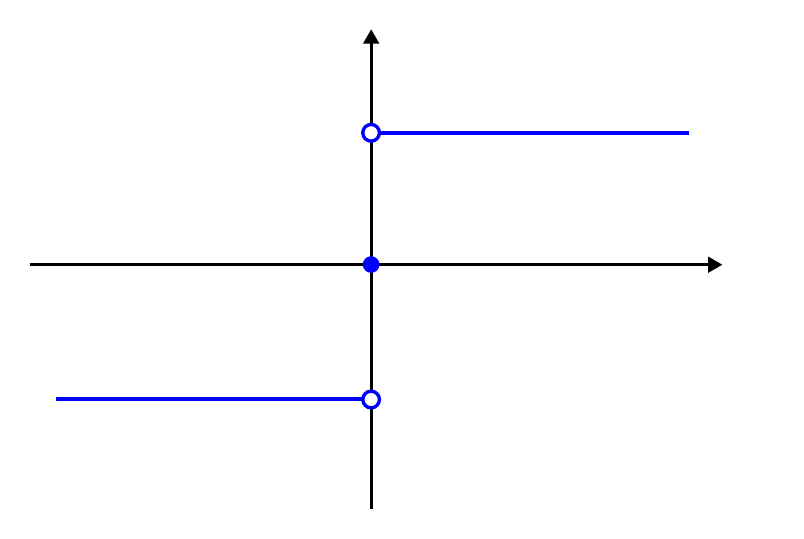}}%
    \put(0.90780645,0.32885765){\color[rgb]{0,0,0}\makebox(0,0)[lt]{\lineheight{0}\smash{\begin{tabular}[t]{l}$x$\end{tabular}}}}%
    \put(0.48519317,0.60919805){\color[rgb]{0,0,0}\makebox(0,0)[lt]{\lineheight{0}\smash{\begin{tabular}[t]{l}$\sgn(x)$\end{tabular}}}}%
    \put(0.41250621,0.48919816){\color[rgb]{0,0,0}\makebox(0,0)[lt]{\lineheight{0}\smash{\begin{tabular}[t]{l}$1$\end{tabular}}}}%
    \put(0.48085387,0.15881381){\color[rgb]{0,0,0}\makebox(0,0)[lt]{\lineheight{0}\smash{\begin{tabular}[t]{l}$-1$\end{tabular}}}}%
  \end{picture}%
\endgroup%

    \caption{Signal function.}%
    \label{sec:sotomayor_teixeira:subfig:signal}
  \end{subfigure}
  \qquad
  \begin{subfigure}[b]{0.45\linewidth}
    \def\svgwidth{\linewidth}
\begingroup%
  \makeatletter%
  \providecommand\color[2][]{%
    \errmessage{(Inkscape) Color is used for the text in Inkscape, but the package 'color.sty' is not loaded}%
    \renewcommand\color[2][]{}%
  }%
  \providecommand\transparent[1]{%
    \errmessage{(Inkscape) Transparency is used (non-zero) for the text in Inkscape, but the package 'transparent.sty' is not loaded}%
    \renewcommand\transparent[1]{}%
  }%
  \providecommand\rotatebox[2]{#2}%
  \newcommand*\fsize{\dimexpr\f@size pt\relax}%
  \newcommand*\lineheight[1]{\fontsize{\fsize}{#1\fsize}\selectfont}%
  \ifx\svgwidth\undefined%
    \setlength{\unitlength}{385.64645842bp}%
    \ifx\svgscale\undefined%
      \relax%
    \else%
      \setlength{\unitlength}{\unitlength * \real{\svgscale}}%
    \fi%
  \else%
    \setlength{\unitlength}{\svgwidth}%
  \fi%
  \global\let\svgwidth\undefined%
  \global\let\svgscale\undefined%
  \makeatother%
  \begin{picture}(1,0.67061215)%
    \lineheight{1}%
    \setlength\tabcolsep{0pt}%
    \put(0,0){\includegraphics[width=\unitlength,page=1]{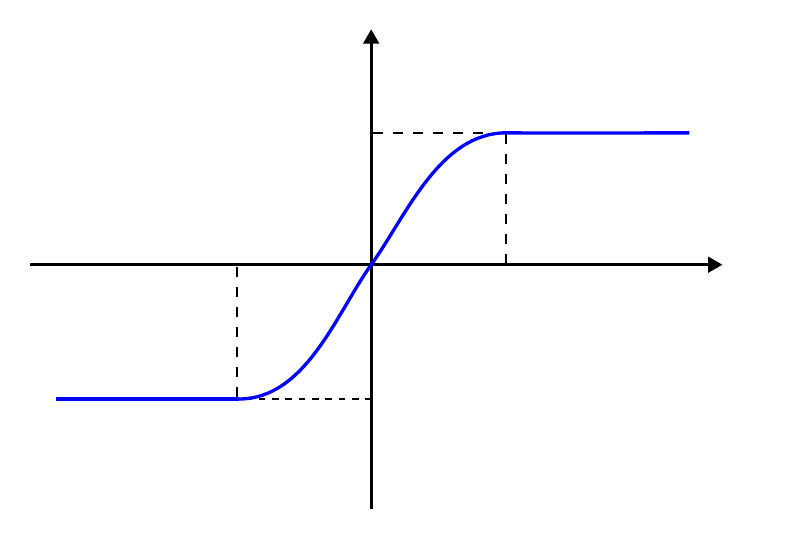}}%
    \put(0.90780645,0.32885765){\color[rgb]{0,0,0}\makebox(0,0)[lt]{\lineheight{0}\smash{\begin{tabular}[t]{l}$x$\end{tabular}}}}%
    \put(0.48519317,0.60919805){\color[rgb]{0,0,0}\makebox(0,0)[lt]{\lineheight{0}\smash{\begin{tabular}[t]{l}$\varphi(x)$\end{tabular}}}}%
    \put(0.41250621,0.48919816){\color[rgb]{0,0,0}\makebox(0,0)[lt]{\lineheight{0}\smash{\begin{tabular}[t]{l}$1$\end{tabular}}}}%
    \put(0.48085387,0.15881381){\color[rgb]{0,0,0}\makebox(0,0)[lt]{\lineheight{0}\smash{\begin{tabular}[t]{l}$-1$\end{tabular}}}}%
    \put(0.19146557,0.28116075){\color[rgb]{0,0,0}\makebox(0,0)[lt]{\lineheight{0}\smash{\begin{tabular}[t]{l}$-1$\end{tabular}}}}%
    \put(0.61801708,0.28128831){\color[rgb]{0,0,0}\makebox(0,0)[lt]{\lineheight{0}\smash{\begin{tabular}[t]{l}$1$\end{tabular}}}}%
  \end{picture}%
\endgroup%

    \caption{Transition function.}%
    \label{sec:sotomayor_teixeira:subfig:transition}
  \end{subfigure}

  \bigskip

  \begin{subfigure}[b]{0.4\linewidth}
    \def\svgwidth{\linewidth}
\begingroup%
  \makeatletter%
  \providecommand\color[2][]{%
    \errmessage{(Inkscape) Color is used for the text in Inkscape, but the package 'color.sty' is not loaded}%
    \renewcommand\color[2][]{}%
  }%
  \providecommand\transparent[1]{%
    \errmessage{(Inkscape) Transparency is used (non-zero) for the text in Inkscape, but the package 'transparent.sty' is not loaded}%
    \renewcommand\transparent[1]{}%
  }%
  \providecommand\rotatebox[2]{#2}%
  \newcommand*\fsize{\dimexpr\f@size pt\relax}%
  \newcommand*\lineheight[1]{\fontsize{\fsize}{#1\fsize}\selectfont}%
  \ifx\svgwidth\undefined%
    \setlength{\unitlength}{710.10460561bp}%
    \ifx\svgscale\undefined%
      \relax%
    \else%
      \setlength{\unitlength}{\unitlength * \real{\svgscale}}%
    \fi%
  \else%
    \setlength{\unitlength}{\svgwidth}%
  \fi%
  \global\let\svgwidth\undefined%
  \global\let\svgscale\undefined%
  \makeatother%
  \begin{picture}(1,0.87964354)%
    \lineheight{1}%
    \setlength\tabcolsep{0pt}%
    \put(0,0){\includegraphics[width=\unitlength,page=1]{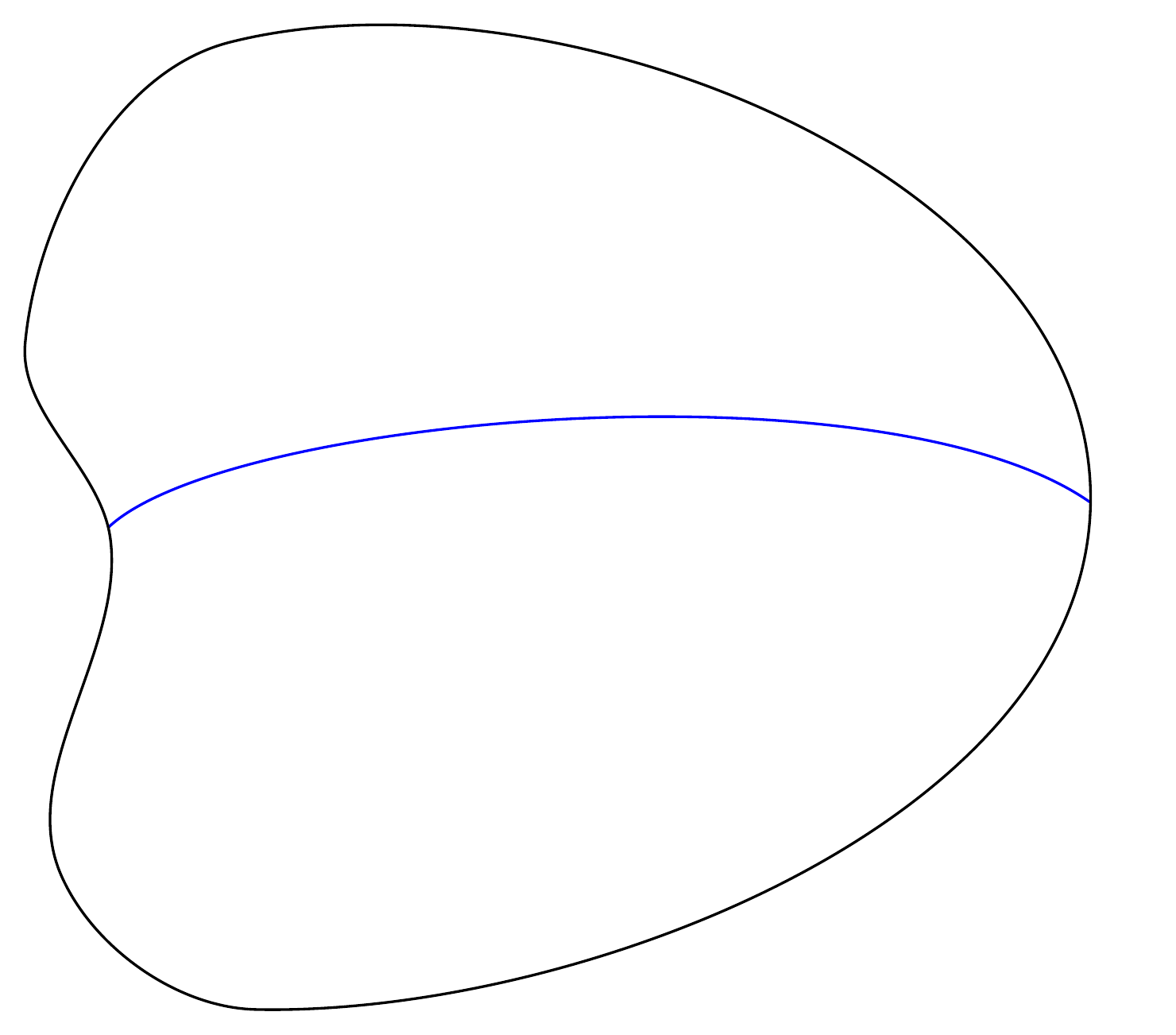}}%
    \put(0.84871925,0.40574074){\color[rgb]{0,0,0}\makebox(0,0)[lt]{\lineheight{0}\smash{\begin{tabular}[t]{l}$\Sigma$\end{tabular}}}}%
    \put(0.46728659,0.63749737){\color[rgb]{0,0,0}\makebox(0,0)[lt]{\lineheight{0}\smash{\begin{tabular}[t]{l}$\mathbf{F_+}$\end{tabular}}}}%
    \put(0.45850608,0.33939887){\color[rgb]{0,0,0}\makebox(0,0)[lt]{\lineheight{0}\smash{\begin{tabular}[t]{l}$\mathbf{F_-}$\end{tabular}}}}%
    \put(0.09466691,0.61846134){\color[rgb]{0,0,0}\makebox(0,0)[lt]{\lineheight{0}\smash{\begin{tabular}[t]{l}$\Sigma_+$\end{tabular}}}}%
    \put(0.24518648,0.5618546){\color[rgb]{0,0,0}\makebox(0,0)[lt]{\lineheight{0}\smash{\begin{tabular}[t]{l} \end{tabular}}}}%
    \put(0.12417874,0.2482945){\color[rgb]{0,0,0}\makebox(0,0)[lt]{\lineheight{0}\smash{\begin{tabular}[t]{l}$\Sigma_-$\end{tabular}}}}%
    \put(0.02142275,0.76980316){\color[rgb]{0,0,0}\makebox(0,0)[lt]{\lineheight{0}\smash{\begin{tabular}[t]{l}$U$\end{tabular}}}}%
  \end{picture}%
\endgroup%

    \caption{Piecewise field.}%
    \label{sec:sotomayor_teixeira:subfig:piecewise}
  \end{subfigure}
  \qquad
  \begin{subfigure}[b]{0.4\linewidth}
    \def\svgwidth{\linewidth}
    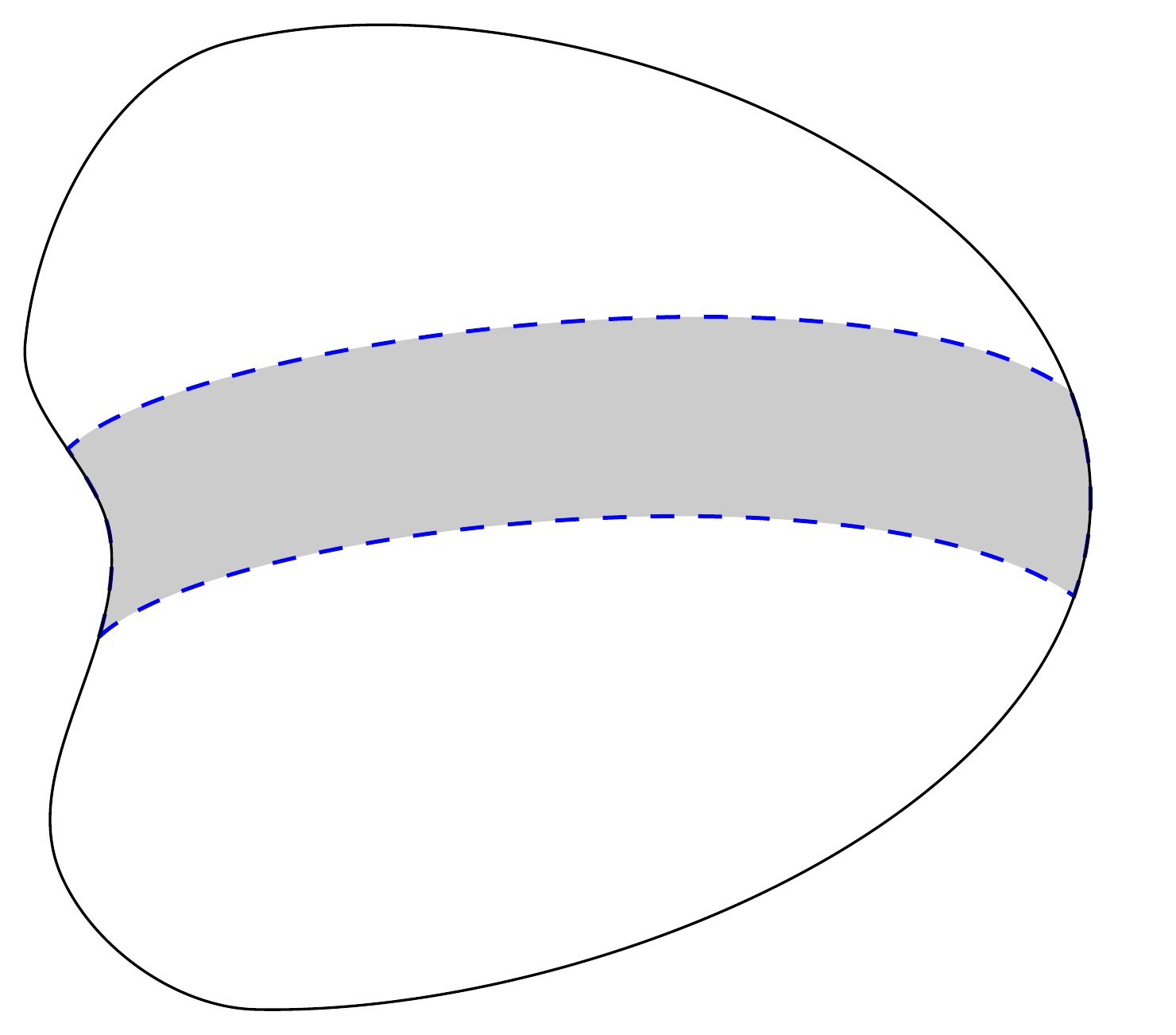
    \caption{Regularized field.}%
    \label{sec:sotomayor_teixeira:subfig:regularization}
  \end{subfigure}
  \caption{Grid representation of the Sotomayor-Teixeira's regularization with
    the signal function (a) associated to the piecewise smooth vector field (c)
    and the transition function (b) associated to the regularized vector field
    (d).}%
  \label{sec:sotomayor_teixeira:fig:regularization_grid}
\end{figure}

Let $\mathbf{F} = (\mathbf{F}_+,\mathbf{F}_-) \in \mathcal{R}^k(U)$ be a
piecewise smooth vector field as defined above with a switching manifold $\Sigma
  = h^{-1}(0)$. A Sotomayor-Teixeira regularization of $\mathbf{F}$, as described
at~\cite{Sotomayor1996}, is a 1-parameter family of smooth vector fields
$\mathbf{F}^{\varepsilon}$ that converges pointwisely to $\mathbf{F}$ as
$\varepsilon \to 0$. More precisely, for $\mathbf{x} \in U \setminus \Sigma$,
observe that the field $\mathbf{F}$ can be written in the form

\begin{equation}
  \label{sec:sotomayor_teixeira:eq:F_0}
  \mathbf{F}(\mathbf{x}) =
  \left[\frac{1 + \sgn(h(\mathbf{x}))}{2}\right] \mathbf{F}_+(\mathbf{x}) +
  \left[\frac{1 - \sgn(h(\mathbf{x}))}{2}\right] \mathbf{F}_-(\mathbf{x}),
\end{equation}

\noindent where $\sgn:\mathbb{R} \to \mathbb{R}$ is the \textbf{signal function}
given by

\begin{equation*}
  \sgn(x) = \begin{cases}
    -1, & \text{ if } x < 0, \\
    0,  & \text{ if } x = 0, \\
    1,  & \text{ if } x > 0,
  \end{cases}
\end{equation*}

\noindent which is a discontinuous function whose graph if represented at
\cref{sec:sotomayor_teixeira:subfig:signal}.

In order to approximate the piecewise smooth vector $\mathbf{F}$ with a
1-parameter family of smooth vector fields, we approximate the signal function
at \cref{sec:sotomayor_teixeira:eq:F_0} with a certain type of
smooth function. More precisely:

\begin{definition}%
  \label{sec:sotomayor_teixeira:defn:transition}
  We say that a smooth function $\varphi:\mathbb{R} \to \mathbb{R}$ is a
  \textbf{monotonous transition function\footnote{A study on the regularization
      process with non-monotonous transition functions can be found in chapter 6
      of \cite{Novaes2015}.}} if

  \begin{equation*}
    \varphi(x) = \begin{cases}
      -1, & \text{ if } x \le -1, \\
      1,  & \text{ if } x \ge 1,
    \end{cases}
  \end{equation*}

  \noindent and $\varphi'(x) > 0$ for $-1 < x < 1$.
\end{definition}

The graph of a typical transition function is represented at
\cref{sec:sotomayor_teixeira:subfig:transition}. Observe that, if
we define $\varphi^{\varepsilon}(x) =
  \varphi\left(\frac{x}{\varepsilon}\right)$, where $\varepsilon > 0$, then
clearly $\varphi^{\varepsilon} \to \sgn$ pointwisely when $\varepsilon \to 0$,
as long as their domains are restricted to the set $\mathbb{R}\setminus\{0\}$.
In particular, if we define

\begin{equation}
  \label{sec:sotomayor_teixeira:eq:F_epsilon}
  \mathbf{F}^{\varepsilon}(\mathbf{x}) =
  \left[\frac{1 + \varphi^{\varepsilon}(h(\mathbf{x}))}{2}\right]
  \mathbf{F}_+(\mathbf{x}) +
  \left[\frac{1 - \varphi^{\varepsilon}(h(\mathbf{x}))}{2}\right]
  \mathbf{F}_-(\mathbf{x}),
\end{equation}

\noindent then we get a 1-parameter family of vector fields
$\mathbf{F}^{\varepsilon} \in C^k(U)$ such that $\mathbf{F}^{\varepsilon} \to
  \mathbf{F}$ pointwisely when $\varepsilon \to 0$, as long as their domains are
restricted to the set $\mathbb{R}\setminus\{0\}$.

\begin{definition}%
  \label{sec:sotomayor_teixeira:defn:regularization} Let
  $\varphi:\mathbb{R} \to \mathbb{R}$ be a monotonous transition function.  We
  say that \cref{sec:sotomayor_teixeira:eq:F_epsilon} is a
  $\varphi^{\varepsilon}$-\textbf{regularization} of
  \cref{sec:sotomayor_teixeira:eq:F_0}.
\end{definition}

Observe that the regularization $\mathbf{F}^{\varepsilon}$ coincides with
$\mathbf{F}$ outside the rectangle given by $-\varepsilon < h(\mathbf{x}) <
  \varepsilon$. In fact,

\begin{equation*}
  \mathbf{F}^{\varepsilon}(\mathbf{x}) = \begin{cases}
    \mathbf{F}_+(\mathbf{x}), & \text{ if } h(\mathbf{x}) \ge \varepsilon,  \\
    \mathbf{F}_-(\mathbf{x}), & \text{ if } h(\mathbf{x}) \le -\varepsilon,
  \end{cases}
\end{equation*}

\noindent as represented at
\cref{sec:sotomayor_teixeira:subfig:regularization}. In
particular, it is clear that $\mathbf{F}^{\varepsilon}$ recovers the smooth
component of the Filippov dynamics given by $\mathbf{F}$, i.e., that associated
to the region $U\setminus\Sigma$, as long as we take $\varepsilon > 0$ small
enough.  As described in the next section, $\mathbf{F}^{\varepsilon}$ also
recovers the non-smooth component of the Filippov dynamics, i.e., that
associated to the region $\Sigma$.

\subsection{Geometrical Singular Perturbation Theory}
\label{sec:perturbacao_singular}

Let $W \subset \mathbb{R}^{m+n}$ be an open set whose elements are represented
by $(\mathbf{x},\mathbf{y})$. Let also $\mathbf{f}:W \times [0,1] \to
  \mathbb{R}^m$ and $\mathbf{g}:W \times [0,1] \to \mathbb{R}^n$ be vector fields
of class $C^k$ with $k \ge 1$. Given $0 < \xi < 1$, consider the system of
differential equations

\begin{equation}
  \label{sec:perturbacao_singular:eq:fast}
  \left\{\begin{matrix*}[l]
    \mathbf{x}' = \mathbf{f}(\mathbf{x}, \mathbf{y}, \xi)\\
    \mathbf{y}' = \xi \mathbf{g}(\mathbf{x}, \mathbf{y}, \xi)
  \end{matrix*}\right.,
\end{equation}

\noindent where $\square' = \sfrac{d\square}{d\tau}$, $\mathbf{x} =
  \mathbf{x}(\tau)$ and $\mathbf{y} = \mathbf{y}(\tau)$. Applying at the
previous system the time rescaling given by $t = \xi \tau$, we obtain the new
system

\begin{equation}
  \label{sec:perturbacao_singular:eq:slow}
  \left\{\begin{matrix*}[l]
    \xi \dot{\mathbf{x}} = \mathbf{f}(\mathbf{x}, \mathbf{y}, \xi)\\
    \dot{\mathbf{y}} = \mathbf{g}(\mathbf{x}, \mathbf{y}, \xi)
  \end{matrix*}\right.,
\end{equation}

\noindent where $\dot{\square} = \sfrac{d\square}{dt}$, $\mathbf{x} =
  \mathbf{x}(t)$ and $\mathbf{y} = \mathbf{y}(t)$.

As $0 < \xi < 1$, then \cref{sec:perturbacao_singular:eq:fast} and
\cref{sec:perturbacao_singular:eq:slow} have exactly the same phase portrait,
except for the trajectories speed, which is greater for first system and smaller
for the second. Therefore, the following definition makes sense:

\begin{definition}%
  \label{sec:perturbacao_singular:defn:slow-fast}
  We say that \cref{sec:perturbacao_singular:eq:fast} and
  \cref{sec:perturbacao_singular:eq:slow} form a $(m,n)$-\textbf{slow-fast
    system} with \textbf{fast system} given by
  \cref{sec:perturbacao_singular:eq:fast} and \textbf{slow system} given by
  \cref{sec:perturbacao_singular:eq:slow}.
\end{definition}

Taking $\xi \to 0$ in \cref{sec:perturbacao_singular:eq:fast}, we get the
so-called \textbf{layer system}

\begin{equation}
  \label{sec:perturbacao_singular:eq:layer}
  \left\{\begin{matrix*}[l]
    \mathbf{x}' = \mathbf{f}(\mathbf{x}, \mathbf{y}, 0)\\
    \mathbf{y}' = \mathbf{0}
  \end{matrix*}\right.,
\end{equation}

\noindent which has dimension $m$. Taking $\xi \to 0$ in
\cref{sec:perturbacao_singular:eq:slow}, we get the so-called \textbf{reduced
  system}

\begin{equation}
  \label{sec:perturbacao_singular:eq:reduced}
  \left\{\begin{matrix*}[l]
    \mathbf{0} = \mathbf{f}(\mathbf{x}, \mathbf{y}, 0)\\
    \dot{\mathbf{y}} = \mathbf{g}(\mathbf{x}, \mathbf{y}, 0)
  \end{matrix*}\right.,
\end{equation}

\noindent which has dimension $n$. Beyond that, we say that the set

\begin{equation*}
  \mathcal{M} = \left\{(\mathbf{x},\mathbf{y}) \in W;~ \mathbf{f}(\mathbf{x},
  \mathbf{y}, 0) = \mathbf{0}\right\}
\end{equation*}

\noindent is the \textbf{slow manifold}. Observe that, on one hand,
$\mathcal{M}$ represents the set of singularities of the layer system; on the
other hand, $\mathcal{M}$ represents the manifold over which the dynamics of the
reduced system takes place.

The main idea of Geometrical Singular Perturbation Theory, or GSP-Theory for
short, established by Fenichel in \cite{Fenichel1979}, consists of combining the
dynamics of the limit systems (layer and reduced) to recover the dynamics of the
initial system (slow-fast) with $\xi > 0$ small. In fact, considering $\xi$ as
an additional variable of the slow system
\cref{sec:perturbacao_singular:eq:slow} we get the new one

\begin{equation}
  \label{sec:perturbacao_singular:eq:fast_complete}
  \left\{\begin{matrix*}[l]
    \mathbf{x}' = \mathbf{f}(\mathbf{x}, \mathbf{y}, \xi)\\
    \mathbf{y}' = \xi \mathbf{g}(\mathbf{x}, \mathbf{y}, \xi)\\
    \xi' = 0
  \end{matrix*}\right.,
\end{equation}

\noindent whose Jacobian matrix at $\left( \mathbf{x}_0, \mathbf{y}_0, 0 \right)
  \in \mathcal{M} \times \left\{ 0 \right\}$ is

\begin{equation}
  \label{sec:perturbacao_singular:eq:jacobian_layer}
  \mathbf{J}_{\text{fast}} =
  \begin{bmatrix}
    \mathbf{f}_{\mathbf{x}} & \mathbf{f}_{\mathbf{y}} & 0 \\
    \mathbf{0}              & \mathbf{0}              & 0 \\
    \mathbf{0}              & \mathbf{0}              & 0
  \end{bmatrix},
\end{equation}

\noindent where $\mathbf{f}_{\mathbf{x}}$ and $\mathbf{f}_{\mathbf{y}}$
represent the partial derivatives calculated at the point $\left( \mathbf{x}_0,
  \mathbf{y}_0, 0 \right)$. The matrix above has the trivial eigenvalue $\lambda =
  0$ with algebraic multiplicity $n+1$. The remaining eigenvalues, called
\textbf{non-trivial}, are divided in three categories: negative, zero or
positive real parts; we denote the number of such eigenvalues by $k^s$, $k^c$
and $k^u$, respectively.

\begin{definition}
  \label{sec:perturbacao_singular:defn:normal_hyperbolic}
  We say that $\left( \mathbf{x}_0, \mathbf{y}_0, 0 \right) \in \mathcal{M}
    \times \left\{ 0 \right\}$ is \textbf{normally hyperbolic} if every
  non-trivial eigenvalue of \cref{sec:perturbacao_singular:eq:jacobian_layer}
  have non-zero real part, i.e., $k^c = 0$.
\end{definition}

Fenichel, in \cite{Fenichel1979}, proved that normal hyperbolicity allows the
persistence of invariant compact parts of the slow manifold under singular
perturbation, i.e., the dynamical structure of such parts with $\xi = 0$
persists for $\xi > 0$ small. Even more, with predictable stability. More
precisely:

\begin{theorem}[Retrieved from \cite{Teixeira2012}, page 1953]
  \label{sec:perturbacao_singular:thm:fenichel}
  Let $\mathcal{N}$ be a normally hyperbolic compact invariant $j$-dimensional
  submanifold of $\mathcal{M}$. Suppose that the stable and unstable manifolds
  of $\mathcal{N}$, with respect to the reduced system, have dimensions $j+j^s$
  and $j+j^u$, respectively. Then, there exists a 1-parameter family of
  invariant submanifolds $\left\{ \mathcal{N}_{\xi} ;~ \xi \sim 0 \right\}$ such
  that $\mathcal{N}_{0} = \mathcal{N}$ and $\mathcal{N}_{\xi}$ has stable and
  unstable manifolds with dimensions $j+j^s+k^s$ and $j+j^u+k^u$, respectively.
\end{theorem}

The reverse idea of GSP-Theory can also be used to recover the non-smooth
component of the Filippov dynamics, given by the piecewise vector field
($\varepsilon = 0$), from its regularization ($\varepsilon > 0$). In fact, let
$\mathbf{F} = (\mathbf{F}_+,\mathbf{F}_-) \in \mathcal{R}^k(U,h)$ be a piecewise
smooth vector field with switching manifold $\Sigma = h^{-1}(\{0\})$. Let also
$\varphi:\mathbb{R} \to \mathbb{R}$ be a monotonous transition function and
$\mathbf{F}^{\varepsilon}$ the $\varphi^{\varepsilon}$-regularization of
$\mathbf{F}$.

We need to transform $\mathbf{F}^{\varepsilon}$ in a slow-fast system. In order
to do so, observe that, as $0$ is a regular value of $h$, then from the Local
Normal Form for Submersions follows that, without loss of generality, we can
admit that $h(x_1,\ldots,x_n) = x_1$ in a neighborhood of a given point
$\mathbf{x} \in \Sigma$. Therefore, if we write $\mathbf{F}_+ =
  (f_1^+,\ldots,f_n^+)$ and $\mathbf{F}_- = (f_1^-,\ldots,f_n^-)$, then follows
that $\mathbf{F}^{\varepsilon}$ can be written as

\begin{align*}
  \dot{x}_i = &
  \left[\frac{1+\varphi^{\varepsilon}(x_1)}{2}\right]f_i^+(x_1,\ldots,x_n) +              \\
  +           & \left[\frac{1+\varphi^{\varepsilon}(x_1)}{2}\right]f_i^-(x_1,\ldots,x_n),
\end{align*}

\noindent where $i \in \{1, \ldots, n\}$. Now, applying to the system above the
polar blow-up given by $x_1 = \xi \cos{\theta}$ and $\varepsilon = \xi
  \sin{\theta}$, where $\xi \ge 0$ and $\theta \in [0,\pi]$, we obtain a
$(1,n-1)$-slow-fast system given by

\begin{equation}
  \label{sec:perturbacao_singular:eq:regul_slow-fast}
  \left\{\begin{matrix*}[l]
    \xi \dot{\theta} = \alpha_1(\theta,x_2,\ldots,x_n,\xi)\\
    \dot{x}_i = \alpha_i(\theta,x_2,\ldots,x_n,\xi)
  \end{matrix*}\right.,
\end{equation}

\noindent where $i \in \{2,\ldots,n \}$.

Observe that, for $\xi = 0$, we have $x_1 = 0$ and $\varepsilon = 0$, i.e., we
are at the non-regularized system $\mathbf{F}$ over the manifold $\Sigma$.  In
the other hand, for $\xi > 0$ and $\theta \in (0,\pi)$, we have $-\xi < x_1 <
  \xi$ and $0 < \varepsilon < \xi$, i.e., we are at the regularized system
$\mathbf{F}^{\varepsilon}$ over the rectangle where it does not coincide to
$\mathbf{F}$, see \cref{sec:sotomayor_teixeira:subfig:regularization}. The
authors of~\cite{Teixeira2012} then proved the result below:

\begin{theorem}[Retrieved from \cite{Teixeira2012}, page 1950]%
  \label{sec:perturbacao_singular:theorem:regular}
  Consider the piecewise smooth vector field $\mathbf{F}$ and the slow-fast
  system \cref{sec:perturbacao_singular:eq:regul_slow-fast}. The sliding region
  $\Sigma^{s}$ is homeomorphic to the slow manifold given by

  \begin{equation*}
    \alpha_1(\theta,x_2,\ldots,x_n,0) = 0
  \end{equation*}

  \noindent and the dynamics of the sliding vector field $\mathbf{F}^{s}$ over
  $\Sigma^{s}$ is topologically equivalent to that of the reduced system given by

  \begin{equation*}
    \left\{\begin{matrix*}[l]
      0 = \alpha_1(\theta,x_2,\ldots,x_n,0)\\
      \dot{x}_i = \alpha_i(\theta,x_2,\ldots,x_n,0)
    \end{matrix*}\right.,
  \end{equation*}

  \noindent where $i \in \{2,\ldots,n \}$.
\end{theorem}

Concisely, the Filippov dynamics of $\mathbf{F}$ is completely recovered by its
regularization $\mathbf{F}^{\varepsilon}$. In order to do so, the following
steps, described in details above, are necessary:

\begin{enumerate}
  \item Normalization of the switching manifold.
  \item Regularization of the piecewise smooth vector field.
  \item Polar blow-up of the regularization.
  \item Analysis of the resulting limit systems (layer and reduced).
\end{enumerate}

\section{Statement of the Problem}%
\label{sec:problem}

One of the fundamental hypotheses in the theory described in
\cref{sec:introduction} is the fact that $0 \in \mathbb{R}$ is a \textbf{regular
  value} of the function $h:\mathbb{R} \to \mathbb{R}$ and, therefore, the
switching manifold $\Sigma = h^{-1}(\{0\})$ is a regular surface. In that case,
as we have seen, there exists at least one well-defined and established dynamics
associated: the Filippov dynamics. A natural question to ask then is: can a
Filippov-like dynamics be defined for the case when $0 \in \mathbb{R}$ is a
\textbf{singular value} of the function $h:\mathbb{R} \to \mathbb{R}$, i.e.,
when the switching manifold is not a regular surface?

\begin{figure}[ht]
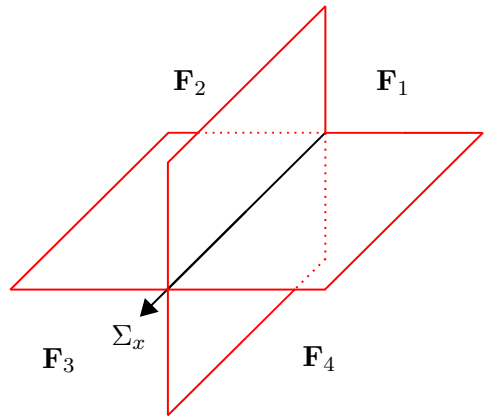

  \centering
  \def\svgwidth{0.6\linewidth}
  
  \caption{Double discontinuity.}%
  \label{sec:problem:fig:switching_manifold}
\end{figure}

In the next sections, we would like to study the particular case known as the
\textbf{double discontinuity}. This particular configuration of the switching
manifold is the simplest one between the four singular configurations (known as
Gutierrez-Sotomayor or simple manifolds) that, according to
\cite{Gutierrez1982}, breaks the regularity condition in a dynamically stable
manner. The double discontinuity is described in detail below.

Let $\mathbf{F}_i: \mathbb{R}^3 \to \mathbb{R}^3$ be vector fields of class
$C^k(\mathbb{R}^3)$ with $i \in \{1,2,3,4 \}$. The piecewise smooth vector field
$\mathbf{F}: \mathbb{R}^3 \to \mathbb{R}^3$ given by

\begin{equation}
  \label{sec:problem:eq:piecewise_field}
  \mathbf{F}(x,y,z) = \begin{cases}
    \mathbf{F}_1(x,y,z),~ & \text{if}~y \ge 0~\text{and}~z \ge 0, \\
    \mathbf{F}_2(x,y,z),~ & \text{if}~y \le 0~\text{and}~z \ge 0, \\
    \mathbf{F}_3(x,y,z),~ & \text{if}~y \le 0~\text{and}~z \le 0, \\
    \mathbf{F}_4(x,y,z),~ & \text{if}~y \ge 0~\text{and}~z \le 0,
  \end{cases}
\end{equation}

\noindent and denoted by $\mathbf{F} = \left(\mathbf{F}_1, \mathbf{F}_2,
  \mathbf{F}_3, \mathbf{F}_4\right)$ is said to have a \textbf{double
  discontinuity} as switching manifold, see
\cref{sec:problem:fig:switching_manifold}. The set of all vector fields
$\mathbf{F}$ defined as above will be denoted by

\begin{equation*}
  \mathcal{D}_3^k \equiv
  C^k(\mathbb{R}^3) \times C^k(\mathbb{R}^3) \times
  C^k(\mathbb{R}^3) \times C^k(\mathbb{R}^3)
\end{equation*}

\noindent and equipped with the Whitney product topology.

The double discontinuity, as defined above, consists of the planes $xy$ and
$xz$ perpendicularly intersecting at the $x$-axis, $\Sigma_x = \{(x,0,0) ;~ x
  \in \mathbb{R}\}$. For points in $\Sigma \setminus \Sigma_{x}$, the ordinary
Filippov dynamics described in \cref{sec:introduction} can be locally
applied. However, for points $(x,0,0) \in \Sigma_x$ that theory cannot be
directly applied.  In fact, $\Sigma = h^{-1}(\{0\})$, where $h:\mathbb{R}^3 \to
  \mathbb{R}$ given by $h(x,y,z) = yz$ has $0 \in \mathbb{R}$ as a singular value,
since $Dh(x,0,0)$ is not a surjective map for $(x,0,0) \in \Sigma_{x}$.

Therefore, we state the problem: given $\mathbf{F} \in \mathcal{D}_3^k$, can we
define a Filippov-like dynamics over $\Sigma_{x}$? How does it generally behave
there? In the next section, we present a methodology based on
\cite{Buzzi2012,Llibre2015,Panazzolo2017,Teixeira2012} to approach this
problem.

\section{Methodology}%
\label{sec:framework}

\begin{figure}[ht]
  \centering
  \begin{subfigure}[b]{0.45\linewidth}
    \def\svgwidth{\linewidth}
    
    \caption{Slice.}%
    \label{sec:framework:sub:blowup:subfig:slice}
  \end{subfigure}
  \qquad
  \begin{subfigure}[b]{0.45\linewidth}
    \def\svgwidth{\linewidth}
\begingroup%
  \makeatletter%
  \providecommand\color[2][]{%
    \errmessage{(Inkscape) Color is used for the text in Inkscape, but the package 'color.sty' is not loaded}%
    \renewcommand\color[2][]{}%
  }%
  \providecommand\transparent[1]{%
    \errmessage{(Inkscape) Transparency is used (non-zero) for the text in Inkscape, but the package 'transparent.sty' is not loaded}%
    \renewcommand\transparent[1]{}%
  }%
  \providecommand\rotatebox[2]{#2}%
  \newcommand*\fsize{\dimexpr\f@size pt\relax}%
  \newcommand*\lineheight[1]{\fontsize{\fsize}{#1\fsize}\selectfont}%
  \ifx\svgwidth\undefined%
    \setlength{\unitlength}{234.06073623bp}%
    \ifx\svgscale\undefined%
      \relax%
    \else%
      \setlength{\unitlength}{\unitlength * \real{\svgscale}}%
    \fi%
  \else%
    \setlength{\unitlength}{\svgwidth}%
  \fi%
  \global\let\svgwidth\undefined%
  \global\let\svgscale\undefined%
  \makeatother%
  \begin{picture}(1,1)%
    \lineheight{1}%
    \setlength\tabcolsep{0pt}%
    \put(0,0){\includegraphics[width=\unitlength,page=1]{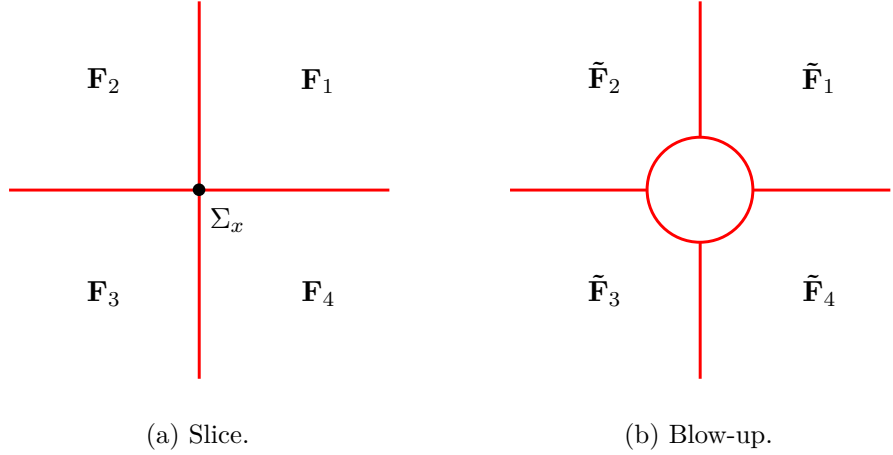}}%
    \put(0.73501035,0.73549126){\color[rgb]{0,0,0}\makebox(0,0)[lt]{\lineheight{0}\smash{\begin{tabular}[t]{l}$\mathbf{\tilde{F}}_1$\end{tabular}}}}%
    \put(0.2403361,0.24061482){\color[rgb]{0,0,0}\makebox(0,0)[lt]{\lineheight{0}\smash{\begin{tabular}[t]{l}$\mathbf{\tilde{F}}_3$\end{tabular}}}}%
    \put(0.73646542,0.24075031){\color[rgb]{0,0,0}\makebox(0,0)[lt]{\lineheight{0}\smash{\begin{tabular}[t]{l}$\mathbf{\tilde{F}}_4$\end{tabular}}}}%
    \put(0.2403361,0.73681095){\color[rgb]{0,0,0}\makebox(0,0)[lt]{\lineheight{0}\smash{\begin{tabular}[t]{l}$\mathbf{\tilde{F}}_2$\end{tabular}}}}%
  \end{picture}%
\endgroup%

    \caption{Blow-up.}%
    \label{sec:framework:sub:blowup:subfig:slice_blowup}
  \end{subfigure}
  \caption{Framework process at slice-level.}%
  \label{sec:framework:sub:blowup:fig:framework_process}
\end{figure}

The first step consists of the application of a polar blow-up at the origin of
the slice represented at \cref{sec:framework:sub:blowup:subfig:slice} or, in
other words, a \textbf{cylindrical blow-up} at $\Sigma_{x}$. More specifically,
assuming that the components of $\mathbf{F} \in \mathcal{D}_3^k$ can be written as

\begin{equation*}
  \mathbf{F}_i = (w_i, p_i, q_i),
\end{equation*}

\noindent we apply the blow-up $\phi_1: \mathbb{R} \times S^1 \times \mathbb{R}^+ \to \mathbb{R}^3$ given by

\begin{equation*}
  \phi_1(x,\theta,r) = (x, r\cos{\theta}, r\sin{\theta}),
\end{equation*}

\noindent which induces $\mathbf{\tilde{F}} = [(\phi_1)_*^{-1}\mathbf{F}] \circ
  \phi_1$ whose components are given by

\begin{equation*}
  \mathbf{\tilde{F}}_i =
  \left(w_i,
  \frac{q_i \cos{\theta} - p_i \sin{\theta}}{r},
  p_i \cos{\theta} + q_i \sin{\theta} \right),
\end{equation*}

\noindent where $w_i$, $p_i$ and $q_i$ must be evaluated at the point
$\phi_1(x,\theta,r)$. We then define the set

\begin{equation*}
  \tilde{\mathcal{D}}_3^k =
  \left\{\mathbf{\tilde{F}} = [(\phi_1)_*^{-1}\mathbf{F}] \circ \phi_1 ;~
  \mathbf{F} \in \mathcal{D}_3^k
  \right\}
\end{equation*}

\noindent of all blow-up induced vector fields.

An extremely important observation at this point consists in the theorem stated
below with minor adaptations to our notation relative to the original one found
in \cite{Llibre2015}\footnote{Up to our knowledge, this theorem where actually
  first stated in \cite[p.~449]{Buzzi2012}. However, \cite{Llibre2015} also
  provides analogous results for the triple, cone, and Whitney discontinuities.
  See \cref{sec:problem:fig:manifolds}. Regarding the double discontinuity,
  similar versions of the theorem can also be found in \cite{Panazzolo2017},
  within the context of foliations, and in \cite{Teixeira2012}, which is actually
  a survey. Besides that, as raised in \cref{sec:introduction}, this theorem (with
  some sparse examples) represents the state of the art. \textbf{In other words,
    results following the statement are novelties on a sparsely explored
    territory.}}

\begin{rtheorem}[Retrieved from \cite{Llibre2015}, page 498]%
  \label{sec:framework:sub:blowup:thm:regular_discontinuities}
  The map \\ $\phi_1: \mathbb{R} \times S^1 \times \mathbb{R}^+ \to \mathbb{R}^3$
  given by \[\phi_1(x,\theta,r) = (x, r\cos{\theta}, r\sin{\theta})\] induces a
  vector field $\mathbf{\tilde{F}}$ satisfying that any discontinuity $q \in
    \tilde{\Sigma} = \phi_1^{-1}\left( \Sigma \right)$ is regular.
\end{rtheorem}

\noindent Hence, since the induced vector field $\mathbf{\tilde{F}}$ has only
\textbf{regular discontinuities}, then classical Filippov theory, as presented
at \cref{sec:introduction}, is enough for its analysis. More precisely, we have
now a piecewise smooth vector field $\mathbf{\tilde{F}}$ given by the four
smooth vector fields $\mathbf{\tilde{F}}_i$, which induces the four
\textbf{slow-fast systems}

\begin{equation}
  \label{sec:framework:sub:blowup:eq:slow-fast}
  \left\{\begin{matrix*}[l]
    \dot{x} = w_i\\
    r \dot{\theta} = q_i \cos{\theta} - p_i \sin{\theta}\\
    \dot{r} = p_i \cos{\theta} + q_i \sin{\theta}
  \end{matrix*}\right.,
\end{equation}

\noindent where $\dot{\square} = \sfrac{d\square}/{dt}$; $w_i$, $p_i$ and $q_i$
must be calculated at the point $\phi_1(x,\theta,r)$; and $r$ is the time
rescaling factor.

The study of the dynamics of~\eqref{sec:problem:eq:piecewise_field} has
therefore been reduced to the study of the slow-fast systems
\eqref{sec:framework:sub:blowup:eq:slow-fast}. In particular, the dynamics over
$\Sigma_x$, previously undefined, can now be associated with
\cref{sec:framework:sub:blowup:eq:slow-fast} at $r = 0$, which is given by the
combination of the dynamics of the \textbf{reduced system}

\begin{equation}
  \label{sec:framework:blowup:eq:slow_system}
  \left\{\begin{matrix*}[l]
    \dot{x} = w_i\\
    0 = q_i \cos{\theta} - p_i \sin{\theta}\\
    \dot{r} = p_i \cos{\theta} + q_i \sin{\theta}
  \end{matrix*}\right.
\end{equation}

\noindent and the dynamics of the \textbf{layer system}

\begin{equation}
  \label{sec:framework:blowup:eq:fast_system}
  \left\{\begin{matrix*}[l]
    x' = 0\\
    \theta' = q_i \cos{\theta} - p_i \sin{\theta}\\
    r' = 0
  \end{matrix*}\right.,
\end{equation}

\noindent where $\square' = \sfrac{d\square}/{d\tau}$ with $t = r \tau$; and the
components $w_i$, $p_i$ and $q_i$ must be calculated at the point
$\phi_1(x,\theta,0) = (x,0,0)$.

\begin{figure}[ht]
  \centering
  \def\svgwidth{0.95\linewidth}
  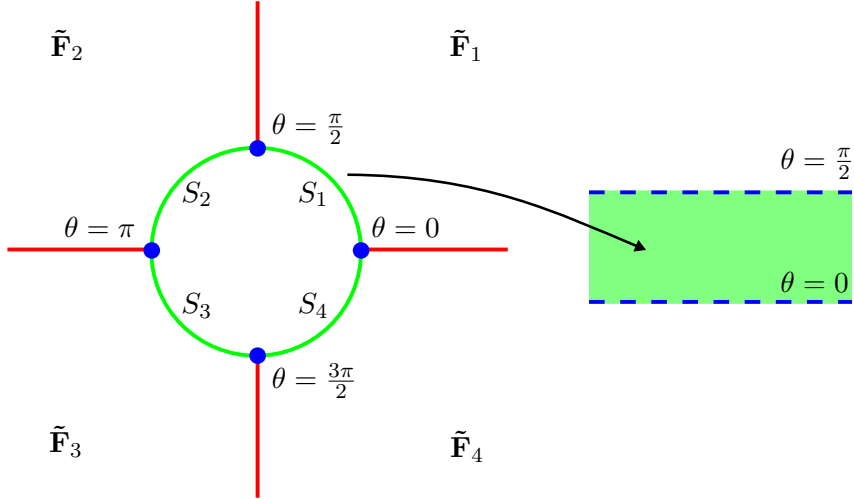
  \caption{Green cylinder $C$ divided in the four stripes $S_i$. A scheme of
    the stripe $S_1$ is also put in evidence.}
  \label{sec:framework:sec:blowup:fig:stripes}
\end{figure}

Geometrically, the dynamics over $\Sigma_{x}$ in
\cref{sec:problem:eq:piecewise_field} can now be associated to the dynamics over
the cylinder $C = \mathbb{R} \times S^1$ divided in the four infinite stripes

\begin{align*}
   & S_2 = \mathbb{R} \times [\sfrac{\pi}{2}, \pi],
   & S_1                                             & = \mathbb{R} \times [0,\sfrac{\pi}{2}],      \\
   & S_3 = \mathbb{R} \times [\pi, \sfrac{3\pi}{2}],
   & S_4                                             & = \mathbb{R} \times [\sfrac{3\pi}{2}, 2\pi],
\end{align*}

\noindent as represented at \cref{sec:framework:sec:blowup:fig:stripes}, where
the slow-fast systems given by
\cref{sec:framework:blowup:eq:slow_system,sec:framework:blowup:eq:fast_system}
acts, respectively. As we previously stated at
\cref{sec:framework:sub:blowup:thm:regular_discontinuities}, the four lines
where these stripes intersect admits at most regular discontinuities. Finally,
the analysis of the dynamics on each stripe $S_i$ can then be carried out using
GSP-Theory.

In particular, the first two equations of the system
\eqref{sec:framework:blowup:eq:slow_system} are independent of $r$ and,
therefore, it can be decoupled as

\begin{equation}
  \label{sec:framework:blowup:eq:slow_system_cyl}
  \left\{\begin{matrix*}[l]
    \dot{x} = w_i\\
    0 = q_i \cos{\theta} - p_i \sin{\theta}\\
  \end{matrix*}\right.,
\end{equation}

\noindent which gives the \textbf{reduced dynamics} over $S_i$; and

\begin{equation}
  \label{sec:framework:blowup:eq:slow_system_rad}
  \dot{r} = p_i \cos{\theta} + q_i \sin{\theta},
\end{equation}

\noindent which gives the respective \textbf{slow radial dynamics} or, in other
words, it indicates how the external dynamics communicates with the dynamics
\eqref{sec:framework:blowup:eq:slow_system_cyl} over the cylinder: entering
($\dot{r} > 0$), leaving ($\dot{r} < 0$) or staying ($\dot{r} = 0$) at $S_i$.

Analogously, the first two equations of the system
\eqref{sec:framework:blowup:eq:fast_system} are independent of $r$ and,
therefore, it can also be decoupled as

\begin{equation}
  \label{sec:framework:blowup:eq:fast_system_cyl}
  \left\{\begin{matrix*}[l]
    x' = 0\\
    \theta' = q_i \cos{\theta} - p_i \sin{\theta}
  \end{matrix*}\right.,
\end{equation}

\noindent which gives the \textbf{layer dynamics} over $S_i$; and

\begin{equation}
  \label{sec:framework:blowup:eq:fast_system_rad}
  r' = 0,
\end{equation}

\noindent which gives the respective \textbf{fast radial dynamics} over the
cylinder.

Summarizing, we conclude that the dynamics over $\Sigma_x$ behaves as described
in the \textbf{fundamental lemma} below, whose proof consists in the analysis
done above.

\LemmaFundamental

In order to perform a deeper analysis of the dynamics given by
\cref{sec:framework:blowup:thm:dynamics} with GSP-Theory as described at
\cref{sec:perturbacao_singular}, let $S_i$ be one of the cylinder's stripe and
let

\begin{equation*}
  \mathcal{M}_i = \left\{ (x, \theta) \in S_i \subset \mathbb{R} \times S^1 ;~
  f_i(x, \theta, 0) = 0 \right\}
\end{equation*}

\noindent be its slow manifold, where $f_i(x, \theta, 0) = q_i \cos{\theta} -
  p_i \sin{\theta}$.

Given $(x_0, \theta_0, 0) \in \mathcal{M}_i \times \left\{ 0 \right\}$, the
Jacobian matrix of the complete layer system
\cref{sec:framework:blowup:eq:fast_system} over this point is

\begin{equation*}
  \mathbf{J}_{\text{fast}} =
  \begin{bmatrix}
    0       & 0              & 0 \\
    (f_i)_x & (f_i)_{\theta} & 0 \\
    0       & 0              & 0
  \end{bmatrix},
\end{equation*}

\noindent where $(f_i)_x$ and $(f_i)_{\theta}$ represents the partial
derivatives calculated at $\left( x_0, \theta_0, 0 \right)$. The eigenvalues of
this matrix are the elements of the set $\left\{0, 0, (f_i)_{\theta} \right\}$
and, therefore, $(x_0, \theta_0)$ is normally hyperbolic if, and only if,
$(f_i)_{\theta} \neq 0$. However, we observe that, since we are over the slow
manifold, then $(f_i)_{\theta} = 0$ leads to the homogeneous linear system

\begin{align*}
  \left\{\begin{matrix*}[r]
    f_i = 0 \\
    (f_i)_{\theta} = 0
  \end{matrix*}\right.
  \quad & \sim \quad
  \left\{\begin{matrix}
    q_i \cos{\theta} - p_i \sin{\theta} = 0 \\
    q_i \sin{\theta} + p_i \cos{\theta} = 0
  \end{matrix}\right.
  \quad \sim \quad   \\
        & \sim \quad
  \begin{bmatrix}
    \cos{\theta} & -\sin{\theta} \\
    \sin{\theta} & \cos{\theta}
  \end{bmatrix}
  \begin{bmatrix}
    q_i \\
    p_i
  \end{bmatrix}
  =
  \begin{bmatrix}
    0 \\
    0
  \end{bmatrix}
\end{align*}

\noindent whose unique solution is the trivial, $p_i = q_i = 0$, since the
trigonometrical matrix above is invertible ($\text{det} \equiv 1$) for every
$\theta \in S^1$ and, therefore, we conclude that $(f_i)_{\theta} \neq 0$
whenever

\begin{equation}
  \tag{WFH}
  \label{sec:framework:sub:blowup:eq:weak_fund_hypo}
  p_i \neq 0 \quad \text{or} \quad q_i \neq 0,
\end{equation}

\noindent henceforth, called \textbf{weak fundamental hypothesis}, or WFH for
short. We also observe that

\begin{equation*}
  (f_i)_x = (q_i)_x \cos{\theta} - (p_i)_x \sin{\theta}
\end{equation*}

\noindent which, as above, supposing $(f_i)_x = 0$ leads to the homogeneous
linear system

\begin{align*}
  \left\{\begin{matrix*}[r]
    q_i \cos{\theta} - p_i \sin{\theta} = 0 \\
    (q_i)_x \cos{\theta} - (p_i)_x \sin{\theta} = 0
  \end{matrix*}\right.
  \quad \sim \quad \\
  \quad \sim \quad
  \begin{bmatrix*}[c]
    q_i     & p_i     \\
    (q_i)_x & (p_i)_x
  \end{bmatrix*}
  \begin{bmatrix}
    \cos{\theta} \\
    \sin{\theta}
  \end{bmatrix}
  =
  \begin{bmatrix}
    0 \\
    0
  \end{bmatrix}
\end{align*}

\noindent which only admits the absurd solution $\cos{\theta} = \sin{\theta} =
  0$ if the matrix above is invertible. Hence, we can ensure $(f_i)_x \neq 0$ by
imposing this absurd, i.e.,

\begin{equation}
  \tag{SFH}
  \label{sec:framework:sub:blowup:eq:strong_fund_hypo}
  0 \neq
  \det{\begin{bmatrix*}[c]
      q_i     & p_i     \\
      (q_i)_x & (p_i)_x
    \end{bmatrix*}}
  =
  q_i (p_i)_x - p_i (q_i)_x
\end{equation}

\noindent which always implies the weak fundamental hypothesis and, therefore,
will be called \textbf{strong fundamental hypothesis}, or SFH for short.

\TheoremTransversality

\begin{proof}
  The first part of the statement is assured by
  \cref{sec:framework:blowup:thm:dynamics}. For the second part, just observe
  that $\dot{r} = - (f_i)_{\theta} \neq 0$ under
  \cref{sec:framework:sub:blowup:eq:weak_fund_hypo}.
\end{proof}

\TheoremGraph

\begin{proof}
  The first part is assured by the usual Implicit Function Theorem applied to
  $f_i(x_0, \theta_0, 0) = 0$ over $\mathcal{M}_i$, since under
  \cref{sec:framework:sub:blowup:eq:weak_fund_hypo} we have
  $\norm{(f_i)_{\theta}} > 0$. Analogously, the second part is assured by the
  Global Implicit Function Theorem found in \cite[p.~253]{Zhang2006}, which
  requires a stronger hypothesis.
\end{proof}

\TheoremNormalHiperbolicity

\begin{proof}
  Just observe that $(f_i)_{\theta}$, the only non-trivial eigenvalue, is
  non-zero under \cref{sec:framework:sub:blowup:eq:weak_fund_hypo}.
\end{proof}

\TheoremHiperbolicSingularities

\begin{proof}
  Let $\mathbf{P} = \left( x_0, \theta_0 \right) \in \mathcal{M}_i$ be a
  hyperbolic singularity of the reduced system, i.e., $w_i(x_0, 0, 0) = 0$ with
  eigenvalue $\lambda_1 = (w_i)_{x}(x_0, 0, 0) \neq 0$. We have two
  possibilities:

  \begin{itemize}
    \item $\lambda_1 > 0 \Rightarrow (j^s, j^u) = (0,1)$; or
    \item $\lambda_1 < 0 \Rightarrow (j^s, j^u) = (1,0)$,
  \end{itemize}

  \noindent where $j^s$ and $j^u$ are the dimensions of the stable and unstable
  manifolds of $\mathbf{P}$ with respect to the reduced system, respectively.

  On the other hand, under \cref{sec:framework:sub:blowup:eq:weak_fund_hypo} we
  also have the non-trivial eigenvalue $\lambda_2 = (f_i)_{\theta}(x_0,
    \theta_0, 0) \neq 0$ for the layer system and, therefore, the two
  possibilities:

  \begin{itemize}
    \item $\lambda_2 > 0 \Rightarrow (k^s, k^u) = (0,1)$; or
    \item $\lambda_2 < 0 \Rightarrow (k^s, k^u) = (1,0)$,
  \end{itemize}

  \noindent where $k^s$ and $k^u$ are the dimensions of the stable and unstable
  manifolds of $\mathbf{P}$ with respect to the layer system, respectively.

  Hence, observing that $j = \dim{\mathbf{P}} = 0$ and remembering
  \cref{sec:perturbacao_singular:thm:fenichel}, any combination of the signs of
  $\lambda_1$ and $\lambda_2$ leads to the total sum of dimensions

  \begin{equation*}
    (j^s + k^s) + (j^u + k^u) = 2 = \dim{S_i},
  \end{equation*}

  \noindent and, therefore, $\mathbf{P}$ acts as a hyperbolic singularity of
  $S_i$. Finally, the saddle-node duality comes from the fact that both
  non-trivial eigenvalues above have no imaginary parts.
\end{proof}

In other words, under \cref{sec:framework:sub:blowup:eq:weak_fund_hypo}, the
slow manifold $\mathcal{M}_i$ is, at the very least, locally a graph. More than
that, it is the entry-point for the external dynamics to the cylinder. Besides
that, it is normally hyperbolic at its full extension, assuring then not only
persistence and well-behaved stability for its invariant compact parts, but also
that $\mathcal{M}_i$ is always attracting or repelling the surrounding (layer)
dynamics. All these nice properties come at the low cost of
\cref{sec:framework:sub:blowup:eq:weak_fund_hypo}. Therefore, it is not a
surprise that, for every system studied below, we require at least
\cref{sec:framework:sub:blowup:eq:weak_fund_hypo}, but also always test for
\cref{sec:framework:sub:blowup:eq:strong_fund_hypo}, whose importance will
become clear when studying affine systems.

\section{Constant Dynamics}%
\label{sec:constant}

Let $\mathcal{C}_3 \subset \mathcal{D}_3^k$ be the set of all piecewise smooth
vector fields $\mathbf{F}$ with a double discontinuity given by constant vector
fields

\begin{equation}
  \label{sec:constant:eq:system}
  \mathbf{F}_i(x,y,z) = (d_{i1}, d_{i2}, d_{i3}),
\end{equation}

\noindent where $d_{ij} \in \mathbb{R}$ for all $i$ and $j$. According to the
Fundamental \cref{sec:framework:blowup:thm:dynamics}, the dynamics over
$\Sigma_x$ of such a field is blow-up associated to the following fundamental
dynamics over the cylinder $C = \mathbb{R} \times S^1 = S_1 \cup \ldots \cup
  S_4$: over each stripe $S_i$ acts a slow-fast dynamics whose reduced dynamics is
given by

\begin{equation}
  \label{sec:constant:eq:reduced}
  \left\{\begin{matrix*}[l]
    \dot{x} = d_{i1}\\
    0 = d_{i3} \cos{\theta} - d_{i2} \sin{\theta}\\
  \end{matrix*}\right.,
\end{equation}

\noindent with radial slow dynamics $\dot{r} = d_{i2} \cos{\theta} + d_{i3}
  \sin{\theta}$; and layer dynamics given by

\begin{equation}
  \label{sec:constant:eq:layer}
  \left\{\begin{matrix*}[l]
    x' = 0\\
    \theta' = d_{i3} \cos{\theta} - d_{i2} \sin{\theta}
  \end{matrix*}\right.,
\end{equation}

\noindent with radial fast dynamics $r' = 0$.

Besides that, for \cref{sec:constant:eq:system}, we have $p_i = d_{i2}$ and $q_i
  = d_{i3}$ so that \cref{sec:framework:sub:blowup:eq:weak_fund_hypo} is
satisfied as long as

\begin{equation}
  \label{sec:constant:eq:wfh}
  d_{i2} \neq 0 \quad \text{or} \quad d_{i3} \neq 0,
\end{equation}

\noindent whereas \cref{sec:framework:sub:blowup:eq:strong_fund_hypo} is
\textbf{never} satisfied, since $(p_i)_x = (q_i)_x = 0$.

Therefore, our goal at this section is to fully describe the fundamental
dynamics of \cref{sec:constant:eq:system} over the cylinder $C$ under the
hypothesis \cref{sec:constant:eq:wfh}. In order to do so, we are going to
systematically analyze the slow-fast systems
\eqref{sec:constant:eq:reduced}--\eqref{sec:constant:eq:layer} for the cases
suggested by \cref{sec:constant:eq:wfh}. This analysis takes place in
\cref{sec:constant:sub:d2_nonzero,sec:constant:sub:d2_zero}, resulting in
\cref{sec:constant:thm:dynamics} stated and exemplified at
\cref{sec:constant:sub:theorem}.

\subsection{Case \texorpdfstring{$d_{i2} \neq 0$}{d2=/=0}}
\label{sec:constant:sub:d2_nonzero}

\begin{figure}[ht]
  \centering
  \def\svgwidth{0.85\linewidth}
\begingroup%
  \makeatletter%
  \providecommand\color[2][]{%
    \errmessage{(Inkscape) Color is used for the text in Inkscape, but the package 'color.sty' is not loaded}%
    \renewcommand\color[2][]{}%
  }%
  \providecommand\transparent[1]{%
    \errmessage{(Inkscape) Transparency is used (non-zero) for the text in Inkscape, but the package 'transparent.sty' is not loaded}%
    \renewcommand\transparent[1]{}%
  }%
  \providecommand\rotatebox[2]{#2}%
  \newcommand*\fsize{\dimexpr\f@size pt\relax}%
  \newcommand*\lineheight[1]{\fontsize{\fsize}{#1\fsize}\selectfont}%
  \ifx\svgwidth\undefined%
    \setlength{\unitlength}{518.12127686bp}%
    \ifx\svgscale\undefined%
      \relax%
    \else%
      \setlength{\unitlength}{\unitlength * \real{\svgscale}}%
    \fi%
  \else%
    \setlength{\unitlength}{\svgwidth}%
  \fi%
  \global\let\svgwidth\undefined%
  \global\let\svgscale\undefined%
  \makeatother%
  \begin{picture}(1,0.62823176)%
    \lineheight{1}%
    \setlength\tabcolsep{0pt}%
    \put(0,0){\includegraphics[width=\unitlength,page=1]{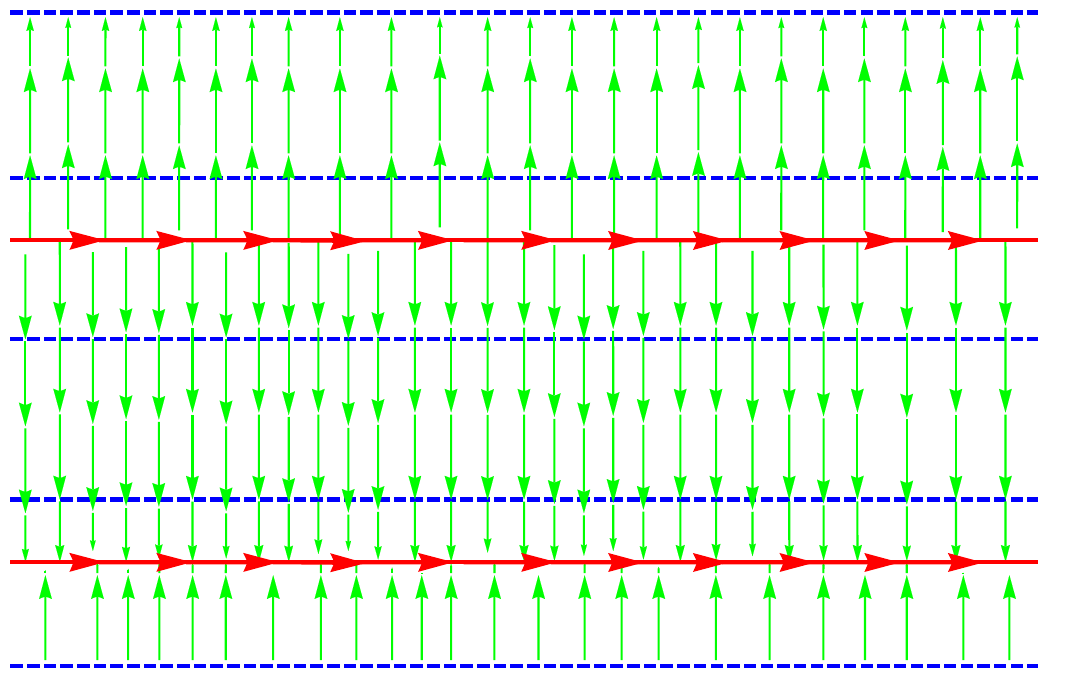}}%
    \put(0.97107729,0.0010595){\color[rgb]{0,0,0}\makebox(0,0)[lt]{\lineheight{1.25}\smash{\begin{tabular}[t]{l}$0$\end{tabular}}}}%
    \put(0.97107729,0.30436762){\color[rgb]{0,0,0}\makebox(0,0)[lt]{\lineheight{1.25}\smash{\begin{tabular}[t]{l}$\pi$\end{tabular}}}}%
    \put(0.97107729,0.60660139){\color[rgb]{0,0,0}\makebox(0,0)[lt]{\lineheight{1.25}\smash{\begin{tabular}[t]{l}$2\pi$\end{tabular}}}}%
    \put(0.97107729,0.39555488){\color[rgb]{0,0,0}\makebox(0,0)[lt]{\lineheight{1.25}\smash{\begin{tabular}[t]{l}$L_{i}^{\pi}$\end{tabular}}}}%
    \put(0.97107729,0.09724682){\color[rgb]{0,0,0}\makebox(0,0)[lt]{\lineheight{1.25}\smash{\begin{tabular}[t]{l}$L_{i}$\end{tabular}}}}%
  \end{picture}%
\endgroup%

  \caption{Constant double discontinuity dynamics for $d_{i1} = 1 > 0$, $d_{i2}
      = 0.7 > 0$ and $d_{i3} = 1 > 0$. At this example we have $\theta_i =
      \arctan{\frac{1}{0.7}} \approx 0.96$. Therefore, for example, $S_1$ has
    $\theta = \theta_i$ as an attracting visible part of the slow manifold;
    whereas $S_2$ has none.}
  \label{sec:constant:fig:d2_nonzero}
\end{figure}
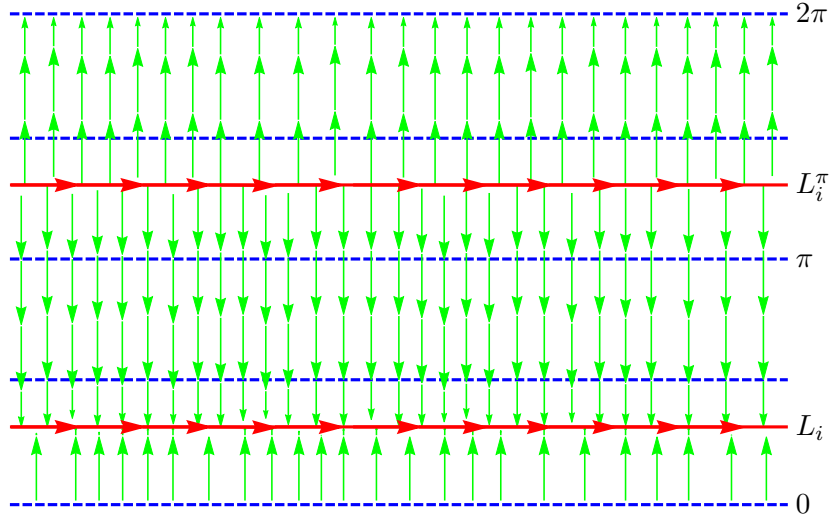

In order to explicitly define the slow manifold $\mathcal{M}_i$, observe that
whenever $\cos{\theta} \neq 0$ the second equation of
\eqref{sec:constant:eq:reduced} gives us

\begin{align*}
  0 & = d_{i3} \cos{\theta} - d_{i2} \sin{\theta}
  \Leftrightarrow \tan{\theta} = \frac{d_{i3}}{d_{i2}} \Leftrightarrow                               \\
    & \Leftrightarrow \theta = \arctan{\left(\frac{d_{i3}}{d_{i2}}\right)} + n\pi = \theta_i + n\pi,
\end{align*}

\noindent where $n \in \mathbb{Z}$. Therefore, without loss of generality, the
slow manifold can be written as $\mathcal{M}_i = L_{i} \cup L_{i}^{\pi}$, where

\begin{align*}
  L_{i}       & = \left\{(x,\theta) \in \mathbb{R} \times \left[0,2\pi\right] ;~
  \theta = \theta_i \right\} \text{ and}                                         \\
  L_{i}^{\pi} & = \left\{(x,\theta) \in \mathbb{R} \times \left[0,2\pi\right] ;~
  \theta = \theta_i + \pi \right\},
\end{align*}

\noindent which consists of two straight lines inside the cylinder $C =
  \mathbb{R} \times \left[0,2\pi\right]$, as the red part of
\cref{sec:constant:fig:d2_nonzero}. In fact, \emph{a priori}, $\mathcal{M}_i$
is a subset of the particular stripe $S_i$. However, since the subjacent
vector fields, \cref{sec:constant:eq:system}, are defined for every point of
$\mathbb{R}^3$, then, without any mathematical restriction or weakness, we can
consider $\mathcal{M}_i$ as a subset of the whole cylinder $C$, not restricted
to the particular stripe $S_i$, in order to study its properties. Once this
global analysis is done, we can then focus on the particular stripe of
interest.

In particular, since $\theta_i \in \left(-\frac{\pi}{2}, \frac{\pi}{2}\right)$
and $\theta_i + \pi \in \left(\frac{\pi}{2}, \frac{3\pi}{2}\right)$, then either
$L_{i} \subset S_1$ and $L_{i}^{\pi} \subset S_3$ or $L_{i} \subset S_4$ and
$L_{i}^{\pi} \subset S_2$. In other words, this straight lines are always at
intercalated stripes. Therefore, a given stripe $S_i$ might or might not contain
one of this straight lines, depending exclusively on the value of
$\theta_i$.\footnote{In particular, when $d_{i3} = 0$ we have $\theta_i = 0$
  and, therefore, the straight lines $L_{i}$ and $L_{i}^{\pi}$ are given by
  $\theta = 0$ and $\theta = \pi$, respectively, which are part of the stripes'
  boundary.} This completes the qualitative analysis of the shape of the slow
manifold.

Over both the straight lines $\mathcal{M}_i = L_{i} \cup L_{i}^{\pi}$, we have
the one-dimensional dynamics given by the first equation of
\eqref{sec:constant:eq:reduced}, i.e., $\dot{x} = d_{i1}$. Analyzing this
equation we observe that, considering the usual growth direction of the
$x$-axis, the dynamics over $\mathcal{M}_i$ is increasing if $d_{i1} > 0$ and
decreasing if $d_{i1} < 0$. This completes the qualitative analysis of the
reduced dynamics.

Regarding the layer dynamics, we have the layer system
\eqref{sec:constant:eq:layer} which says that for each fixed value of $x \in
  \mathbb{R}$, we have a one-dimensional dynamics given by the second equation of
\eqref{sec:constant:eq:layer}. In particular, assuming that $\cos{\theta} > 0$
and $d_{i2} > 0$, then

\begin{align*}
  \theta' > 0 & \Leftrightarrow
  d_{i3} \cos{\theta} - d_{i2} \sin{\theta} > 0 \Leftrightarrow
  \tan{\theta} < \frac{d_{i3}}{d_{i2}} \Leftrightarrow                                           \\
              & \Leftrightarrow \theta < \arctan{\left(\frac{d_{i3}}{d_{i2}}\right)} = \theta_i,
\end{align*}

\noindent since the arctangent function is strictly increasing. Likewise and
under the same conditions, we have that

\begin{equation*}
  \theta' < 0 \Leftrightarrow
  \theta > \arctan{\left(\frac{d_{i3}}{d_{i2}}\right)} = \theta_i
\end{equation*}

\noindent and, therefore, we conclude that for $d_{i2} > 0$, the straight line
$L_{i}$ is attractor of surrounding layer dynamics and, therefore, $L_{i}^{\pi}$
is a repellor, as the green part of \cref{sec:constant:fig:d2_nonzero}. An
analogous study for $d_{i2} < 0$ allows us to reach the results summarized in
\cref{sec:constant:tab:fast}.

\begin{table}[ht]
  \centering
  \caption{Layer dynamics around the straight lines $L_{i}$ and $L_{i}^{\pi}$
    that compose the slow manifold $\mathcal{M}_i = L_{i} \cup L_{i}^{\pi}$.}
  \begin{tabular}{lll}
    \hline
                  & $d_{i2} < 0$ & $d_{i2} > 0$ \\
    \hline
    $L_{i}$       & repellor     & attractor    \\
    $L_{i}^{\pi}$ & attractor    & repellor     \\
    \hline
  \end{tabular}
  \label{sec:constant:tab:fast}
\end{table}

Finally, at $\cos{\theta} = 0$ with $d_{i2} \neq 0$ the reduced system
\eqref{sec:constant:eq:reduced} tells us that $\mathcal{M}_i = \emptyset$ and,
therefore, there is only the fast dynamics \eqref{sec:constant:eq:layer} which
reduces to

\begin{equation*}
  \left\{\begin{matrix*}[l]
    x' & = & 0 \\
    \theta' & = & - d_{i2}
  \end{matrix*}\right.
  \quad \text{ and } \quad
  \left\{\begin{matrix*}[l]
    x' & = & 0 \\
    \theta' & = & d_{i2}
  \end{matrix*}\right.
\end{equation*}

\noindent for $\theta = \frac{\pi}{2}$ and $\theta = \frac{3\pi}{2}$,
respectively, whose dynamics is consistent with \cref{sec:constant:tab:fast}.
This completes the qualitative analysis of the layer dynamics and, therefore, the
qualitative analysis of this case. See
\cref{sec:constant:exmp:no_singularities}.

\subsection{Case \texorpdfstring{$d_{i2} = 0$}{d2=0}}%
\label{sec:constant:sub:d2_zero}

Now, the reduced system \eqref{sec:constant:eq:reduced} can be written as

\begin{equation}
  \label{sec:constant:sub:d2_zero:eq:reduced}
  \left\{\begin{matrix*}[l]
    \dot{x} & = & d_{i1} \\
    0 & = & d_{i3} \cos{\theta}
  \end{matrix*}\right.,
\end{equation}

\noindent whose slow manifold $\mathcal{M}_i$ is implicitly given by the
equation $0 = d_{i3} \cos{\theta}$ which actually means $0 = \cos{\theta}$,
since we are under \cref{sec:framework:sub:blowup:eq:weak_fund_hypo} and,
therefore, $d_{i3} \neq 0$. In other words, $\mathcal{M}_i = L_{i} \cup
L_{i}^{\pi}$ with $L_{i}$ and $L_{i}^{\pi}$ being the straight lines given by
$\theta = \frac{\pi}{2}$ and $\theta = \frac{3\pi}{2}$,
respectively.\footnote{Here, again, the straight lines $L_{i}$ and $L_{i}^{\pi}$
are part of the boundary of the stripes.} The dynamics over and around
$\mathcal{M}_i$ behaves exactly as in the case $d_{i2} \neq 0$, but exchanging
$d_{i2}$ with $d_{i3}$ at \cref{sec:constant:tab:fast}.

\subsection{Theorem and Examples}%
\label{sec:constant:sub:theorem}

Summarizing, we conclude that the dynamics over $\Sigma_x$ for constant fields
behaves as described in the theorem below, whose proof consists in the analysis
done above in \cref{sec:constant:sub:d2_nonzero,sec:constant:sub:d2_zero}.

\begin{rtheorem}[Constant Dynamics]
  \label{sec:constant:thm:dynamics}
  Given $\mathbf{F} \in \mathcal{C}_3$ with constant components $\mathbf{F}_i =
    (d_{i1},d_{i2},d_{i3})$ such that $d_{i2} \neq 0$ or $d_{i3} \neq 0$, let
  $\mathbf{\tilde{F}} \in \tilde{\mathcal{C}_3}$ be the vector field induced
  by the blow-up $\phi_1(x,\theta,r) = (x, r\cos{\theta}, r\sin{\theta})$.
  Then, this blow-up associates the dynamics over $\Sigma_x$ with the
  following fundamental dynamics over the cylinder $C = \mathbb{R} \times S^1
    = S_1 \cup \ldots \cup S_4$: over each stripe $S_i$ acts a slow-fast
  dynamics whose slow manifold is given by $\mathcal{M}_i = L_{i} \cup
    L_{i}^{\pi}$, where $L_{i}^{\pi}$ is a $\pi$-translation of $L_{i}$ in
  $\theta$ and

  \begin{enumerate}
    \item \label{sec:constant:thm:dynamics:d2_nonzero} case $d_{i2} \neq 0$,
          then

          \begin{equation*}
            L_{i} = \left\{(x,\theta) \in \mathbb{R} \times \left[0,2\pi\right] ;~
            \theta = \arctan{\left(\frac{d_{i3}}{d_{i2}}\right)} \right\};
          \end{equation*}

    \item \label{sec:constant:thm:dynamics:d2_zero_d3_nonzero} case $d_{i2} = 0$
          and $d_{i3} \neq 0$, then

          \begin{equation*}
            L_{i} = \left\{(x,\theta) \in \mathbb{R} \times \left[0,2\pi\right] ;~
            \theta = \frac{\pi}{2} \right\};
          \end{equation*}
  \end{enumerate}

  \noindent which, in both cases, consists of two straight lines inside the
  cylinder $C$, possibly invisible relative to $S_i$. Over this straight
  lines acts the reduced dynamics $\dot{x} = d_{i1}$ and, around then,
  acts the layer dynamics described in \cref{sec:constant:tab:fast}, but
  exchanging $d_{i2}$ with $d_{i3}$ if $d_{i2} = 0$.
\end{rtheorem}

\begin{figure}[ht]
  \centering
  \begin{subfigure}[b]{0.45\linewidth}
    \def\svgwidth{\linewidth}
\begingroup%
  \makeatletter%
  \providecommand\color[2][]{%
    \errmessage{(Inkscape) Color is used for the text in Inkscape, but the package 'color.sty' is not loaded}%
    \renewcommand\color[2][]{}%
  }%
  \providecommand\transparent[1]{%
    \errmessage{(Inkscape) Transparency is used (non-zero) for the text in Inkscape, but the package 'transparent.sty' is not loaded}%
    \renewcommand\transparent[1]{}%
  }%
  \providecommand\rotatebox[2]{#2}%
  \newcommand*\fsize{\dimexpr\f@size pt\relax}%
  \newcommand*\lineheight[1]{\fontsize{\fsize}{#1\fsize}\selectfont}%
  \ifx\svgwidth\undefined%
    \setlength{\unitlength}{251.47663988bp}%
    \ifx\svgscale\undefined%
      \relax%
    \else%
      \setlength{\unitlength}{\unitlength * \real{\svgscale}}%
    \fi%
  \else%
    \setlength{\unitlength}{\svgwidth}%
  \fi%
  \global\let\svgwidth\undefined%
  \global\let\svgscale\undefined%
  \makeatother%
  \begin{picture}(1,0.92736389)%
    \lineheight{1}%
    \setlength\tabcolsep{0pt}%
    \put(0,0){\includegraphics[width=\unitlength,page=1]{slice.pdf}}%
    \put(0.48325888,0.39247149){\color[rgb]{0,0,0}\makebox(0,0)[lt]{\lineheight{1.25}\smash{\begin{tabular}[t]{l}$\Sigma_x$\end{tabular}}}}%
    \put(0.77900636,0.78367183){\color[rgb]{0,0,0}\makebox(0,0)[lt]{\lineheight{0}\smash{\begin{tabular}[t]{l}$\mathbf{F}_1$\end{tabular}}}}%
    \put(0.12457511,0.13116122){\color[rgb]{0,0,0}\makebox(0,0)[lt]{\lineheight{0}\smash{\begin{tabular}[t]{l}$\mathbf{F}_3$\end{tabular}}}}%
    \put(0.78036067,0.1228519){\color[rgb]{0,0,0}\makebox(0,0)[lt]{\lineheight{0}\smash{\begin{tabular}[t]{l}$\mathbf{F}_4$\end{tabular}}}}%
    \put(0.12668394,0.78490013){\color[rgb]{0,0,0}\makebox(0,0)[lt]{\lineheight{0}\smash{\begin{tabular}[t]{l}$\mathbf{F}_2$\end{tabular}}}}%
  \end{picture}%
\endgroup%

    \caption{Before the blow-up.}%
    \label{sec:constant:exmp:no_singularities:subfig:slice}
  \end{subfigure}
  \qquad
  \begin{subfigure}[b]{0.45\linewidth}
    \def\svgwidth{\linewidth}
    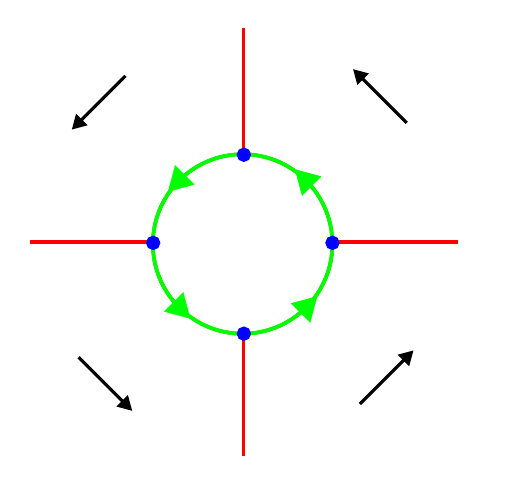
    \caption{After the blow-up.}%
    \label{sec:constant:exmp:no_singularities:subfig:slice_blowup}
  \end{subfigure}
  \caption{Slices of the system studied at
    \cref{sec:constant:exmp:no_singularities}.}
  \label{sec:constant:exmp:no_singularities:fig:slices}
\end{figure}

\begin{example}
  \label{sec:constant:exmp:no_singularities}
  Let $\mathbf{F} \in \mathcal{C}_3$ be given by the constant vector fields

  \begin{align*}
    \mathbf{F}_2(x,y,z) & = (1,-1,-1), & \mathbf{F}_1(x,y,z) & = (1,-1,1), \\
    \mathbf{F}_3(x,y,z) & = (1,1,-1),  & \mathbf{F}_4(x,y,z) & = (1,1,1),
  \end{align*}

  \noindent that behaves as represented at
  \cref{sec:constant:exmp:no_singularities:subfig:slice}. Using
  \cref{sec:constant:thm:dynamics} we can verify that, over the cylinder $C$
  given by the blow-up of $\Sigma_x$, this system behaves as expected, i.e., as
  represented at \cref{sec:constant:exmp:no_singularities:subfig:slice_blowup}.

  For instance, over the stripe $S_1 = \mathbb{R} \times
    \left[0,\sfrac{\pi}{2}\right]$ we have

  \begin{equation*}
    (d_{11},d_{12},d_{13}) = \mathbf{F}_1(x,y,z) = (1,-1,1)
  \end{equation*}

  \noindent such that, according to \cref{sec:constant:thm:dynamics}, induces
  over $S_1$ a slow-fast system with $L_{1} \subset \mathcal{M}_1$ given by

  \begin{equation*}
    \theta = \theta_1 = \arctan{\left(\frac{d_{13}}{d_{12}}\right)} =
    \arctan{\left(\frac{1}{-1}\right)} =
    -\frac{\pi}{4},
  \end{equation*}

  \noindent and, therefore, the slow manifold $\mathcal{M}_1$ consists of the
  straight lines $L_{1} \subset S_4$ and $L_{1}^{\pi} \subset S_2$ given by
  $\theta = \theta_1 = - \frac{\pi}{4}$ and $\theta = \theta_1 + \pi =
    \frac{3\pi}{4}$, respectively. In particular, none of these lines are visible
  at $S_1$. Over these lines acts the reduced dynamics $\dot{x} = d_{11} = 1$.
  Finally, since $d_{12} = -1 < 0$, then $L_{1}$ is repellor and $L_{1}^{\pi}$
  is attractor of surrounding layer dynamics, according to
  \cref{sec:constant:tab:fast}.

  Therefore, we conclude that the dynamics generated by $\mathbf{F}_1$ over the
  whole cylinder $C$ behaves as represented in
  \cref{sec:constant:exmp:no_singularities:fig:stripe}. In particular, the
  dynamics over the stripe $S_1$ behaves as represented in
  \cref{sec:constant:exmp:no_singularities:subfig:slice_blowup}. The dynamics
  over the other stripes can be similarly verified to be as represented. \qed

  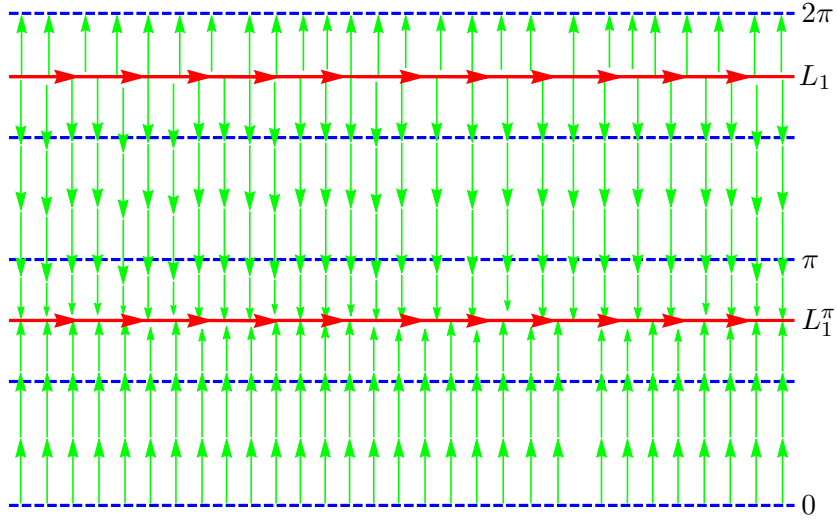
\begin{figure}[ht]
    \centering
    \def\svgwidth{0.85\linewidth}
\begingroup%
  \makeatletter%
  \providecommand\color[2][]{%
    \errmessage{(Inkscape) Color is used for the text in Inkscape, but the package 'color.sty' is not loaded}%
    \renewcommand\color[2][]{}%
  }%
  \providecommand\transparent[1]{%
    \errmessage{(Inkscape) Transparency is used (non-zero) for the text in Inkscape, but the package 'transparent.sty' is not loaded}%
    \renewcommand\transparent[1]{}%
  }%
  \providecommand\rotatebox[2]{#2}%
  \newcommand*\fsize{\dimexpr\f@size pt\relax}%
  \newcommand*\lineheight[1]{\fontsize{\fsize}{#1\fsize}\selectfont}%
  \ifx\svgwidth\undefined%
    \setlength{\unitlength}{513.10205078bp}%
    \ifx\svgscale\undefined%
      \relax%
    \else%
      \setlength{\unitlength}{\unitlength * \real{\svgscale}}%
    \fi%
  \else%
    \setlength{\unitlength}{\svgwidth}%
  \fi%
  \global\let\svgwidth\undefined%
  \global\let\svgscale\undefined%
  \makeatother%
  \begin{picture}(1,0.62999164)%
    \lineheight{1}%
    \setlength\tabcolsep{0pt}%
    \put(0,0){\includegraphics[width=\unitlength,page=1]{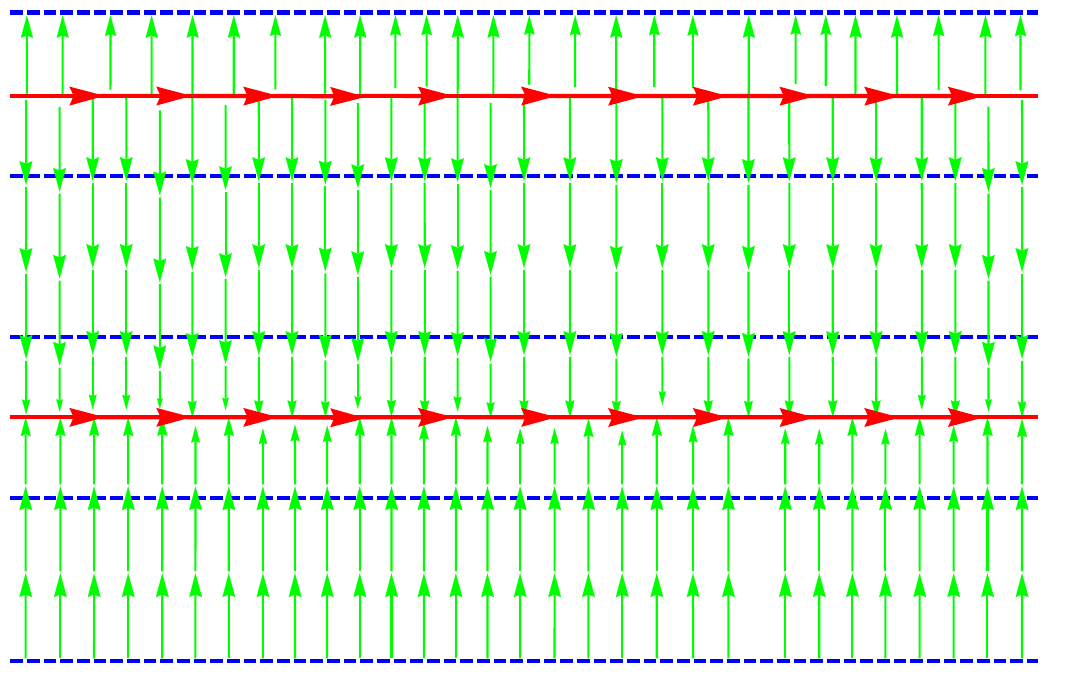}}%
    \put(0.97937981,0.00106963){\color[rgb]{0,0,0}\makebox(0,0)[lt]{\lineheight{1.25}\smash{\begin{tabular}[t]{l}$0$\end{tabular}}}}%
    \put(0.97937981,0.30467353){\color[rgb]{0,0,0}\makebox(0,0)[lt]{\lineheight{1.25}\smash{\begin{tabular}[t]{l}$\pi$\end{tabular}}}}%
    \put(0.97937981,0.60814992){\color[rgb]{0,0,0}\makebox(0,0)[lt]{\lineheight{1.25}\smash{\begin{tabular}[t]{l}$2\pi$\end{tabular}}}}%
    \put(0.97605132,0.53003964){\color[rgb]{0,0,0}\makebox(0,0)[lt]{\lineheight{1.25}\smash{\begin{tabular}[t]{l}$L_{1}$\end{tabular}}}}%
    \put(0.97918943,0.22944949){\color[rgb]{0,0,0}\makebox(0,0)[lt]{\lineheight{1.25}\smash{\begin{tabular}[t]{l}$L_{1}^{\pi}$\end{tabular}}}}%
  \end{picture}%
\endgroup%

    \caption{Dynamics over $C$ generated by the field $\mathbf{F}_1$ studied at
      \cref{sec:constant:exmp:no_singularities}. The dynamics over $S_1$ behaves
      as represented in
      \cref{sec:constant:exmp:no_singularities:subfig:slice_blowup}.}
    \label{sec:constant:exmp:no_singularities:fig:stripe}
  \end{figure}
\end{example}

\section{Affine Dynamics}%
\label{sec:affine}

Let $\mathcal{A}_3 \subset \mathcal{D}_3^k$ be the set of all piecewise smooth
vector fields $\mathbf{F}$ with a double discontinuity given by the affine
vector fields

\begin{equation}
  \label{sec:affine:eq:system}
  \begin{aligned}
    \mathbf{F}_i(x,y,z) =
    ( & a_{i1}x + b_{i1}y + c_{i1}z + d_{i1},  \\
      & a_{i2}x + b_{i2}y + c_{i2}z + d_{i2},  \\
      & a_{i3}x + b_{i3}y + c_{i3}z + d_{i3}),
  \end{aligned}
\end{equation}

\noindent where $a_{ij}, b_{ij}, c_{ij}, d_{ij} \in \mathbb{R}$ for all $i$ and
$j$. According to the Fundamental \cref{sec:framework:blowup:thm:dynamics}, the
dynamics over $\Sigma_x$ of such a field is blow-up associated to the following
fundamental dynamics over the cylinder $C = \mathbb{R} \times S^1 = S_1 \cup
  \ldots \cup S_4$: over each stripe $S_i$ acts a slow-fast dynamics whose reduced
dynamics is given by

\begin{equation}
  \label{sec:affine:eq:reduced}
  \left\{\begin{matrix*}[l]
    \dot{x} = a_{i1}x + d_{i1}\\
    0 = (a_{i3}x + d_{i3}) \cos{\theta} - (a_{i2}x + d_{i2}) \sin{\theta}\\
  \end{matrix*}\right.,
\end{equation}

\noindent with radial slow dynamics $\dot{r} = (a_{i2}x + d_{i2}) \cos{\theta} +
  (a_{i3}x + d_{i3}) \sin{\theta}$; and layer dynamics given by

\begin{equation}
  \label{sec:affine:eq:layer}
  \left\{\begin{matrix*}[l]
    x' = 0\\
    \theta' = (a_{i3}x + d_{i3}) \cos{\theta} - (a_{i2}x + d_{i2}) \sin{\theta}
  \end{matrix*}\right.,
\end{equation}

\noindent with radial fast dynamics $r' = 0$.

Besides that, for \cref{sec:affine:eq:system}, we have $p_i = a_{i2}x + d_{i2}$
and $q_i = a_{i3}x + d_{i3}$ so that
\cref{sec:framework:sub:blowup:eq:weak_fund_hypo} is satisfied as long as

\begin{equation}
  \label{sec:affine:eq:wfh}
  a_{i2}x + d_{i2} \neq 0 \quad \text{or} \quad a_{i3}x + d_{i3} \neq 0,
\end{equation}

\noindent whereas, since $(p_i)_x = a_{i2}$ and $(q_i)_x = a_{i3}$, then
\cref{sec:framework:sub:blowup:eq:strong_fund_hypo} is satisfied as long as

\begin{equation}
  \label{sec:affine:eq:sfh}
  \begin{aligned}
    0 & \neq p_i(q_i)_x - q_i(p_i)_x =                          \\
      & = (a_{i2}x + d_{i2})a_{i3} - (a_{i3}x + d_{i3})a_{i2} = \\
      & = a_{i3}d_{i2} - a_{i2}d_{i3} \eqqcolon \gamma_i,
  \end{aligned}
\end{equation}

\noindent which not only assures the fundamental hypothesis but also avoids the
already studied constant case, as we will see below.

As in the constant case, our goal at this section is to fully describe the
fundamental dynamics of \cref{sec:affine:eq:system} over the cylinder $C$ under
the hypothesis \cref{sec:affine:eq:sfh}. In order to do so, we are going to
systematically analyze the slow-fast systems
\eqref{sec:affine:eq:reduced}--\eqref{sec:affine:eq:layer} for the cases
suggested by \cref{sec:affine:eq:wfh} and outlined at
\cref{sec:affine:tab:cases}.

\begin{table}[ht]
  \centering
  \caption{Division \eqref{sec:affine:eq:system} dynamics in study cases.}
  \begin{tabular}{ccc}
    \hline
                    & $a_{i2}x + d_{i2} \neq 0$ & $a_{i2}x + d_{i2} = 0$ \\
    \hline
    $a_{i2} \neq 0$ & A                         & B                      \\
    $a_{i2} = 0$    & C                         & D                      \\
    \hline
  \end{tabular}%
  \label{sec:affine:tab:cases}
\end{table}

Observe that case (B) actually complements case (A). Moreover, observe that at
case (D) we have $a_{i2} = 0$ and $d_{i2} = 0$ which implies the absurd
$\gamma_i = 0$. Therefore, cases (A) and (B) complement each other and it will
be studied at \cref{sec:affine:sub:a2_nonzero}; case (C) will be studied at
\cref{sec:affine:sub:a2_zero}. The resulting \cref{sec:affine:thm:dynamics} is
stated and exemplified at \cref{sec:affine:sub:theorem}.

\subsection{Case \texorpdfstring{$a_{i2} \neq 0$}{a2=/=0}}
\label{sec:affine:sub:a2_nonzero}

Lets start with case (A), i.e., assume that $a_{i2} \neq 0$ and $a_{i2}x +
  d_{i2} \neq 0$.  In order to explicitly define $\mathcal{M}_i$, observe that
whenever $\cos{\theta} \neq 0$ the second equation of
\eqref{sec:affine:eq:reduced} gives us

\begin{equation*}
  \begin{aligned}[l]
    0 & = (a_{i3}x + d_{i3}) \cos{\theta} - (a_{i2}x + d_{i2}) \sin{\theta}
    \Leftrightarrow                                                         \\
      & \Leftrightarrow
    \tan{\theta} = \frac{a_{i3}x + d_{i3}}{a_{i2}x + d_{i2}} \eqqcolon h(x)
    \Leftrightarrow                                                         \\
      & \Leftrightarrow
    \theta = \arctan{\left(\frac{a_{i3}x + d_{i3}}{a_{i2}x + d_{i2}}\right)}
    + n\pi = \theta_i\left(x \right) + n\pi,
  \end{aligned}%
\end{equation*}

\noindent where $n \in \mathbb{Z}$. As in the constant case, since the subjacent
vector fields, \cref{sec:affine:eq:system}, are defined for every point of
$\mathbb{R}^3$, then we can consider $\mathcal{M}_i$ as a subset of the whole
cylinder $C$, not restricted to the particular stripe $S_i$. Therefore, without
loss of generality, the slow manifold can be written as $\mathcal{M}_i = H_{i}
  \cup H_{i}^{\pi}$, where

\begin{align*}
  H_{i}       & = \left\{(x,\theta) \in \mathbb{R} \times \left[0,2\pi\right] ;~
  \theta = \theta_i(x) \right\} \text{ and}                                      \\
  H_{i}^{\pi} & = \left\{(x,\theta) \in \mathbb{R} \times \left[0,2\pi\right] ;~
  \theta = \theta_i(x) + \pi \right\},
\end{align*}

\noindent which consists of two arctangent-normalized hyperboles inside the
cylinder $C = \mathbb{R} \times S^1$. In fact, since $a_{i2} \neq 0$, then
$h(x)$ is a hyperbole such that

\begin{align*}
  \frac{d}{dx}h(x) & =
  \frac{d}{dx}\left[\frac{a_{i3}x + d_{i3}}{a_{i2}x + d_{i2}}\right] =
  \frac{a_{i3}d_{i2} - d_{i3}a_{i2}}{(a_{i2}x + d_{i2})^2} = \\
                   & = \frac{\gamma_i}{(a_{i2}x + d_{i2})^2}
\end{align*}

\noindent or, in other words, it is an increasing hyperbole if $\gamma_i > 0$
and decreasing if $\gamma_i < 0$\footnote{If $\gamma_i = 0$, then $h(x)$ is a
  constant function and, therefore, $H_{i}$ and $H_{i}^{\pi}$ are straight lines.
  In other words, the constant case is recovered.}.

Besides that, observe that $h(x)$ has a vertical asymptote at

\begin{equation*}
  a_{i2}x + d_{i2} = 0
  \Leftrightarrow
  x = -\frac{d_{i2}}{a_{i2}} \eqqcolon \alpha_i
\end{equation*}

\noindent which satisfies

\begin{equation*}
  \lim_{x \to \alpha_i^{\pm}} h(x) = \mp\infty
  \quad \text{and} \quad
  \lim_{x \to \alpha_i^{\pm}} h(x) = \pm\infty
\end{equation*}

\noindent if $\gamma_i > 0$ and $\gamma_i < 0$, respectively; and $h(x)$ has a
horizontal asymptote at

\begin{equation*}
  \lim_{x \to \pm\infty} h(x) =
  \lim_{x \to \pm\infty} \left(\frac{a_{i3}x + d_{i3}}
  {a_{i2}x + d_{i2}} \right) = \frac{a_{i3}}{a_{i2}}.
\end{equation*}

\noindent Translating the information above about the hyperbole $h(x)$ to the
arctangent-normalized hyperbole $H_{i}$, we get that it

\begin{itemize}
  \item is an increasing curve if $\gamma_i > 0$ and decreasing if $\gamma_i <
          0$;

  \item has a vertical asymptote at $x = \alpha_i$ which satisfies

        \begin{equation*}
          \lim_{x \to \alpha_i^{\pm}} \theta_{i}(x) = \mp\frac{\pi}{2}
          \quad \text{and} \quad
          \lim_{x \to \alpha_i^{\pm}} \theta_{i}(x) = \pm\frac{\pi}{2}
        \end{equation*}

        \noindent if $\gamma_i > 0$ and $\gamma_i < 0$, respectively;

  \item has a horizontal asymptote at $\theta =
          \arctan{\left(\frac{a_{i3}}{a_{i2}} \right)} \eqqcolon \beta_i$.
\end{itemize}

\begin{figure}[ht]
  \centering
  \begin{subfigure}[b]{0.85\linewidth}
    \def\svgwidth{\linewidth}
\begingroup%
  \makeatletter%
  \providecommand\color[2][]{%
    \errmessage{(Inkscape) Color is used for the text in Inkscape, but the package 'color.sty' is not loaded}%
    \renewcommand\color[2][]{}%
  }%
  \providecommand\transparent[1]{%
    \errmessage{(Inkscape) Transparency is used (non-zero) for the text in Inkscape, but the package 'transparent.sty' is not loaded}%
    \renewcommand\transparent[1]{}%
  }%
  \providecommand\rotatebox[2]{#2}%
  \newcommand*\fsize{\dimexpr\f@size pt\relax}%
  \newcommand*\lineheight[1]{\fontsize{\fsize}{#1\fsize}\selectfont}%
  \ifx\svgwidth\undefined%
    \setlength{\unitlength}{502.87051392bp}%
    \ifx\svgscale\undefined%
      \relax%
    \else%
      \setlength{\unitlength}{\unitlength * \real{\svgscale}}%
    \fi%
  \else%
    \setlength{\unitlength}{\svgwidth}%
  \fi%
  \global\let\svgwidth\undefined%
  \global\let\svgscale\undefined%
  \makeatother%
  \begin{picture}(1,0.34579226)%
    \lineheight{1}%
    \setlength\tabcolsep{0pt}%
    \put(0,0){\includegraphics[width=\unitlength,page=1]{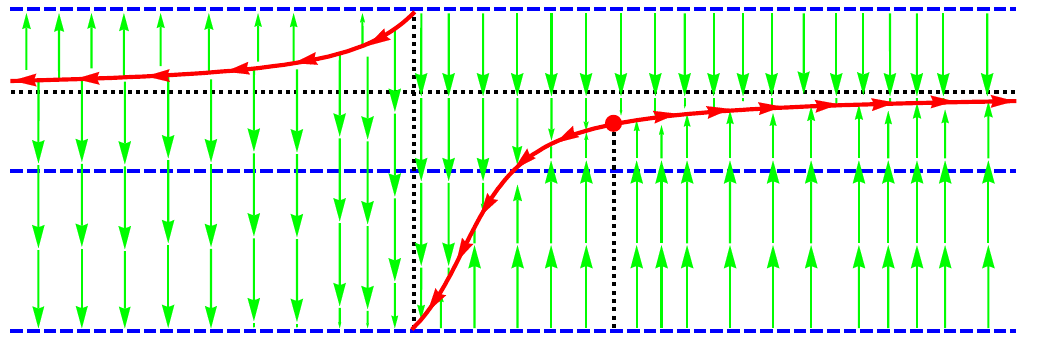}}%
    \put(0.97994631,0.01947414){\color[rgb]{0,0,0}\makebox(0,0)[lt]{\lineheight{1.25}\smash{\begin{tabular}[t]{l}$-\frac{\pi}{2}$\end{tabular}}}}%
    \put(0.97994631,0.17199689){\color[rgb]{0,0,0}\makebox(0,0)[lt]{\lineheight{1.25}\smash{\begin{tabular}[t]{l}$0$\end{tabular}}}}%
    \put(0.97994631,0.32671029){\color[rgb]{0,0,0}\makebox(0,0)[lt]{\lineheight{1.25}\smash{\begin{tabular}[t]{l}$\frac{\pi}{2}$\end{tabular}}}}%
    \put(0.97994631,0.24788548){\color[rgb]{0,0,0}\makebox(0,0)[lt]{\lineheight{1.25}\smash{\begin{tabular}[t]{l}$\beta_i$\end{tabular}}}}%
    \put(0.38595585,-0.00441499){\color[rgb]{0,0,0}\makebox(0,0)[lt]{\lineheight{1.25}\smash{\begin{tabular}[t]{l}$x = \alpha_i$\end{tabular}}}}%
    \put(0.57689652,-0.00441499){\color[rgb]{0,0,0}\makebox(0,0)[lt]{\lineheight{1.25}\smash{\begin{tabular}[t]{l}$x = \delta_i$\end{tabular}}}}%
    \put(0.59777664,0.19543765){\color[rgb]{0,0,0}\makebox(0,0)[lt]{\lineheight{1.25}\smash{\begin{tabular}[t]{l}$\mathbf{P}$\end{tabular}}}}%
  \end{picture}%
\endgroup%

    \caption{Hyperbole $H_{i}$.}
    \label{sec:affine:sub:a2_nonzero:subfig:hyperbole}
  \end{subfigure}
  \quad
  \begin{subfigure}[b]{0.85\linewidth}
    \def\svgwidth{\linewidth}
\begingroup%
  \makeatletter%
  \providecommand\color[2][]{%
    \errmessage{(Inkscape) Color is used for the text in Inkscape, but the package 'color.sty' is not loaded}%
    \renewcommand\color[2][]{}%
  }%
  \providecommand\transparent[1]{%
    \errmessage{(Inkscape) Transparency is used (non-zero) for the text in Inkscape, but the package 'transparent.sty' is not loaded}%
    \renewcommand\transparent[1]{}%
  }%
  \providecommand\rotatebox[2]{#2}%
  \newcommand*\fsize{\dimexpr\f@size pt\relax}%
  \newcommand*\lineheight[1]{\fontsize{\fsize}{#1\fsize}\selectfont}%
  \ifx\svgwidth\undefined%
    \setlength{\unitlength}{503.10632324bp}%
    \ifx\svgscale\undefined%
      \relax%
    \else%
      \setlength{\unitlength}{\unitlength * \real{\svgscale}}%
    \fi%
  \else%
    \setlength{\unitlength}{\svgwidth}%
  \fi%
  \global\let\svgwidth\undefined%
  \global\let\svgscale\undefined%
  \makeatother%
  \begin{picture}(1,0.64397686)%
    \lineheight{1}%
    \setlength\tabcolsep{0pt}%
    \put(0,0){\includegraphics[width=\unitlength,page=1]{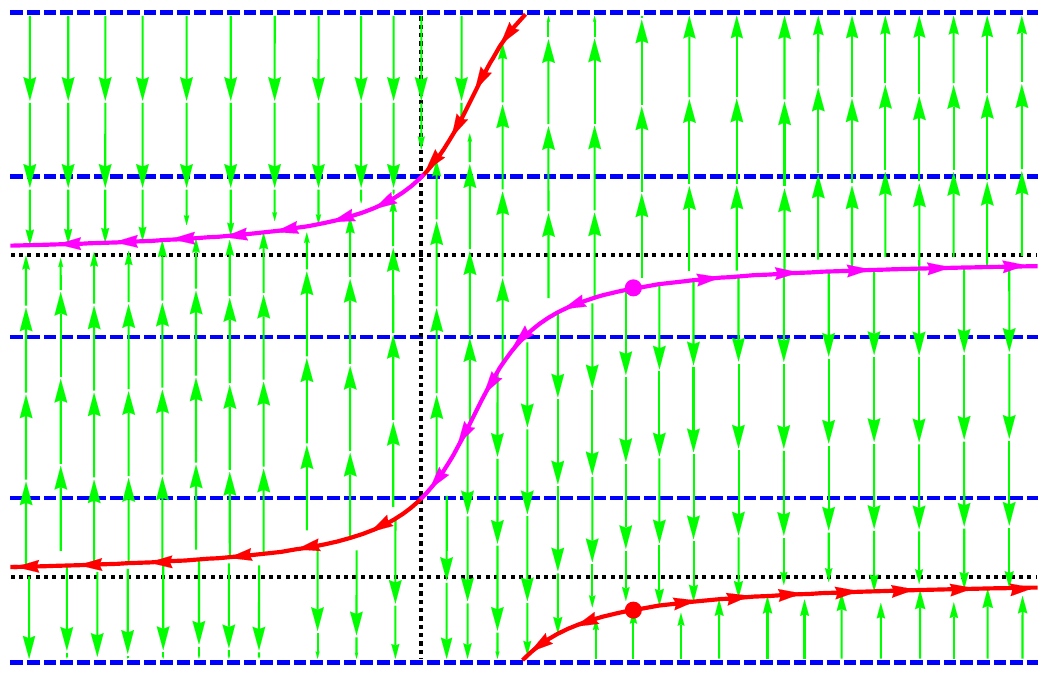}}%
    \put(0.6153852,0.0311463){\color[rgb]{0,0,0}\makebox(0,0)[lt]{\lineheight{1.25}\smash{\begin{tabular}[t]{l}$\mathbf{P}$\end{tabular}}}}%
    \put(0.6153852,0.33525697){\color[rgb]{0,0,0}\makebox(0,0)[lt]{\lineheight{1.25}\smash{\begin{tabular}[t]{l}$\mathbf{P}^{\pi}$\end{tabular}}}}%
    \put(1.00149341,0.00109089){\color[rgb]{0,0,0}\makebox(0,0)[lt]{\lineheight{1.25}\smash{\begin{tabular}[t]{l}$0$\end{tabular}}}}%
    \put(1.00149341,0.08357813){\color[rgb]{0,0,0}\makebox(0,0)[lt]{\lineheight{1.25}\smash{\begin{tabular}[t]{l}$\beta_i$\end{tabular}}}}%
    \put(1.00149341,0.31182909){\color[rgb]{0,0,0}\makebox(0,0)[lt]{\lineheight{1.25}\smash{\begin{tabular}[t]{l}$\pi$\end{tabular}}}}%
    \put(1.00149341,0.39016737){\color[rgb]{0,0,0}\makebox(0,0)[lt]{\lineheight{1.25}\smash{\begin{tabular}[t]{l}$\beta_i + \pi$\end{tabular}}}}%
    \put(1.00149341,0.62170119){\color[rgb]{0,0,0}\makebox(0,0)[lt]{\lineheight{1.25}\smash{\begin{tabular}[t]{l}$2\pi$\end{tabular}}}}%
  \end{picture}%
\endgroup%

    \caption{Hyperboles $H_{i}$ and $H_{i}^{\pi}$ together at the cylinder $C$
      forming the arctangents $A_{i}$ and $A_{i}^{\pi}$.}
    \label{sec:affine:sub:a2_nonzero:subfig:cylinder}
  \end{subfigure}
  \caption{Affine double discontinuity dynamics for $a_{i1} = 1$, $d_{i1} = -1$,
    $a_{i2} = 1$, $d_{i2} = 1$, $a_{i3} = 1$ and $d_{i3} = 0$. At this example
    we have $\alpha_i = -1$, $\beta_i = \frac{\pi}{4}$ and $\delta_i = 1$.
    Therefore, for example, $S_1$ has part of the hyperbole $H_{i}$ as a visible
    part of the slow manifold; whereas $S_2$ has only part of $A_{i}^{\pi}$
    visible.}
  \label{sec:affine:sub:a2_nonzero:fig:cylinder}
\end{figure}
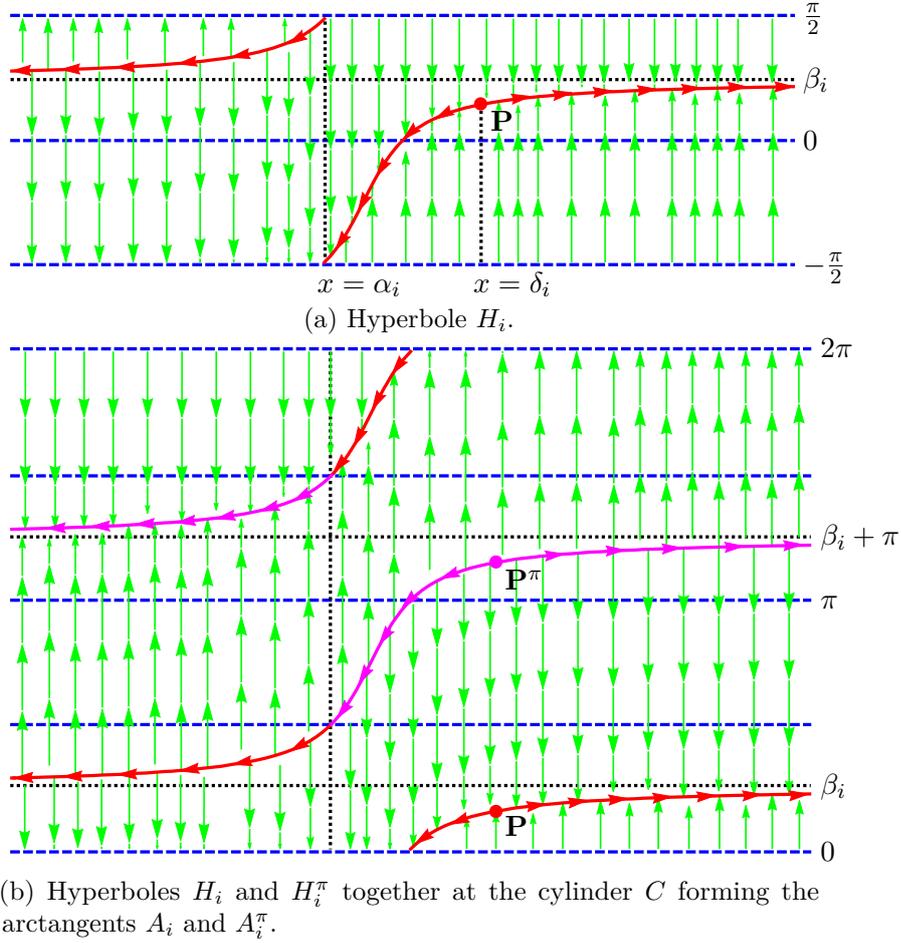

\noindent More precisely, the hyperbole $H_{i}$ behave as the red part of
\cref{sec:affine:sub:a2_nonzero:subfig:hyperbole}. However, putting together the
hyperboles $H_{i}$ and $H_{i}^{\pi}$ we get that they actually behave as two
arctangent-like curves as represented at
\cref{sec:affine:sub:a2_nonzero:subfig:cylinder}.

These arctangent-like curves will be denoted by $A_{i}$ and $A_{i}^{\pi}$. Based
on the analysis done before, we conclude that they are given by

\begin{align*}
  A_{i}       & = \left\{(x,\theta) \in [-\infty,\alpha_i] \times
  \left[0,2\pi\right];~
  \theta = \theta_i(x) + \pi \right\} \cup                                                 \\
              & \cup \left\{(x,\theta) \in [\alpha_i,+\infty] \times \left[0,2\pi\right];~
  \theta = \theta_i(x) \right\},                                                           \\
  A_{i}^{\pi} & = \left\{(x,\theta) \in [-\infty,\alpha_i] \times \left[0,2\pi\right];~
  \theta = \theta_i(x) \right\} \cup                                                       \\
              & \cup \left\{(x,\theta) \in [\alpha_i,+\infty] \times \left[0,2\pi\right];~
  \theta = \theta_i(x) + \pi \right\},
\end{align*}

\noindent and, therefore, on one hand, $A_{i}$ is an arctangent-like curve with
$\theta = \beta_i + \pi$ and $\theta = \beta_i$ as negative and
positive\footnote{Where negative means $x \to -\infty$ and positive means $x \to
    +\infty$.} horizontal asymptotes, respectively; on the other hand, $A_{i}^{\pi}$
is an arctangent-like curve with $\theta = \beta_i$ and $\theta = \beta_i + \pi$
as negative and positive horizontal asymptotes, respectively.\footnote{In
  particular, when $a_{i3} = 0$ we have $\beta_i = 0$ and, therefore, the
  horizontal asymptotes are given by $\theta = 0$ and $\theta = \pi$, which are
  part of the stripes' boundary.} Moreover, because of the very definition of
$\beta_i$, the positioning of the asymptotes inside the cylinder behaves
similarly as the straight lines $L_{i}$ and $L_{i}^{\pi}$ in
\cref{sec:constant}. This completes the qualitative analysis of the shape of the
slow manifold and, from now on we will write $\mathcal{M}_i = A_{i} \cup
  A_{i}^{\pi}$.

Over both the arctangents $\mathcal{M}_i = A_{i} \cup A_{i}^{\pi}$, we have the
one-dimensional dynamics given by the first equation of
\eqref{sec:affine:eq:reduced}, i.e., $\dot{x} = a_{i1}x + d_{i1}$. Analyzing
this equation we observe that, if $a_{i1} \neq 0$, then there are hyperbolic
critical points at

\begin{equation*}
  x = -\frac{d_{i1}}{a_{i1}} \eqqcolon \delta_i,
\end{equation*}

\noindent being these points attractors if $a_{i1} < 0$ and repellers if $a_{i1}
  > 0$, as represented at \cref{sec:affine:sub:a2_nonzero:subfig:cylinder}.
Since we are under \cref{sec:framework:sub:blowup:eq:strong_fund_hypo}, then
\cref{sec:framework:sub:blowup:cor:singularities} tells us that, in this case,
these hyperbolic singularities are actually hyperbolic singularities of the
whole stripe $S_i$. If $a_{i1} = 0$, then there is no critical point and the
dynamics over $\mathcal{M}_i$ is exactly as in the constant case described in
\cref{sec:constant}. This completes the qualitative analysis of the reduced
dynamics.

Regarding the layer dynamics, we have the layer system
\eqref{sec:affine:eq:layer} which says that for each fixed value of $x \in
  \mathbb{R}$, we have a one-dimensional dynamics given by the second equation of
\eqref{sec:affine:eq:layer}. In particular, assuming that $\cos{\theta} > 0$ and
$a_{i2}x + d_{i2} > 0$, then

\begin{equation*}
  \theta' > 0 \Leftrightarrow \theta < \theta_i(x),
\end{equation*}

\noindent since the arctangent function is strictly increasing. Likewise, and
under the same conditions, we have that

\begin{equation*}
  \theta' < 0 \Leftrightarrow \theta > \theta_i(x),
\end{equation*}

\noindent and, therefore, we conclude that for $a_{i2}x + d_{i2} > 0$, the piece
of curve $\theta = \theta_i(x)$ is attractor of the surrounding dynamics and,
therefore, $\theta = \theta_i(x) + \pi$ is repellor. Moreover, if $a_{i2} > 0$,
then $a_{i2}x + d_{i2} > 0$ happens for $x > \alpha_i$; if $a_{i2} < 0$, then
$a_{i2}x + d_{i2} > 0$ happens for $x < \alpha_i$. Completing this analysis and
comparing with the definition of $A_{i}$ and $A_{i}^{\pi}$ we reach the results
summarized at \cref{sec:affine:tab:fast} and represented as the green part of
\cref{sec:affine:sub:a2_nonzero:fig:cylinder}. Moreover, at $\cos{\theta} = 0$
with $a_{i2} \neq 0$ and $a_{i2}x + d_{i2} \neq 0$, \eqref{sec:affine:eq:layer}
give us the layer systems

\begin{equation*}
  \left\{\begin{matrix*}[l]
    x' & = & 0 \\
    \theta' & = & - (a_{i2}x + d_{i2})
  \end{matrix*}\right.
  \quad \text{ and } \quad
  \left\{\begin{matrix*}[l]
    x' & = & 0 \\
    \theta' & = & a_{i2}x + d_{i2}
  \end{matrix*}\right.
\end{equation*}

\noindent for $\theta = \frac{\pi}{2}$ and $\theta = \frac{3\pi}{2}$,
respectively, whose dynamics is consistent with \cref{sec:affine:tab:fast}. This
completes the qualitative analysis of the layer dynamics for case (A).

Now, lets consider the case (B), which complements the case (A) studied above
defining the missing dynamics over $a_{i2}x + d_{i2} = 0$ ($\Leftrightarrow x =
  \alpha_i$) with $a_{i2} \neq 0$. At this case, the reduced system
\eqref{sec:affine:eq:reduced} becomes

\begin{equation*}
  \left\{\begin{matrix*}[l]
    \dot{x} & = & a_{i1}x + d_{i1} \\
    0 & = &
    (a_{i3}x + d_{i3}) \cos{\theta}
  \end{matrix*}\right.,
\end{equation*}

\noindent whose slow manifold $\mathcal{M}_i$ is implicitly given by the
equation $0 = (a_{i3}x + d_{i3}) \cos{\theta}$ which actually means $0 =
  \cos{\theta}$, since we are under SFH and, therefore $a_{i3}x + d_{i3} \neq 0$.
In other words, $\mathcal{M}_i = \left\{\left(\alpha_i, \frac{\pi}{2}\right),
  \left(\alpha_i, \frac{3\pi}{2}\right)\right\}$. Over these points acts the
dynamics $\dot{x} = a_{i1}x + d_{i1}$, which is consistent with case (A).
Regarding the fast dynamics, we have the layer system

\begin{equation*}
  \left\{\begin{matrix*}[l]
    x' & = & 0 \\
    \theta' & = & (a_{i3}x + d_{i3}) \cos{\theta}
  \end{matrix*}\right.
  \quad
  \sim
  \quad
  \left\{\begin{matrix*}[l]
    x' & = & 0 \\
    \theta' & = & -\frac{\gamma_i}{a_{i2}} \cos{\theta}
  \end{matrix*}\right.,
\end{equation*}

\noindent since $x = \alpha_i$, which can be easily verified to be consistent
with the layer dynamics given by \cref{sec:affine:tab:fast} and, therefore, it
is consistent with case (A). Therefore, we conclude that case whole (B) is
consistent with case (A). In other words, the dynamics over the asymptote
$a_{i2}x + d_{i2} = 0$ agrees with the surrounding dynamics.

\begin{table}[ht]
  \centering
  \caption{Layer dynamics around the arctangents $A_{i}$ and $A_{i}^{\pi}$ that
    compose the slow manifold $\mathcal{M}_i = A_{i} \cup A_{i}^{\pi}$.}
  \begin{tabular}{lll}
    \hline
                  & $a_{i2} < 0$ & $a_{i2} > 0$ \\
    \hline
    $A_{i}$       & repellor     & attractor    \\
    $A_{i}^{\pi}$ & attractor    & repellor     \\
    \hline
  \end{tabular}
  \label{sec:affine:tab:fast}
\end{table}

\subsection{Case \texorpdfstring{$a_{i2} = 0$}{a2=0}}
\label{sec:affine:sub:a2_zero}

\begin{figure}[ht]
  \centering
  \def\svgwidth{0.85\linewidth}
\begingroup%
  \makeatletter%
  \providecommand\color[2][]{%
    \errmessage{(Inkscape) Color is used for the text in Inkscape, but the package 'color.sty' is not loaded}%
    \renewcommand\color[2][]{}%
  }%
  \providecommand\transparent[1]{%
    \errmessage{(Inkscape) Transparency is used (non-zero) for the text in Inkscape, but the package 'transparent.sty' is not loaded}%
    \renewcommand\transparent[1]{}%
  }%
  \providecommand\rotatebox[2]{#2}%
  \newcommand*\fsize{\dimexpr\f@size pt\relax}%
  \newcommand*\lineheight[1]{\fontsize{\fsize}{#1\fsize}\selectfont}%
  \ifx\svgwidth\undefined%
    \setlength{\unitlength}{532.47088623bp}%
    \ifx\svgscale\undefined%
      \relax%
    \else%
      \setlength{\unitlength}{\unitlength * \real{\svgscale}}%
    \fi%
  \else%
    \setlength{\unitlength}{\svgwidth}%
  \fi%
  \global\let\svgwidth\undefined%
  \global\let\svgscale\undefined%
  \makeatother%
  \begin{picture}(1,0.64512054)%
    \lineheight{1}%
    \setlength\tabcolsep{0pt}%
    \put(0,0){\includegraphics[width=\unitlength,page=1]{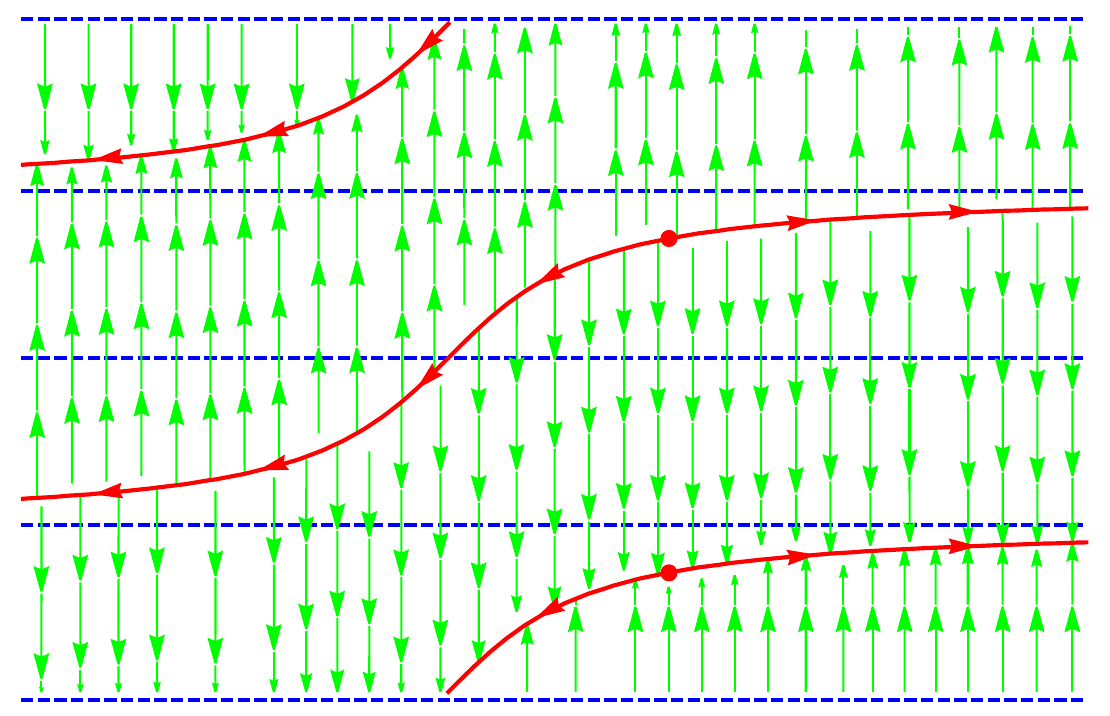}}%
    \put(0.99248405,0.00170737){\makebox(0,0)[lt]{\lineheight{1.25}\smash{\begin{tabular}[t]{l}$0$\end{tabular}}}}%
    \put(0.99248405,0.16509661){\makebox(0,0)[lt]{\lineheight{1.25}\smash{\begin{tabular}[t]{l}$\sigma_{i+}$\end{tabular}}}}%
    \put(0.99248405,0.31565502){\makebox(0,0)[lt]{\lineheight{1.25}\smash{\begin{tabular}[t]{l}$\pi$\end{tabular}}}}%
    \put(0.99248405,0.4669085){\makebox(0,0)[lt]{\lineheight{1.25}\smash{\begin{tabular}[t]{l}$\sigma_{i-}$\end{tabular}}}}%
    \put(0.99248405,0.61621244){\makebox(0,0)[lt]{\lineheight{1.25}\smash{\begin{tabular}[t]{l}$2\pi$\end{tabular}}}}%
    \put(0.61228921,0.3974192){\makebox(0,0)[lt]{\lineheight{1.25}\smash{\begin{tabular}[t]{l}$\mathbf{P}^{\pi}$\end{tabular}}}}%
    \put(0.61228921,0.0988113){\makebox(0,0)[lt]{\lineheight{1.25}\smash{\begin{tabular}[t]{l}$\mathbf{P}$\end{tabular}}}}%
  \end{picture}%
\endgroup%

  \caption{Affine double discontinuity dynamics for $a_{i1} = 1$, $d_{i1} = -1$,
    $a_{i2} = 0$, $d_{i2} = 1$, $a_{i3} = 1$ and $d_{i3} = 1$. At this example
    we have $\delta_i = 1$ and $\sigma_{i\pm} = \pm\frac{\pi}{2}$.}
  \label{sec:affine:sub:a2_zero:fig:cylinder}
\end{figure}
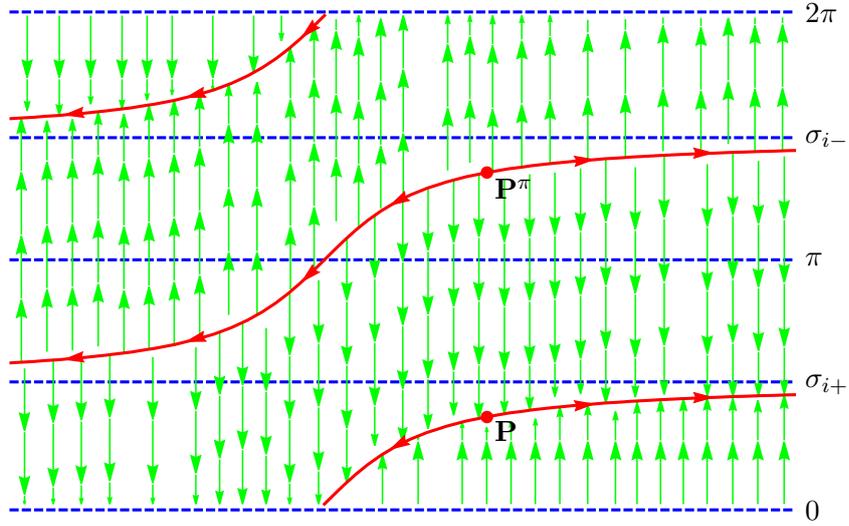

For case (C), remember that we have $a_{i2} = 0$ and $a_{i2}x + d_{i2} \neq 0$
implying $d_{i2} \neq 0$. Therefore, everything at the beginning of
\cref{sec:affine:sub:a2_nonzero} is true. However, whenever $\cos{\theta} \neq
  0$, the explicit expression for the slow manifold $\mathcal{M}_i$ is now

\begin{equation*}
  \theta = \arctan{\left(\frac{a_{i3}x + d_{i3}}{d_{i2}}\right)}
  + n\pi = \theta_i(x) + n\pi,
\end{equation*}

\noindent where $n \in \mathbb{Z}$. Therefore, without loss of generality, the
slow manifold can be written as $\mathcal{M}_i = A_{i} \cup A_{i}^{\pi}$, where

\begin{align*}
  A_{i}       & = \left\{(x,\theta) \in \mathbb{R} \times \left[0,2\pi\right] ;~
  \theta = \theta_i(x) \right\} \text{ and}                                      \\
  A_{i}^{\pi} & = \left\{(x,\theta) \in \mathbb{R} \times \left[0,2\pi\right] ;~
  \theta = \theta_i(x) + \pi \right\},
\end{align*}

\noindent which consists of two arctangent-like curves inside the cylinder $C
  = \mathbb{R} \times S^1$ as the red part of
\cref{sec:affine:sub:a2_zero:fig:cylinder}. In fact, since $a_{i2} = 0$, then
$h(x)$ is a straight line and, therefore,

\begin{equation*}
  \theta = \theta_i(x) = \arctan{\left( h(x) \right)}
\end{equation*}

\noindent is an arctangent curve. Besides that, we have

\begin{equation*}
  \frac{d}{dx}h(x) =
  \frac{\gamma_i}{(a_{i2}x + d_{i2})^2} =
  \frac{\gamma_i}{d_{i2}^{2}}
\end{equation*}

\noindent and, therefore, $A_{i}$ and $A_{i}^{\pi}$ are increasing curves if
$\gamma_i > 0$ and decreasing if $\gamma_i < 0$\footnote{Again, if $\gamma_{i} =
    0$ ($\Leftrightarrow a_{i3} = 0$), then $h(x)$ is a constant function and,
  therefore, $A_{i}$ and $A_{i}^{\pi}$ are straight lines. In other words, the
  constant case is recovered.}. Moreover, since

\begin{align*}
  \lim_{x \to \pm\infty} \theta_i(x)
   & = \lim_{x \to \pm\infty} \left[ \arctan{\left(\frac{a_{i3}x +
  d_{i3}}{d_{i2}} \right)}\right] =                                          \\
   & = \arctan{\left[\lim_{x \to \pm\infty} \left(\frac{a_{i3}x + d_{i3}}
  {d_{i2}} \right) \right]} =                                                \\
   & = \arctan{\left[\pm\sgn{\left(\frac{a_{i3}}{d_{i2}} \right)}\infty
  \right]} =                                                                 \\
   & = \pm\sgn{\left(\frac{a_{i3}}{d_{i2}} \right)}\frac{\pi}{2}
  =  \pm\sgn{\left(\frac{\gamma_i}{d_{i2}^{2}} \right)}\frac{\pi}{2} =       \\
   & = \pm\sgn{\left(\gamma_i \right)}\frac{\pi}{2} \eqqcolon \sigma_{i\pm},
\end{align*}

\noindent then $A_{i}$ has $\sigma_{i-}$ and $\sigma_{i+}$ as negative and
positive horizontal asymptote, respectively; while $A_{i}^{\pi}$ has
$\sigma_{i+}$ and $\sigma_{i-}$ as negative and positive horizontal asymptote,
respectively. This completes the qualitative analysis of the shape of the slow
manifold.

Over both the arctangents $\mathcal{M}_i = A_{i} \cup A_{i}^{\pi}$, we have the
one-dimensional dynamics given by the first equation of
\eqref{sec:affine:eq:reduced}, i.e., $\dot{x} = a_{i1}x + d_{i1}$ which behaves
as described in \cref{sec:affine:sub:a2_nonzero}. This completes the qualitative
analysis of the reduced dynamics.

Regarding the layer dynamics, a completely analogous analysis such as that made
for the previous cases allows us to conclude that it behaves as described in
\cref{sec:affine:tab:fast}, including the case $\cos{\theta} = 0$, but
exchanging $a_{i2}$ with $d_{i2}$.

\subsection{Theorem and Examples}
\label{sec:affine:sub:theorem}

Summarizing, we conclude that the dynamics over $\Sigma_x$ for affine fields
behaves as described in the theorem below, whose proof consists in the analysis
done above.

\TheoremAffineDynamics

\begin{example}
  \label{sec:affine:exmp:slow_cycle}
  Let $\mathbf{F} \in \mathcal{A}_3$ be given by affine vector fields such that

  \begin{align*}
     & \mathbf{F}_2 :
    \begin{bmatrix}
      a_{21} & d_{21} \\
      a_{22} & d_{22} \\
      a_{23} & d_{23}
    \end{bmatrix}
    =
    \begin{bmatrix*}[r]
      1 & -2 \\
      -1  & 1 \\
      -1  & 0
    \end{bmatrix*},
    \quad
     & \mathbf{F}_1 :
    \begin{bmatrix}
      a_{11} & d_{11} \\
      a_{12} & d_{12} \\
      a_{13} & d_{13}
    \end{bmatrix}
    =
    \begin{bmatrix*}[r]
      -1 & 2 \\
      -1  & 1 \\
      1  & 0
    \end{bmatrix*}, \\
    \\
     & \mathbf{F}_3 :
    \begin{bmatrix}
      a_{31} & d_{31} \\
      a_{32} & d_{32} \\
      a_{33} & d_{33}
    \end{bmatrix}
    =
    \begin{bmatrix*}[r]
      1 & -2 \\
      1  & 1 \\
      -1  & 0
    \end{bmatrix*},
    \quad
     & \mathbf{F}_4 :
    \begin{bmatrix}
      a_{41} & d_{41} \\
      a_{42} & d_{42} \\
      a_{43} & d_{43}
    \end{bmatrix}
    =
    \begin{bmatrix*}[r]
      -1 & 2 \\
      1  & 1 \\
      1  & 0
    \end{bmatrix*},
  \end{align*}

  \noindent with parameters $c_{ij}$'s and $d_{ij}$'s arbitrary since, according
  to \cref{sec:affine:thm:dynamics}, they only affect the dynamics outside the
  cylinder. Using this theorem we can also verify that, over the cylinder $C$
  given by the blow-up of $\Sigma_x$, the system has a single \textbf{slow
    cycle} as represented at \cref{sec:affine:sub:theorem:fig:affine_slow_cycle}.

  For instance, according to \cref{sec:affine:thm:dynamics}, the field
  $\mathbf{F}_1$ induces a slow-fast system whose slow manifold $\mathcal{M}_1 =
    A_1 \cup A_1^{\pi}$ consists of arctangents with horizontal asymptotes

  \begin{equation*}
    \theta = \beta_1 = \arctan{\left( \frac{a_{13}}{a_{12}} \right)}
    = \arctan{(-1)} = -\frac{\pi}{4}
  \end{equation*}

  \noindent at $S_4$ and $\theta = \beta_1 + \pi = \frac{3\pi}{4}$ at $S_2$.
  Besides that, since

  \begin{equation*}
    \gamma_1 = a_{13}d_{12} - a_{12}d_{13} = 1
  \end{equation*}

  \noindent then these arctangents are increasing. Therefore, we conclude that
  $\mathcal{M}_1 \cap S_1 \subset A_1^{\pi}$, and it transversally crosses $S_1$
  as represented at the lowest stripe of
  \cref{sec:affine:sub:theorem:fig:affine_slow_cycle} from $\mathbf{R}_1$ to
  $\mathbf{Q}_1$, where the point $\mathbf{Q}_1$ is given by

  \begin{equation*}
    \frac{\pi}{2} = \theta_1(x) =
    \arctan{\left(\frac{a_{i3}x + d_{i3}}{a_{i2}x + d_{i2}}\right)} =
    \arctan{\left(\frac{x}{-x+1} \right)},
  \end{equation*}

  \noindent which happens when $x \to 1^-$; and the point $\mathbf{R}_1$ is given
  by

  \begin{equation*}
    0 = \theta_1(x) =
    \arctan{\left(\frac{a_{i3}x + d_{i3}}{a_{i2}x + d_{i2}}\right)} =
    \arctan{\left(\frac{x}{-x+1} \right)},
  \end{equation*}

  \noindent which happens when $x \to 0^+$. Dynamically it also goes
  $\mathbf{R}_1 \to \mathbf{Q}_1$, since over $\mathcal{M}_1 \cap S_1$ acts the
  reduced dynamics $\dot{x} = -x+2$, which has $x = 2$ as a stable singularity.
  Finally, since $a_{12} = -1 < 0$ and $\mathcal{M}_1 \cap S_1 \subset
    A_1^{\pi}$, then $\mathcal{M}_1 \cap S_1$ attracts the surrounding layer
  dynamics, according to \cref{sec:affine:tab:fast}.

  \begin{figure}[ht]
    \centering
    \def\svgwidth{0.85\linewidth}
\begingroup%
  \makeatletter%
  \providecommand\color[2][]{%
    \errmessage{(Inkscape) Color is used for the text in Inkscape, but the package 'color.sty' is not loaded}%
    \renewcommand\color[2][]{}%
  }%
  \providecommand\transparent[1]{%
    \errmessage{(Inkscape) Transparency is used (non-zero) for the text in Inkscape, but the package 'transparent.sty' is not loaded}%
    \renewcommand\transparent[1]{}%
  }%
  \providecommand\rotatebox[2]{#2}%
  \newcommand*\fsize{\dimexpr\f@size pt\relax}%
  \newcommand*\lineheight[1]{\fontsize{\fsize}{#1\fsize}\selectfont}%
  \ifx\svgwidth\undefined%
    \setlength{\unitlength}{521.375bp}%
    \ifx\svgscale\undefined%
      \relax%
    \else%
      \setlength{\unitlength}{\unitlength * \real{\svgscale}}%
    \fi%
  \else%
    \setlength{\unitlength}{\svgwidth}%
  \fi%
  \global\let\svgwidth\undefined%
  \global\let\svgscale\undefined%
  \makeatother%
  \begin{picture}(1,0.6396698)%
    \lineheight{1}%
    \setlength\tabcolsep{0pt}%
    \put(0,0){\includegraphics[width=\unitlength,page=1]{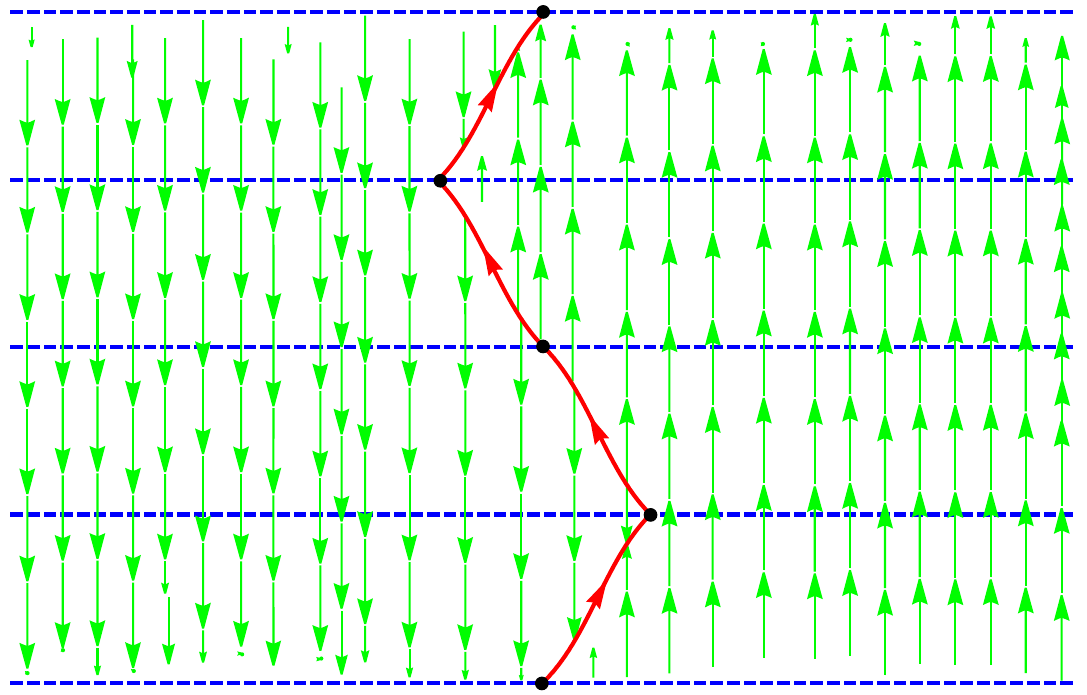}}%
    \put(0.59981922,0.13280372){\makebox(0,0)[lt]{\lineheight{1.25}\smash{\begin{tabular}[t]{l}$\mathbf{Q}_1$\end{tabular}}}}%
    \put(0.45411456,0.02876955){\makebox(0,0)[lt]{\lineheight{1.25}\smash{\begin{tabular}[t]{l}$\mathbf{R}_1$\end{tabular}}}}%
    \put(0.9990719,0.00078487){\makebox(0,0)[lt]{\lineheight{1.25}\smash{\begin{tabular}[t]{l}$0$\end{tabular}}}}%
    \put(0.9990719,0.31060828){\makebox(0,0)[lt]{\lineheight{1.25}\smash{\begin{tabular}[t]{l}$\pi$\end{tabular}}}}%
    \put(0.9990719,0.61973676){\makebox(0,0)[lt]{\lineheight{1.25}\smash{\begin{tabular}[t]{l}$2\pi$\end{tabular}}}}%
  \end{picture}%
\endgroup%

    \caption{Dynamics over $C$ generated by the field $\mathbf{F}$ studied at
      \cref{sec:affine:exmp:slow_cycle}.}
    \label{sec:affine:sub:theorem:fig:affine_slow_cycle}
  \end{figure}
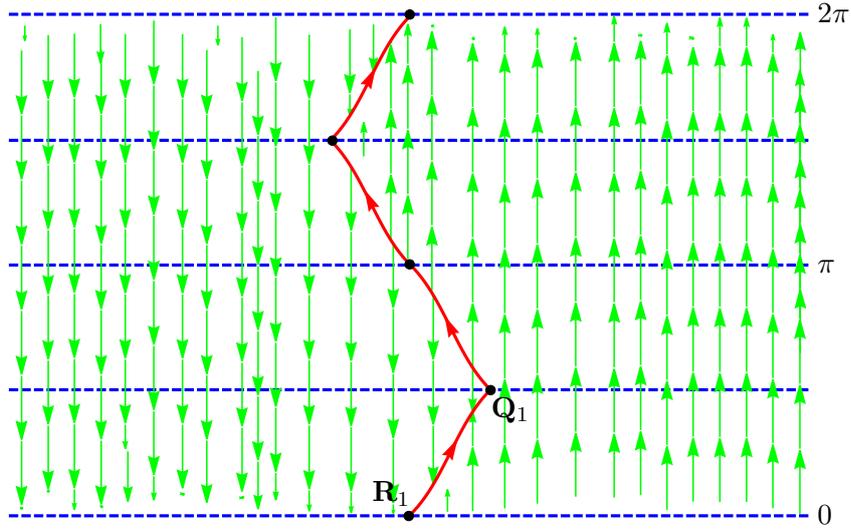

  Therefore, we conclude that the dynamics generated by $\mathbf{F}_1$ over the
  stripe $S_1$, in fact, behave as represented at
  \cref{sec:affine:sub:theorem:fig:affine_slow_cycle}. The dynamics over the
  other stripes can be similarly verified to be as represented. \qed
\end{example}

\begin{corollary}
  Every $\mathbf{F} \in \mathcal{A}_3$ with $\gamma_i \neq 0$ can induce at most
  one slow cycle over the cylinder.
\end{corollary}

\begin{proof}
  Given a stripe $S_i$, according to \cref{sec:affine:thm:dynamics} the
  arctangents that forms the slow manifold $\mathcal{M}_i$ can either have a
  horizontal asymptote inside $S_i$ or not.

  If a horizontal asymptote is inside $S_i$, then a slow cycle construction is
  impossible, even if the asymptote is at one of the borders of $S_i$, since
  $\mathcal{M}_i$ does not cross transversally both borders of $S_i$.

  However, if no horizontal asymptote is inside $S_i$, then a construction
  similar to that realized at \cref{sec:affine:exmp:slow_cycle} can occur.
  Finally, no more than one slow cycle can occur, since the arctangents are
  strictly monotonous and, therefore, transversally crosses $S_i$ at most once.
\end{proof}

\section{Structural Stability}%
\label{sec:stability}

Let $\mathbf{F} \in \mathcal{D}_3^k$ be a piecewise smooth vector field with a
double discontinuity given by affine vector fields \eqref{sec:affine:eq:system}.
The theorems obtained in the previous sections fully describe the constant and
affine fundamental dynamics over the cylinder $C$ of the induced vector field
$\mathbf{\tilde{F}} \in \tilde{\mathcal{D}}_3^k$. As an application, we would
like to leverage this knowledge to study its structural stability. The first
step in this process consists of defining a concept of structural stability
which fits the systems we are studying. In order to do so, we are going to mimic
the classical definition for the regular case, $\mathcal{R}^k(U)$, presented in
\cite[p.~18]{Teixeira1990}, which can be easily extended to $\mathcal{D}_3^k$.
In fact, on one hand, systems in $\mathcal{R}^k(U)$ have a single subset that
should be kept invariant, $\Sigma = h^{-1}(0)$; on the other hand, systems in
$\mathcal{D}_3^k$ have a set of subsets

\begin{equation*}
  \mathcal{I} = \left\{
  \Sigma_{12}, \Sigma_{23}, \Sigma_{34}, \Sigma_{14}, \Sigma_{x}
  \right\}
\end{equation*}

\noindent which should be kept invariants by topological equivalence. Therefore,
a direct substitution gives us the following definition:

\begin{definition}
  \label{sec:stability:defn:equivalence_D}
  Let $\mathbf{F}, \mathbf{G} \in \mathcal{D}_3^k$. We say that $\mathbf{F}$ and
  $\mathbf{G}$ are \textbf{topologically equivalent} and denote $\mathbf{F} \sim
    \mathbf{G}$ if, and only if, there exists a homeomorphism $\varphi:
    \mathbb{R}^3 \to \mathbb{R}^3$ \noindent that keeps every $I \in \mathcal{I}$
  invariant and takes orbits of $\mathbf{F}$ into orbits of $\mathbf{G}$
  preserving the orientation of time. From this definition the concept of
  structural stability in $\mathcal{D}_3^k$ is naturally obtained.
\end{definition}

For the blow-up induced vector fields, $\tilde{\mathcal{D}}_3^k$, the set of
invariant subsets are given by

\begin{equation*}
  \tilde{\mathcal{I}} = \left\{
  \tilde{\Sigma}_{12}, \tilde{\Sigma}_{23}, \tilde{\Sigma}_{34},
  \tilde{\Sigma}_{14}, C  \right\}
\end{equation*}

\noindent and, therefore, we define:

\begin{definition}
  \label{sec:stability:defn:equivalence_D_ind}
  Let $\mathbf{\tilde{F}}, \mathbf{\tilde{G}} \in \tilde{\mathcal{D}}_3^k$. We say
  that $\mathbf{\tilde{F}}$ and $\mathbf{\tilde{G}}$ are \textbf{topologically
    equivalent} and denote $\mathbf{\tilde{F}} \sim \mathbf{\tilde{G}}$ if, and
  only if, there exists a homeomorphism $\tilde{\varphi}: \mathbb{R} \times S^1
    \times \mathbb{R}^+ \to \mathbb{R} \times S^1 \times \mathbb{R}^+$ that keeps
  every $I \in \tilde{\mathcal{I}}$ invariant and takes orbits of
  $\mathbf{\tilde{F}}$ into orbits of $\mathbf{\tilde{G}}$ preserving the
  orientation of time. From this definition the concept of structural stability
  in $\tilde{\mathcal{D}}_3^k$ is naturally obtained.
\end{definition}

Now, let $\mathbf{\tilde{F}}, \mathbf{\tilde{G}} \in \tilde{\mathcal{D}}_3^k$ be
topologically equivalent by a homeomorphism $\tilde{\varphi}$. In this case, we
have that $\restr{\tilde{\varphi}}{I}$ with $I \in \tilde{\mathcal{I}}$ are also
homeomorphisms taking orbits into orbits and preserving the orientation of time.
In other words, the existence of these homeomorphisms is a necessary condition
for the topological equivalence. More precisely:

\begin{proposition}
  \label{sec:stability:prop:equivalence_invariants}
  If $\mathbf{\tilde{F}} \sim \mathbf{\tilde{G}}$, then
  $\restr{\mathbf{\tilde{F}}}{I} \sim \restr{\mathbf{\tilde{G}}}{I}$ for every
  $I \in \tilde{\mathcal{I}}$.
\end{proposition}

We are interested on the dynamics over the cylinder $C$. Therefore, given
$\mathbf{F} \in \mathcal{D}_3^k$, we look for necessary and/or sufficient
conditions for the structural stability of $\restr{\mathbf{\tilde{F}}}{C}$.
Beyond the intrinsic interest, given in
\cref{sec:stability:prop:equivalence_invariants} above, such conditions shall
also reveal relevant information on the structural stability of
$\mathbf{\tilde{F}}$ and, therefore, on the structural stability of
$\mathbf{F}$. In fact, from \cref{sec:stability:prop:equivalence_invariants} it
follows the result below.

\begin{corollary}
  \label{sec:stability:cor:stability_invariants}
  If $\mathbf{\tilde{F}}$ is structurally stable, then
  $\restr{\mathbf{\tilde{F}}}{I}$ is structurally stable for every $I \in
    \tilde{\mathcal{I}}$.
\end{corollary}

\begin{proof}
  Given $I \in \tilde{\mathcal{I}}$, let $\tilde{\mathcal{W}} \subset
    \mathcal{D}_3^k$ be an open neighborhood of $\mathbf{\tilde{F}}$. Observe that

  \begin{equation*}
    \restr{\tilde{\mathcal{W}}}{I} =
    \left\{ \restr{\mathbf{\tilde{H}}}{I} ;~
    \mathbf{\tilde{H}} \in \tilde{\mathcal{W}} \right\}
  \end{equation*}

  \noindent is an open neighborhood of $\restr{\mathbf{\tilde{F}}}{I}$.

  Therefore, if $\restr{\mathbf{\tilde{F}}}{I}$ was not structurally stable,
  would exist $\restr{\mathbf{\tilde{G}}}{I} \in \restr{\tilde{\mathcal{W}}}{I}$
  such that $\restr{\mathbf{\tilde{F}}}{I} \nsim \restr{\mathbf{\tilde{G}}}{I}$
  and, therefore, from \cref{sec:stability:prop:equivalence_invariants} would
  follow that $\mathbf{\tilde{G}} \nsim \mathbf{\tilde{F}}$, then implying that
  $\mathbf{\tilde{F}}$ would not be structurally stable.
\end{proof}

Thus, from now on, we will exclusively study conditions for the structural
stability of $\restr{\mathbf{\tilde{F}}}{C}$. In order to do so, remember that
over $C$ acts a regular Filippov dynamics whose switching manifold is formed by
the elements of

\begin{equation*}
  \tilde{\mathcal{I}}_C = \left\{
  \Sigma_{0}, \Sigma_{\frac{\pi}{2}}, \Sigma_{\pi}, \Sigma_{\frac{3\pi}{2}} \right\},
\end{equation*}

\noindent where $\Sigma_{\theta} = \left\{ (x,\theta) ;~ x \in \mathbb{R}
  \right\}$. Therefore, without loss of generality for the previous results, it
is natural to adopt the following definitions of equivalence and stability for
$C$:

\begin{definition}
  \label{sec:stability:defn:equivalence_cyl}
  Let $\mathbf{\tilde{F}}, \mathbf{\tilde{G}} \in \tilde{\mathcal{D}}_3^k$. We say
  that $\mathbf{\tilde{F}}$ and $\mathbf{\tilde{G}}$ are
  \textbf{$C$-topologically equivalent} and denote $\mathbf{\tilde{F}} \sim_{c}
    \mathbf{\tilde{G}}$ if, and only if, there exists a homeomorphism
  $\tilde{\varphi}: C \to C$ that keeps every $I \in \tilde{\mathcal{I}}_C$
  invariant and takes orbits of $\restr{\mathbf{\tilde{F}}}{C}$ into orbits of
  $\restr{\mathbf{\tilde{G}}}{C}$ preserving the orientation of time. From this
  definition the concept of $C$-structural stability is naturally obtained.
\end{definition}

Although global and naturally derived from the regular case, $C$-structural
stability, as presented above, is still a fairly complex property to prove and, in
fact, to the best of the author's knowledge, it is an open problem to characterize
it through simple conditions and, therefore, shall be treated in future works.

However, many of the difficulties found at characterizing $C$-structural
stability comes from its global aspect. In fact, conditions for a semi-local
approach can be found in \cite{Broucke2001} and, in order to apply these
results, a \emph{regular} and \emph{compact} Filippov section of the cylinder
$C$ must be taken.

\begin{figure}[ht]
  \centering
  \def\svgwidth{0.85\linewidth}
\begingroup%
  \makeatletter%
  \providecommand\color[2][]{%
    \errmessage{(Inkscape) Color is used for the text in Inkscape, but the package 'color.sty' is not loaded}%
    \renewcommand\color[2][]{}%
  }%
  \providecommand\transparent[1]{%
    \errmessage{(Inkscape) Transparency is used (non-zero) for the text in Inkscape, but the package 'transparent.sty' is not loaded}%
    \renewcommand\transparent[1]{}%
  }%
  \providecommand\rotatebox[2]{#2}%
  \newcommand*\fsize{\dimexpr\f@size pt\relax}%
  \newcommand*\lineheight[1]{\fontsize{\fsize}{#1\fsize}\selectfont}%
  \ifx\svgwidth\undefined%
    \setlength{\unitlength}{399.94170422bp}%
    \ifx\svgscale\undefined%
      \relax%
    \else%
      \setlength{\unitlength}{\unitlength * \real{\svgscale}}%
    \fi%
  \else%
    \setlength{\unitlength}{\svgwidth}%
  \fi%
  \global\let\svgwidth\undefined%
  \global\let\svgscale\undefined%
  \makeatother%
  \begin{picture}(1,0.50127929)%
    \lineheight{1}%
    \setlength\tabcolsep{0pt}%
    \put(0,0){\includegraphics[width=\unitlength,page=1]{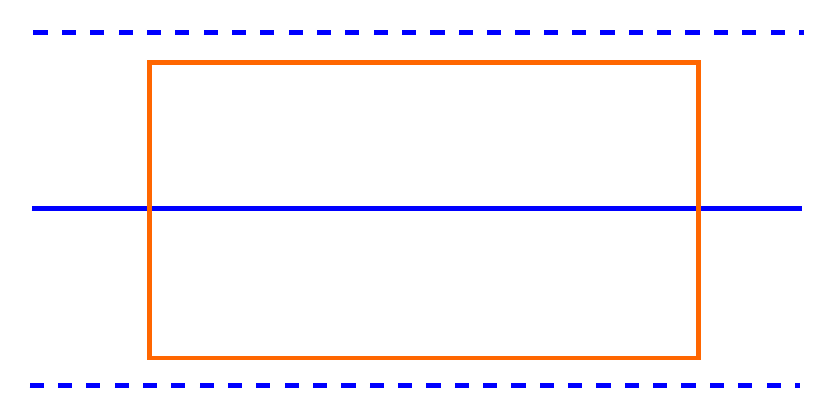}}%
    \put(0.04077242,0.33401711){\color[rgb]{0,0,0}\makebox(0,0)[lt]{\lineheight{1.25}\smash{\begin{tabular}[t]{l}$S_+$\end{tabular}}}}%
    \put(0.04077242,0.14820714){\color[rgb]{0,0,0}\makebox(0,0)[lt]{\lineheight{1.25}\smash{\begin{tabular}[t]{l}$S_-$\end{tabular}}}}%
    \put(0.45954528,0.31718465){\color[rgb]{0,0,0}\makebox(0,0)[lt]{\lineheight{1.25}\smash{\begin{tabular}[t]{l}$X_+$\end{tabular}}}}%
    \put(0.45954528,0.16503959){\color[rgb]{0,0,0}\makebox(0,0)[lt]{\lineheight{1.25}\smash{\begin{tabular}[t]{l}$X_-$\end{tabular}}}}%
    \put(0.98168396,0.24277963){\makebox(0,0)[lt]{\lineheight{1.25}\smash{\begin{tabular}[t]{l}$\Sigma_{\theta_0}$\end{tabular}}}}%
    \put(0.79425292,0.38841683){\makebox(0,0)[lt]{\lineheight{1.25}\smash{\begin{tabular}[t]{l}$K$\end{tabular}}}}%
  \end{picture}%
\endgroup%

  \caption{Regular Filippov system $\mathbf{X} = (\mathbf{X}_{-},
      \mathbf{X}_{+})$ defined at a convex compact set $K \subset C_{+} \cup
      C_{-}$ with switching manifold $\Sigma_{\theta_0}$.}
  \label{sec:stability:fig:compact}
\end{figure}
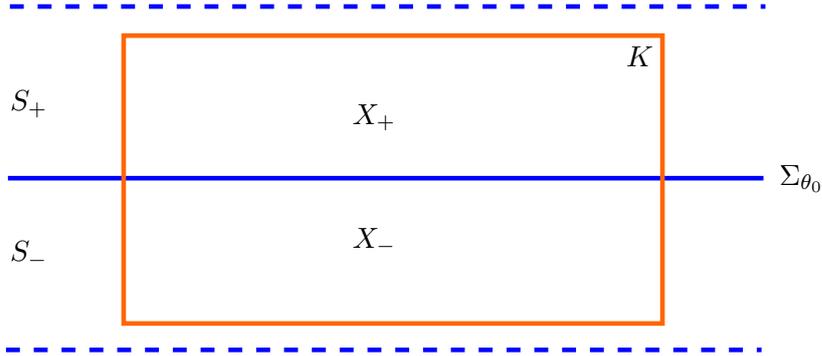

More precisely, given $\mathbf{F} \in \mathcal{D}_3^k$ and two consecutive
stripes $C_{+}$ and $C_{-}$ meeting at a straight line $\Sigma_{\theta_0} \in
  \tilde{\mathcal{I}}_C$, let $\mathbf{X}_{+}$ and $\mathbf{X}_{-}$ be the smooth
vector fields induced over $C_{+} \cap K$ and $C_{-} \cap K$, respectively, as
described at the previous sections and where $K \subset C_{+} \cup C_{-}$ is a
convex compact set, see \cref{sec:stability:fig:compact}. Observe that
$\mathbf{X} = (\mathbf{X}_{-}, \mathbf{X}_{+})$ is a \emph{regular} and
\emph{compact} Filippov system with a connected (due to the convexity of $K$)
switching manifold $\Sigma_{\theta_0}$.

Then, a direct application of Theorem B from \cite[p.~5]{Broucke2001} and the
Proposition 1.1 from \cite[p.~122]{Palis1982} give us the following result:

\begin{proposition}
  \label{sec:stability:prop:conditions}
  Given $\mathbf{F} \in \mathcal{D}_3^k$, two consecutive stripes $C_{+}$ and
  $C_{-}$ and a convex compact set $K \subset C_{+} \cup C_{-}$, then the
  induced Filippov system $\mathbf{X} = (\mathbf{X}_{-}, \mathbf{X}_{+})$ is
  structurally stable inside $K$ if, and only if, the following sets of
  conditions are satisfied:

  \begin{enumerate}[label=(\Roman*)]
    \item\label[condition]{sec:stability:conditions:ms} $\mathbf{X}_{+}$ and
          $\mathbf{X}_{-}$ are robustly\footnote{In other words, the property is
            stable under small perturbations.} Morse-Smale, i.e., they have:

          \begin{enumerate}[label=(C.\arabic*)]
            \item\label[condition]{sec:stability:conditions:ms_hyperbolic}
                  finitely many critical elements\footnote{Singularities and
                    periodic orbits.}, all hyperbolic;

            \item\label[condition]{sec:stability:conditions:ms_separatrix} no
                  saddle-connections;

            \item\label[condition]{sec:stability:conditions:ms_wandering} only
                  critical elements as non-wandering points;
          \end{enumerate}

    \item\label[condition]{sec:stability:conditions:sg} $\mathbf{X}_{+}$ and
          $\mathbf{X}_{-}$ robustly satisfies that:

          \begin{enumerate}[resume, label=(C.\arabic*)]
            \item\label[condition]{sec:stability:conditions:sg_zeros} none of
                  them vanishes at a point of $\Sigma_{\theta_0}$;

            \item\label[condition]{sec:stability:conditions:sg_tangency} they
                  are tangent to $\Sigma_{\theta_0}$ at only finitely many
                  points with both never tangent at the same point;

            \item\label[condition]{sec:stability:conditions:sg_colinear} they
                  are colinear at only finitely many points;
          \end{enumerate}

    \item\label[condition]{sec:stability:conditions:cr} $\mathbf{X}$ have:

          \begin{enumerate}[resume, label=(C.\arabic*)]
            \item\label[condition]{sec:stability:conditions:cr_periodic} only
                  hyperbolic periodic orbits;

            \item\label[condition]{sec:stability:conditions:cr_separatrix} no
                  separatrix-connections or relations\footnote{Unstable
                    separatrices arriving at the same point are said to be
                    \textbf{related}.};

            \item\label[condition]{sec:stability:conditions:cr_recurrent} only
                  trivial recurrent orbits.
          \end{enumerate}
  \end{enumerate}
\end{proposition}

Observe that \cref{sec:stability:conditions:ms} refers only to the usual
dynamics of $\mathbf{X}_{+}$ and $\mathbf{X}_{-}$ over the smooth parts. On the
other hand, \cref{sec:stability:conditions:sg} considers only the values of
$\mathbf{X}_{+}$ and $\mathbf{X}_{-}$ over the switching manifold
$\Sigma_{\theta_0}$. Finally, only \cref{sec:stability:conditions:cr} refers to
the actual Filippov dynamics of $\mathbf{X}$. With that in mind, over the next,
and final sections, we will apply \cref{sec:constant:thm:dynamics} and
\cref{sec:affine:thm:dynamics} to analyze this conditions for the particular
cases of constant and affine double discontinuities, respectively, and therefore
derive semi-local structural stability theorems or, more precisely:

\begin{definition}
  \label{sec:stability:defn:semilocal-stability}
  We say that $\mathbf{F} \in \mathcal{D}_3^k$ is \textbf{$(I, K)$-semi-local
    structurally stable} if, and only if, the induced Filippov system
  $\mathbf{X} = (\mathbf{X}_{-}, \mathbf{X}_{+})$ is structurally stable
  inside a convex compact set $K \subset C_{+} \cup C_{-}$, where $C_{+}$ and
  $C_{-}$ are two consecutive stripes meeting at $I \in
    \tilde{\mathcal{I}}_C$.
\end{definition}

In fact, given the bifurcation described below, it is natural to study the
constant and affine cases separately, since the first is always structurally
unstable inside the last one. More precisely:

\begin{proposition}
  \label{sec:stability:prop:bifurcation}
  Every $\mathbf{F} \in \mathcal{C}_3$ is structurally unstable as an element of
  $\mathcal{A}_3$.
\end{proposition}

\begin{proof}
  Let $\mathbf{F} \in \mathcal{C}_3 \subset \mathcal{D}_3^k$ be a piecewise smooth
  vector field with a double discontinuity given by constant vector fields

  \begin{equation*}
    \mathbf{F}_i(x,y,z) = (d_{i1}, d_{i2}, d_{i3}),
  \end{equation*}

  \noindent where $d_{ij} \in \mathbb{R}$ for all $i$ and $j$.

  Assume, without loss of generality, that $d_{i1} > 0$. Then, according to
  \cref{sec:constant:thm:dynamics}, over the slow manifold we have the dynamics
  $\dot{x} = d_{i1}$. As $d_{i1} > 0$, then it is strictly increasing and, in
  particular, has no singularities.

  However, considering $\mathbf{F}$ as an element of $\mathcal{A}_3 \subset
    \mathcal{D}_3^k$ and, in particular, perturbing $\mathbf{F}_i$ inside
  $\mathcal{A}_3$ with $a_{i1} \neq 0$, then we would now have the dynamics
  $\dot{x} = a_{i1}x + d_{i1}$ over the slow manifold. As $d_{i1} > 0$ and
  $a_{i1} \neq 0$, then it does now have a single singularity at $x =
    \delta_{i}$ and, besides that, half of its stability was inverted when
  compared with the unperturbed dynamics.

  In other words, $\mathbf{F}$ as an element of $\mathcal{A}_3$ violates the
  robustness of condition \cref{sec:stability:conditions:ms_hyperbolic} of
  \cref{sec:stability:prop:conditions} and, therefore, is structurally unstable.
\end{proof}

\subsection{Constant Dynamics}%
\label{sec:stability:constant}

Let $\mathbf{F} \in \mathcal{C}_3$ be a piecewise smooth vector field with a
double discontinuity given by constant vector fields

\begin{equation}
  \label{sec:stability:sub:constant:eq:system}
  \mathbf{F}_i(x,y,z) = (d_{i1}, d_{i2}, d_{i3}),
\end{equation}

\noindent with $d_{i2} \neq 0$ or $d_{i3} \neq 0$. Remember that, in this case,
\cref{sec:constant:thm:dynamics} provides a full description of the fundamental
dynamics of \cref{sec:stability:sub:constant:eq:system} and, therefore, we would
like to combine it with \cref{sec:stability:prop:conditions} to derive a
semi-local structural stability theorem.

In order to apply this results, given $\Sigma_{\theta_0} \in
  \tilde{\mathcal{I}}_C$, let $\mathbf{X} = (\mathbf{X}_{-}, \mathbf{X}_{+})$ be
the Filippov system induced by \cref{sec:stability:sub:constant:eq:system} in
a convex compact set $K \subset C_{+} \cup C_{-}$, where $C_{+}$ and $C_{-}$
are two consecutive stripes meeting at $\Sigma_{\theta_0}$ as represented at
\cref{sec:stability:fig:compact}. According to
\cref{sec:constant:thm:dynamics}, the following are the possible categories of
dynamics for a stripe $S_i \in \left\{ C_{+}, C_{-}\right\}$, which we now
analyze against conditions \cref{sec:stability:conditions:ms_hyperbolic} —
\cref{sec:stability:conditions:sg_tangency} of
\cref{sec:stability:prop:conditions} case by case in order to discover those
that can possibly generate structural stable systems, henceforth called
\textbf{candidates}:

\begin{enumerate}
  \item $d_{i2} \neq 0$:

        \begin{enumerate}
          \item $d_{i3} \neq 0$: \\
                One, and only one, of the straight lines $L_{i}$ or
                $L_{i}^{\pi}$ is visible inside the stripe. Hence, if $d_{i1} =
                  0$, then we have a continuum of singularities, i.e., a violation
                of condition \cref{sec:stability:conditions:ms_hyperbolic}.
                However, if $d_{i1} \neq 0$, then no critical elements are
                present and, therefore
                \cref{sec:stability:conditions:ms_hyperbolic} and
                \cref{sec:stability:conditions:ms_separatrix} validates. About
                \cref{sec:stability:conditions:ms_wandering}, since the slow
                manifold acts as $\alpha$ or $\omega$-limit of the surrounding
                dynamics, then it also validates if $d_{i1} \neq 0$. Even more,
                since over the borders of $S_i$ there is only transversal layer
                dynamics, then \cref{sec:stability:conditions:sg_zeros} and
                \cref{sec:stability:conditions:sg_tangency} also validates.
                Finally, observe that, invoking theorems such as
                \emph{continuity theorems} and \emph{Thom Transversality
                  Theorem}, we easily conclude the robustness of the properties
                validated above when perturbing inside $\mathcal{C}_3$.
                Therefore, this case is a \emph{candidate} if, and only if,
                $d_{i1} \neq 0$.

          \item $d_{i3} = 0$: \\
                The only difference between this case and the previous is the
                fact that, now, one of straight lines $L_{i}$ or $L_{i}^{\pi}$
                is over one of the borders of the stripe $S_i$ and, therefore,
                \cref{sec:stability:conditions:sg_tangency} is possibly
                violated, whatever $d_{i1}$. More specifically, if $L_{i}$ or
                $L_{i}^{\pi}$ coincides with $\Sigma_{\theta_0}$, then we have
                instability; otherwise, we have a \emph{candidate}.
        \end{enumerate}

  \item $d_{i2} = 0$ and $d_{i3} \neq 0$: \\
        This case is similar to the previous one ($d_{i2 \neq 0}$ and $d_{i3} =
          0$): whatever $d_{i1}$, if $L_{i}$ or $L_{i}^{\pi}$ coincides with
        $\Sigma_{\theta_0}$, then we have instability; otherwise, we have a
        \emph{candidate}.
\end{enumerate}

The analysis of the remaining conditions
\cref{sec:stability:conditions:sg_colinear} —
\cref{sec:stability:conditions:cr_recurrent} requires the combined dynamics of
the stripes $C_{+}$ and $C_{-}$. Therefore, in order to decide stability, we
shall now analyze all the combinations of candidates obtained above, and
summarized at \cref{sec:stability:sub:constant:tab:candidates}, against these
conditions.

\begin{table}[ht]
  \centering
  \caption{Conditions under which the stripe $S_i$ is a semi-local structural
    stability candidate.}
  \begin{tabular}{ccc}
    \hline
                    & $d_{i2} \neq 0$          & $d_{i2} = 0$             \\
    \hline
    $d_{i3} \neq 0$ & $d_{i1} \neq 0$          & $\theta_i \neq \theta_0$ \\
    $d_{i3} = 0$    & $\theta_i \neq \theta_0$ & unstable                 \\
    \hline
  \end{tabular}%
  \label{sec:stability:sub:constant:tab:candidates}
\end{table}

Actually, most of the remaining conditions can be easily dropped. In fact,
according to \cref{sec:constant:thm:dynamics}, none of the candidates have
periodic orbits and, besides that, because of the $\alpha$ and $\omega$-limit
nature of $\mathcal{M}_i$, an orbit that enters $S_i$ never touches the same
border again and, therefore, \cref{sec:stability:conditions:cr_periodic} always
validates, because there is no periodic orbits. Likewise, there is no
singularities, usual or not and, therefore, there is no separatrix-connections
or relations, i.e., \cref{sec:stability:conditions:cr_separatrix} always
validates. Finally, as long as $d_{i1} \neq 0$, Poincaré-Bendixson Theorem
assures that no non-trivial recurrent orbits can happen inside $S_i$ and,
besides that, again because of the $\alpha$ and $\omega$-limit nature of
$\mathcal{M}_i$, neither can they happen thought the switching manifold and,
therefore, \cref{sec:stability:conditions:cr_recurrent} also always validates.
At this point, the following theorem has been proved:

\begin{rtheorem}[Constant Dynamics Stability]
  \label{sec:stability:sub:constant:thm:conditions}
  Let $\mathbf{F} \in \mathcal{C}_3$ be given by
  \cref{sec:stability:sub:constant:eq:system} with $d_{i2} \neq 0$ or $d_{i3}
    \neq 0$. Given $\Sigma_{\theta_0} \in \tilde{\mathcal{I}}_C$, let $\mathbf{X}
    = (\mathbf{X}_{-}, \mathbf{X}_{+})$ be the Filippov system induced around
  $\Sigma_{\theta_0}$ and inside a convex compact set $K \subset C_{+} \cup
    C_{-}$, where $C_{+}$ and $C_{-}$ are two consecutive stripes meeting at
  $\Sigma_{\theta_0}$. Then, $\mathbf{F}$ is $(\Sigma_{\theta_0}, K)$-semi-local
  structurally stable in $\mathcal{C}_3$ if, and only if, $\mathbf{X}_{+}$ and
  $\mathbf{X}_{-}$ satisfies at least one of the conditions

  \begin{enumerate}
    \item $d_{i1}d_{i2}d_{i3} \neq 0$; or

    \item $d_{i1} \neq 0$, $d_{i2}^2 + d_{i3}^2 \neq 0$ and $\theta_i \neq
            \theta_0$;
  \end{enumerate}

  \noindent and, additionally, $\mathbf{X}_{+}$ and $\mathbf{X}_{-}$ are
  non-colinear over $\Sigma_{\theta_0}$, except at finitely many points.
\end{rtheorem}

\begin{example}
  \label{sec:stability:sub:constant:exmp:unstable}
  Lets see an example of instability around the discontinuity manifold
  $\Sigma_{\frac{\pi}{2}} \in \tilde{\mathcal{I}}_C$ of the cylinder generated
  by constant vector fields. More precisely, take $\mathbf{F} \in \mathcal{C}_3$
  with

  \begin{equation*}
    \mathbf{F}_1(x,y,z) = (1,-1,1)
    \quad \text{and} \quad
    \mathbf{F}_2(x,y,z) = (-1,1,1),
  \end{equation*}

  \noindent whose dynamics over the stripes $S_1 \cup S_2$, represented at
  \cref{sec:stability:sub:constant:fig:const_sep} below, can be determined as in
  \cref{sec:constant:exmp:no_singularities} using
  \cref{sec:constant:thm:dynamics}.

  Since $d_{11}d_{12}d_{13} = -1 \neq 0$ and $d_{21}d_{22}d_{23} = -1 \neq 0$,
  then the first part of \cref{sec:stability:sub:constant:thm:conditions} is
  satisfied. However, the induced dynamics $\mathbf{X}_1$ and $\mathbf{X}_2$
  over the stripes $S_1$ and $S_2$, respectively, are colinear over their whole
  intersection, the discontinuity manifold $\Sigma_{\frac{\pi}{2}}$.

  In fact, as represented at \cref{sec:stability:sub:constant:subfig:const_x1},
  for $\mathbf{F}_1$ the slow manifold consists of the straight lines given by
  $\theta = \theta_1 = -\frac{\pi}{4}$ and $\theta = \theta_1 + \pi =
    \frac{3\pi}{4}$; over then acts the increasing dynamics $\dot{x} = 1$. Besides
  that, the first line is repellor and, the second, attractor of the allround
  dynamics. On the other hand, as represented at
  \cref{sec:stability:sub:constant:subfig:const_x2}, for $\mathbf{F}_2$ the slow
  manifold consists of the straight lines given by $\theta = \theta_2 =
    \frac{\pi}{4}$ and $\theta = \theta_1 + \pi = \frac{5\pi}{4}$; over then acts
  the decreasing dynamics $\dot{x} = -1$. Besides that, the first line is
  attractor and, the second, repellor of the surrounding layer dynamics. In
  other words, the only differences between their dynamics is a
  $\pi$-translation in $\theta$ and inverse stability.

  This symmetry assures the colinearity of $\mathbf{X}_1$ and $\mathbf{X}_2$
  over $\Sigma_{\frac{\pi}{2}}$, as represented at
  \cref{sec:stability:sub:constant:subfig:const_x}. Hence, the final part of
  \cref{sec:stability:sub:constant:thm:conditions} is violated and, therefore,
  this configuration is structurally unstable around $\Sigma_{\frac{\pi}{2}}$,
  whatever the convex compact set $K$ considered. Geometrically, the instability
  here comes from the fact that each point of colinearity is associated with a
  pseudo-singularity of the sliding vector field of the Filippov system
  $\mathbf{X} = \left(\mathbf{X}_1, \mathbf{X}_2 \right)$ and, at our
  configuration we have a continuum of them. This whole continuum of
  pseudo-singularities can be easily destroyed by perturbing any of associated
  vector fields. \qed
\end{example}

\begin{figure}[ht]
  \centering
  \begin{subfigure}[b]{0.47\linewidth}
    \def\svgwidth{\linewidth}
    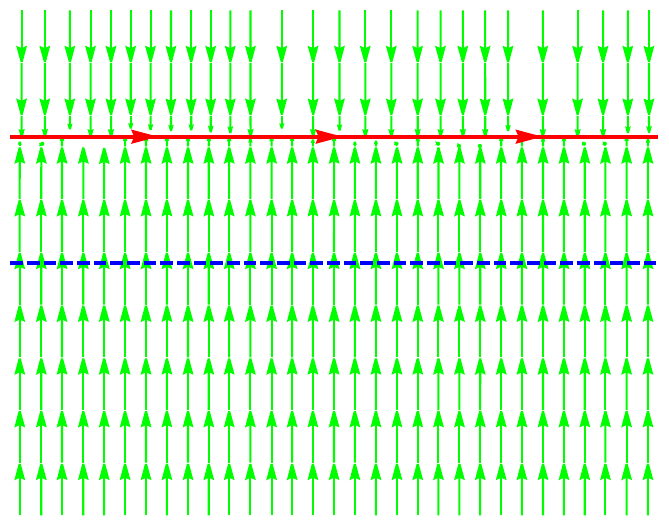
    \caption{$\mathbf{X}_1$}
    \label{sec:stability:sub:constant:subfig:const_x1}
  \end{subfigure}
  \quad
  \begin{subfigure}[b]{0.47\linewidth}
    \def\svgwidth{\linewidth}
    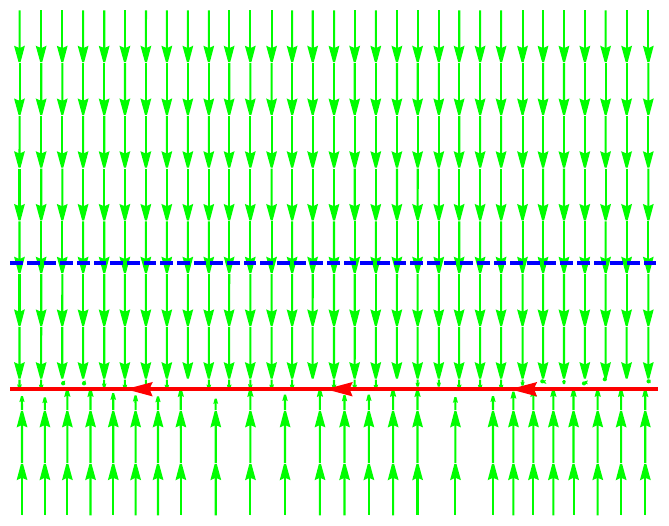
    \caption{$\mathbf{X}_2$}
    \label{sec:stability:sub:constant:subfig:const_x2}
  \end{subfigure}

  \bigskip

  \begin{subfigure}[b]{0.47\linewidth}
    \def\svgwidth{\linewidth}
    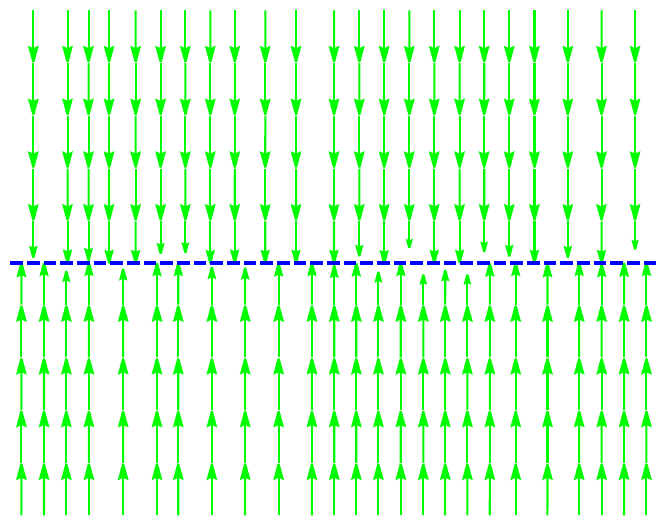
    \caption{$\mathbf{X} = \left( \mathbf{X}_1, \mathbf{X}_2 \right)$}
    \label{sec:stability:sub:constant:subfig:const_x}
  \end{subfigure}
  \quad
  \begin{subfigure}[b]{0.47\linewidth}
    \def\svgwidth{\linewidth}
    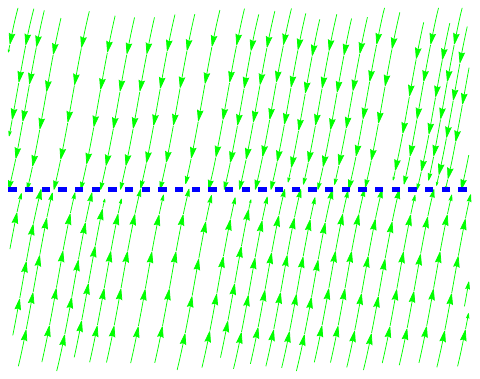
    \caption{$\mathbf{X} = \left( \mathbf{X}_1, \mathbf{X}_2 \right)$ with $r >
        0$}
    \label{sec:stability:sub:constant:subfig:const_x_r}
  \end{subfigure}
  \caption{Dynamics over the stripes $S_1 \cup S_2$ generated by the fields
    studied at \cref{sec:stability:sub:constant:exmp:unstable}.}
  \label{sec:stability:sub:constant:fig:const_sep}
\end{figure}

\subsection{Affine Dynamics}%
\label{sec:stability:affine}

Let $\mathbf{F} \in \mathcal{A}_3$ be a piecewise smooth vector field with a
double discontinuity given by affine vector fields

\begin{equation}
  \label{sec:stability:sub:affine:eq:system}
  \begin{aligned}
    \mathbf{F}_i(x,y,z) =
    ( & a_{i1}x + b_{i1}y + c_{i1}z + d_{i1},  \\
      & a_{i2}x + b_{i2}y + c_{i2}z + d_{i2},  \\
      & a_{i3}x + b_{i3}y + c_{i3}z + d_{i3}),
  \end{aligned}
\end{equation}

\noindent with $\gamma_{i} \neq 0$. Remember that, in this case,
\cref{sec:affine:thm:dynamics} provides a full description of the fundamental
dynamics of \cref{sec:stability:sub:affine:eq:system} and, therefore, as in the
previous section, we would like to combine it with
\cref{sec:stability:prop:conditions} to derive a semi-local structural stability
theorem.

In order to apply this results, given $\Sigma_{\theta_0} \in
  \tilde{\mathcal{I}}_C$, let $\mathbf{X} = (\mathbf{X}_{-}, \mathbf{X}_{+})$ be
the Filippov system induced by \cref{sec:stability:sub:affine:eq:system} in a
convex compact set $K \subset C_{+} \cup C_{-}$, where $C_{+}$ and $C_{-}$
are two consecutive stripes meeting at $\Sigma_{\theta_0}$ as represented at
\cref{sec:stability:fig:compact}. According to \cref{sec:affine:thm:dynamics},
the following are the possible categories of dynamics for a stripe $S_i \in
  \left\{ C_{+}, C_{-}\right\}$, which we now analyze against conditions
\cref{sec:stability:conditions:ms_hyperbolic} —
\cref{sec:stability:conditions:sg_tangency} of
\cref{sec:stability:prop:conditions} case by case in order to discover those
that can possibly generate structural stable systems, i.e., the
\emph{candidates}:

\begin{enumerate}
  \item $a_{i2} \neq 0$:

        \begin{enumerate}
          \item $a_{i3} \neq 0$: \\
                The characterizing property of this case is the fact that
                $\beta_{i} \neq 0$ and, therefore, the horizontal asymptotes
                resides inside the stripes, possibly even $S_i$. As a
                consequence, there is always a visible part of the slow manifold
                inside $S_i$. Hence, if $a_{i1} = 0$ and $d_{i1} = 0$, then we
                have a continuum of singularities; if $a_{i1} = 0$ and $d_{i1}
                  \neq 0$, then we have a similar bifurcation to that described at
                \cref{sec:stability:prop:bifurcation} when perturbing. Either
                way, \cref{sec:stability:conditions:ms_hyperbolic} is violated.
                However, if $a_{i1} \neq 0$, then
                \cref{sec:framework:sub:blowup:cor:singularities} assures the
                existence of at most one robust singularity $\mathbf{P}$, always
                hyperbolic and, therefore,
                \cref{sec:stability:conditions:ms_hyperbolic} and
                \cref{sec:stability:conditions:ms_separatrix} validates, since
                obviously the is no periodic orbits inside $S_i$. As in the
                constant case, the $\alpha$ or $\omega$-limit nature of the slow
                manifold also assures
                \cref{sec:stability:conditions:ms_wandering}. For
                \cref{sec:stability:conditions:sg_zeros} and
                \cref{sec:stability:conditions:sg_tangency}, observe that the
                fast dynamics is always transversal and, therefore, we only need
                the additional condition $\mathbf{P} \not\in \Sigma_{\theta_0}$.
                Finally, as in the constant case, invoking theorems such as
                \emph{continuity theorems} and \emph{Thom Transversality
                  Theorem}, we easily conclude the robustness of the properties
                validated above when perturbing inside $\mathcal{A}_3$. Therefore,
                this case is a \emph{candidate} if, and only if, $a_{i1} \neq 0$
                and $\mathbf{P} \not\in \Sigma_{\theta_0}$.

          \item $a_{i3} = 0$: \\
                The only difference between this case and the previous is the
                fact that $\beta_{i} = 0$ and, therefore, the horizontal
                asymptotes are exactly at the borders $\theta = 0$ and $\theta
                  = \pi$ of the stripes. However, since we are working inside a
                convex compact set $K$, then the same arguments of the
                previous case apply here.
        \end{enumerate}

  \item $a_{i2} = 0$: \\
        Finally, the only difference between this case and the previous ($a_{i2}
          \neq 0$ and $a_{i3} = 0$) is the fact that now the horizontal
        asymptotes are exactly at the borders $\theta = \sfrac{\pi}{2}$ and
        $\theta = \sfrac{3\pi}{2}$ of the stripes. Therefore, the same
        arguments applies.
\end{enumerate}

The analysis of the remaining conditions
\cref{sec:stability:conditions:sg_colinear} —
\cref{sec:stability:conditions:cr_recurrent} requires the combined dynamics of
the stripes $C_{+}$ and $C_{-}$. Therefore, in order to decide stability, we
need to analyze all the combinations of candidates obtained above against these
conditions. Generally, it is fairly easy to perform this analysis given a
specific combination. However, a translation of these final conditions to
parametric ones, although possible, would lead to a relatively large
number\footnote{More specifically, \cref{sec:affine:thm:dynamics} give us a
  normal form with 8 possible dynamics for each stripe. Combining them 2 by 2
  (with repetition) leave us with 36 combinations. Even if half of the
  combinations lead to a repeating condition, we would still be left with 18
  conditions!} of conditions that, worse than that, would carry little to no
geometrical meaning. Hence, leaving these final conditions ``untranslated'' is a
better approach and, therefore, the following theorem has been proved:

\TheoremAffineStability

\begin{example}
  \label{sec:stability:sub:affine:exmp:unstable}
  Lets see an example of instability around the discontinuity manifold
  $\Sigma_{0} \in \tilde{\mathcal{I}}_C$ of the cylinder generated by affine
  vector fields. More precisely, take $\mathbf{F} \in \mathcal{A}_3$ with
  $\mathbf{F}_4$ and $\mathbf{F}_1$ affine vector fields given by
  \cref{sec:stability:sub:affine:eq:system} such that

  \begin{equation*}
    \mathbf{F}_4 :
    \begin{bmatrix}
      a_{41} & d_{41} \\
      a_{42} & d_{42} \\
      a_{43} & d_{43}
    \end{bmatrix}
    =
    \begin{bmatrix*}[r]
      -1 & 1  \\
      0  & -1 \\
      1  & 0
    \end{bmatrix*}
    ~~ \text{and} ~~
    \mathbf{F}_1 :
    \begin{bmatrix}
      a_{11} & d_{11} \\
      a_{12} & d_{12} \\
      a_{13} & d_{13}
    \end{bmatrix}
    =
    \begin{bmatrix*}[r]
      1 & -1  \\
      0  & 1 \\
      1  & 0
    \end{bmatrix*},
  \end{equation*}

  \noindent whose dynamics over the stripes $S_4 \cup S_1$, represented at
  \cref{sec:stability:sub:affine:fig:unstable} below, can be determined as in
  \cref{sec:affine:exmp:slow_cycle} using \cref{sec:affine:thm:dynamics}.

  Regarding $\mathbf{F}_4$, since $a_{42} = 0$ and $\gamma_4 = -1 < 0$, then
  \cref{sec:affine:thm:dynamics} tells us that the slow manifold is a decreasing
  arctangent with horizontal asymptotes $\theta = -\sfrac{\pi}{2}$ and $\theta =
    \sfrac{\pi}{2}$, as represented at
  \cref{sub@sec:stability:sub:affine:subfig:affine_x4}. This manifold crosses
  the line $\theta = \theta_0 = 0$ at $x \in \mathbb{R}$ such that

  \begin{align*}
    0 = \theta_0 = \theta_4(x) & = \arctan{\left(\frac{a_{43}x +
    d_{43}}{d_{42}}\right)} =                                    \\
                               & = \arctan{\left( -x \right)}
    \Leftrightarrow
    x = 0,
  \end{align*}

  \noindent i.e., at the point $\mathbf{Q}_4 = (0,0)$. Besides that, over the
  slow manifold acts the dynamics $\dot{x} = -x+1$ whose only singularity at the
  point

  \begin{equation*}
    \mathbf{P}_4 = \left( \delta_4, \theta_4(\delta_4) \right) =
    \left( 1, \arctan{\left( -1 \right)} \right) =
    \left( 1, -\frac{\pi}{4} \right),
  \end{equation*}

  \noindent is stable, since $a_{41} < 0$. Even more, since $a_{42} = 0$ and
  $d_{42} < 0$ then, according to \cref{sec:affine:tab:fast}, the slow manifold
  repels the layer dynamics around. Therefore, remembering of
  \cref{sec:framework:sub:blowup:cor:singularities} we conclude that
  $\mathbf{P}_4$, as a singularity of $\mathbf{X}_4$, is a hyperbolic
  \emph{saddle}.

  \begin{figure}[ht]
    \centering
    \begin{subfigure}[b]{0.47\linewidth}
      \def\svgwidth{\linewidth}
\begingroup%
  \makeatletter%
  \providecommand\color[2][]{%
    \errmessage{(Inkscape) Color is used for the text in Inkscape, but the package 'color.sty' is not loaded}%
    \renewcommand\color[2][]{}%
  }%
  \providecommand\transparent[1]{%
    \errmessage{(Inkscape) Transparency is used (non-zero) for the text in Inkscape, but the package 'transparent.sty' is not loaded}%
    \renewcommand\transparent[1]{}%
  }%
  \providecommand\rotatebox[2]{#2}%
  \newcommand*\fsize{\dimexpr\f@size pt\relax}%
  \newcommand*\lineheight[1]{\fontsize{\fsize}{#1\fsize}\selectfont}%
  \ifx\svgwidth\undefined%
    \setlength{\unitlength}{315.91592407bp}%
    \ifx\svgscale\undefined%
      \relax%
    \else%
      \setlength{\unitlength}{\unitlength * \real{\svgscale}}%
    \fi%
  \else%
    \setlength{\unitlength}{\svgwidth}%
  \fi%
  \global\let\svgwidth\undefined%
  \global\let\svgscale\undefined%
  \makeatother%
  \begin{picture}(1,0.79134979)%
    \lineheight{1}%
    \setlength\tabcolsep{0pt}%
    \put(0,0){\includegraphics[width=\unitlength,page=1]{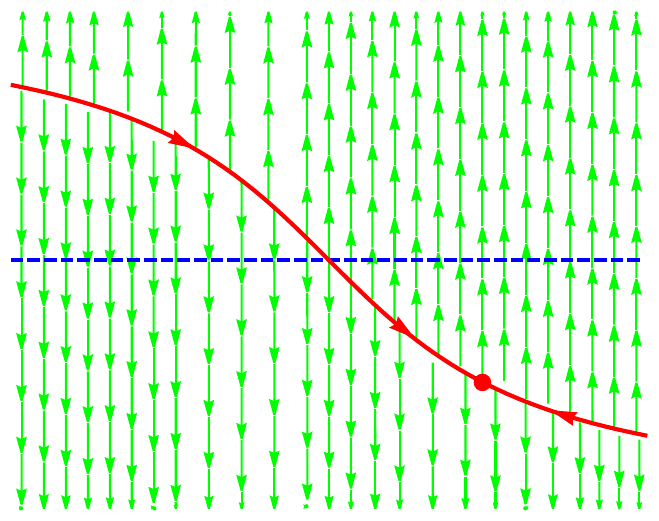}}%
    \put(0.74824717,0.24058241){\makebox(0,0)[lt]{\lineheight{1.25}\smash{\begin{tabular}[t]{l}$\mathbf{P}_4$\end{tabular}}}}%
    \put(0.49659789,0.42303534){\makebox(0,0)[lt]{\lineheight{1.25}\smash{\begin{tabular}[t]{l}$\mathbf{Q}_4$\end{tabular}}}}%
  \end{picture}%
\endgroup%

      \caption{$\mathbf{X}_4$}
      \label{sec:stability:sub:affine:subfig:affine_x4}
    \end{subfigure}
    \quad
    \begin{subfigure}[b]{0.47\linewidth}
      \def\svgwidth{\linewidth}
\begingroup%
  \makeatletter%
  \providecommand\color[2][]{%
    \errmessage{(Inkscape) Color is used for the text in Inkscape, but the package 'color.sty' is not loaded}%
    \renewcommand\color[2][]{}%
  }%
  \providecommand\transparent[1]{%
    \errmessage{(Inkscape) Transparency is used (non-zero) for the text in Inkscape, but the package 'transparent.sty' is not loaded}%
    \renewcommand\transparent[1]{}%
  }%
  \providecommand\rotatebox[2]{#2}%
  \newcommand*\fsize{\dimexpr\f@size pt\relax}%
  \newcommand*\lineheight[1]{\fontsize{\fsize}{#1\fsize}\selectfont}%
  \ifx\svgwidth\undefined%
    \setlength{\unitlength}{315.9156189bp}%
    \ifx\svgscale\undefined%
      \relax%
    \else%
      \setlength{\unitlength}{\unitlength * \real{\svgscale}}%
    \fi%
  \else%
    \setlength{\unitlength}{\svgwidth}%
  \fi%
  \global\let\svgwidth\undefined%
  \global\let\svgscale\undefined%
  \makeatother%
  \begin{picture}(1,0.78632274)%
    \lineheight{1}%
    \setlength\tabcolsep{0pt}%
    \put(0,0){\includegraphics[width=\unitlength,page=1]{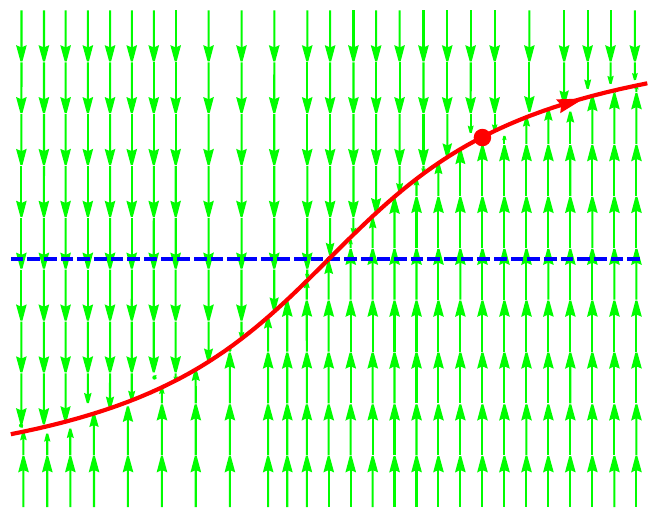}}%
    \put(0.74824784,0.52265299){\makebox(0,0)[lt]{\lineheight{1.25}\smash{\begin{tabular}[t]{l}$\mathbf{P}_1$\end{tabular}}}}%
    \put(0,0){\includegraphics[width=\unitlength,page=2]{affine_x1.pdf}}%
    \put(0.49659837,0.33872808){\makebox(0,0)[lt]{\lineheight{1.25}\smash{\begin{tabular}[t]{l}$\mathbf{Q}_1$\end{tabular}}}}%
    \put(0,0){\includegraphics[width=\unitlength,page=3]{affine_x1.pdf}}%
  \end{picture}%
\endgroup%

      \caption{$\mathbf{X}_1$}
      \label{sec:stability:sub:affine:subfig:affine_x1}
    \end{subfigure}

    \bigskip

    \begin{subfigure}[b]{0.47\linewidth}
      \def\svgwidth{\linewidth}
\begingroup%
  \makeatletter%
  \providecommand\color[2][]{%
    \errmessage{(Inkscape) Color is used for the text in Inkscape, but the package 'color.sty' is not loaded}%
    \renewcommand\color[2][]{}%
  }%
  \providecommand\transparent[1]{%
    \errmessage{(Inkscape) Transparency is used (non-zero) for the text in Inkscape, but the package 'transparent.sty' is not loaded}%
    \renewcommand\transparent[1]{}%
  }%
  \providecommand\rotatebox[2]{#2}%
  \newcommand*\fsize{\dimexpr\f@size pt\relax}%
  \newcommand*\lineheight[1]{\fontsize{\fsize}{#1\fsize}\selectfont}%
  \ifx\svgwidth\undefined%
    \setlength{\unitlength}{315.73532104bp}%
    \ifx\svgscale\undefined%
      \relax%
    \else%
      \setlength{\unitlength}{\unitlength * \real{\svgscale}}%
    \fi%
  \else%
    \setlength{\unitlength}{\svgwidth}%
  \fi%
  \global\let\svgwidth\undefined%
  \global\let\svgscale\undefined%
  \makeatother%
  \begin{picture}(1,0.78944168)%
    \lineheight{1}%
    \setlength\tabcolsep{0pt}%
    \put(0,0){\includegraphics[width=\unitlength,page=1]{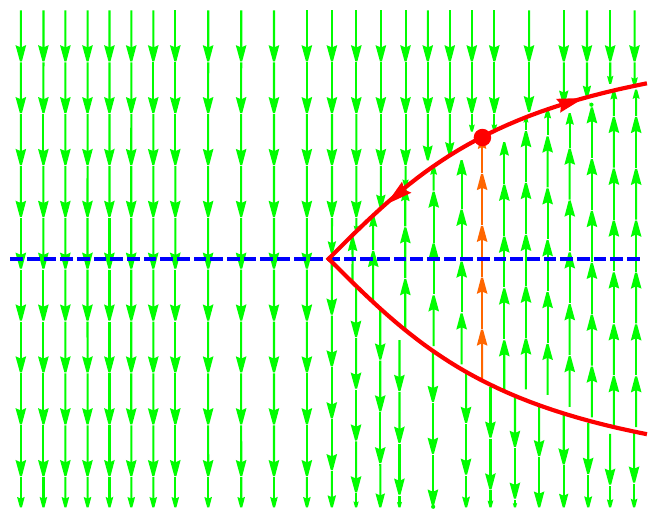}}%
    \put(0.74867515,0.52562137){\makebox(0,0)[lt]{\lineheight{1.25}\smash{\begin{tabular}[t]{l}$\mathbf{P}_1$\end{tabular}}}}%
    \put(0,0){\includegraphics[width=\unitlength,page=2]{affine_x.pdf}}%
    \put(0.74867515,0.22885698){\makebox(0,0)[lt]{\lineheight{1.25}\smash{\begin{tabular}[t]{l}$\mathbf{P}_4$\end{tabular}}}}%
    \put(0,0){\includegraphics[width=\unitlength,page=3]{affine_x.pdf}}%
    \put(0.44440745,0.42091579){\makebox(0,0)[lt]{\lineheight{1.25}\smash{\begin{tabular}[t]{l}$\mathbf{Q}$\end{tabular}}}}%
  \end{picture}%
\endgroup%

      \caption{$\mathbf{X} = \left( \mathbf{X}_4, \mathbf{X}_1 \right)$}
      \label{sec:stability:sub:affine:subfig:affine_x}
    \end{subfigure}
    \quad
    \begin{subfigure}[b]{0.47\linewidth}
      \def\svgwidth{\linewidth}
      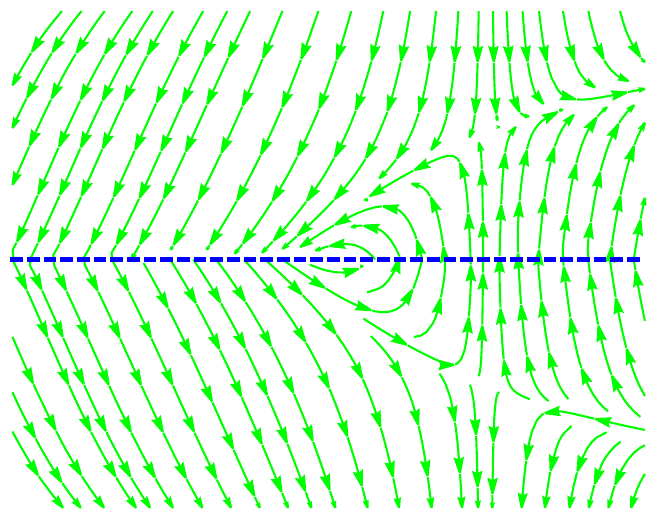
      \caption{$\mathbf{X} = \left( \mathbf{X}_4, \mathbf{X}_1 \right)$ with $r
          > 0$}
      \label{sec:stability:sub:affine:subfig:affine_rot}
    \end{subfigure}
    \caption{Dynamics over the stripes $S_4 \cup S_1$ generated by the fields
      studied at \cref{sec:stability:sub:affine:exmp:unstable}.}
    \label{sec:stability:sub:affine:fig:unstable}
  \end{figure}

  On the other hand, regarding $\mathbf{F}_1$, since $a_{12} = 0$ and $\gamma_1
    = 1 > 0$, then \cref{sec:affine:thm:dynamics} tells us that the slow
  manifold is a decreasing arctangent with horizontal asymptotes $\theta =
    -\sfrac{\pi}{2}$ and $\theta = \sfrac{\pi}{2}$, as represented at
  \cref{sec:stability:sub:affine:subfig:affine_x1}. This manifold crosses the
  line $\theta = \theta_0 = 0$ at $x \in \mathbb{R}$ such that

  \begin{align*}
    0 = \theta_0 = \theta_1(x) & = \arctan{\left(\frac{a_{13}x +
    d_{13}}{d_{12}}\right)}                                      \\
                               & = \arctan{\left( x \right)}
    \Leftrightarrow
    x = 0,
  \end{align*}

  \noindent i.e., also at the point $\mathbf{Q}_1 = (0,0)$. Besides that, over
  the slow manifold acts the dynamics $\dot{x} = x-1$ whose only singularity at
  the point

  \begin{equation*}
    \mathbf{P}_1 = \left( \delta_1, \theta_1(\delta_1) \right) =
    \left( 1, \arctan{\left( 1 \right)} \right) =
    \left( 1, \frac{\pi}{4} \right),
  \end{equation*}

  \noindent is unstable, since $a_{11} > 0$. Even more, since $a_{12} = 0$ and
  $d_{12} > 0$ then, according to \cref{sec:affine:tab:fast}, the slow manifold
  attracts the layer dynamics around. Therefore, remembering
  \cref{sec:framework:sub:blowup:cor:singularities} we conclude that
  $\mathbf{P}_1$, as a singularity of $\mathbf{X}_1$, is also a hyperbolic
  \emph{saddle}.

  Hence, as represented at \cref{sub@sec:stability:sub:affine:subfig:affine_x},
  since $\mathbf{Q}_4 = \mathbf{Q}_1$ with $\mathbf{P}_4$ and $\mathbf{P}_1$
  hyperbolic saddles, then the Filippov system $\mathbf{X} = (\mathbf{X}_4,
    \mathbf{X}_1)$ has a separatrix-connection and, therefore, it violates
  condition \cref{sec:stability:conditions:cr_separatrix} of
  \cref{sec:stability:prop:conditions}, whatever the convex compact set $K$
  considered. In other words, according to
  \cref{sec:stability:sub:affine:thm:conditions}, this configuration is
  structurally unstable around the discontinuity manifold $\Sigma_0$.

  Finally, we observe that, as represented at
  \cref{sec:stability:sub:affine:subfig:affine_x}, there is actually two
  separatrix-connections between the saddles $\mathbf{P}_4$ and $\mathbf{P}_1$.
  These connections enclose a \emph{rotating region}, represented at
  \cref{sec:stability:sub:affine:subfig:affine_rot}. \qed
\end{example}

\section{Conclusion}%
\label{sec:conclusion}

In this work, we tackled the problems stated at \cref{sec:problem}: given
$\mathbf{F}$ with a double discontinuity (two switches) as switching manifold,
can we define Filippov-like dynamics over the singular part, i.e., the
intersection $\Sigma_{x}$? How do these dynamics generally behave there?
Specifically, we leveraged \emph{Buzzi's} blow-up methodology to approach these
problems, resulting initially in the Fundamental
\cref{sec:framework:blowup:thm:dynamics}, which induces the so-called
\emph{fundamental dynamics}: a (regular) discontinuous slow-fast dynamics
happening over a cylinder representing $\Sigma_{x}$ after the blow-up. Many
qualitative properties were derived for the general non-linear case, as long as
one of the so-called \emph{fundamental hypotheses} was satisfied.

Especially, when $\mathbf{F}$ is composed of constant or affine vector fields, we
were able to fully describe the qualitative aspects of the fundamental dynamics
as presented in \cref{sec:constant:thm:dynamics} and
\cref{sec:affine:thm:dynamics}, respectively. Once we had this knowledge,
\cite{Broucke2001, Gomide2020} inspired us to look after the semi-local
structural stability of the fundamental dynamics, resulting in
\cref{sec:stability:sub:constant:thm:conditions} and
\cref{sec:stability:sub:affine:thm:conditions}, which characterize stability
through a set of simple algebraical and geometrical conditions.

To the best of our knowledge, none of the other methodologies (\emph{Jeffrey's}
and \emph{Dieci's}) discussed in \cref{sec:introduction} were able to deliver
similar results, expressing then the effectiveness and practicality of
\emph{Buzzi's} blow-up based methodology as presented and improved here.
However, we acknowledge the beauty of \emph{Jeffrey's} canopy based methodology,
which directly extends the Filippov dynamics to the singular part of the
switching manifold through convex combinations. In fact, we conjecture and look
forward to prove the equivalence of these methodologies, then unifying its
strengths.

Howsoever, on the applicability of these methodologies, we expand the remark in
section 8.1 of \cite[p.~1102]{Jeffrey2014} by Mike Jeffrey as follows:

\begin{quote}
    \emph{``If one is able to find physical laws for the dynamics on
        $\mathcal{D}$, then these supersede the discontinuous model (\ldots),
        and whether these agree with the \sout{ canopy and dummy dynamics}
        \underline{model} is open for experimenters of various disciplines to
        put to the test.''}
\end{quote}

\noindent Nature has the final word.

\section*{Acknowledgments}
\addcontentsline{toc}{section}{Acknowledgments}

The first author was financed in part by the Coordenação de Aperfeiçoamento de
Pessoal de Nível Superior — Brasil (CAPES) — Finance Code 001. The second author
was partially supported by FAPESP grants 2018/03338--6 and 2018/13481--0.

\phantomsection%
\addcontentsline{toc}{section}{\refname}
\bibliographystyle{acm}
\bibliography{references}

\end{document}